\documentclass[a4paper]{amsart}
\usepackage[francais]{babel}
\usepackage{amsmath,amssymb,amsthm}
\usepackage{amsfonts}
\usepackage[latin1]{inputenc}
\usepackage[T1]{fontenc}
\usepackage{cite}
\usepackage{pifont}
\usepackage{euscript}
\usepackage{footnpag}
\RequirePackage{mathrsfs}

\def\noignorespaces{{\catcode`\ =11\ }}
\newcounter{constante}
\newcommand{\cst}{\refstepcounter{constante}\theconstante}

\def\cD{\mathcal{D}}

\DeclareMathOperator{\spec}{Spec}
\DeclareMathOperator{\card}{card} 
\DeclareMathOperator{\codim}{codim}

\DeclareMathOperator{\Proj}{Proj}

\DeclareMathOperator{\lie}{Lie}
\DeclareMathOperator{\specf}{Specf}

\theoremstyle{plain}
\newtheorem{theo}{Th\'eor\`eme}[section]
\newtheorem{fait}[theo]{Fait}
\newtheorem{lemm}[theo]{Lemme}
\newtheorem{prop}[theo]{Proposition}
\newtheorem{coro}[theo]{Corollaire}

\theoremstyle{definition}
\newtheorem{defi}[theo]{D\'efinition}
\newtheorem{rema}[theo]{Remarque}
\newtheorem{remas}[theo]{Remarques}

\title[{\'E}tude du cas rationnel]{{\'E}tude du cas rationnel de la th{\'e}orie des formes lin{\'e}aires de logarithmes}
\author{{\'E}ric Gaudron}
\date{8 septembre 2007}
\begin{document}

\footnotetext{MSC~$2000$: 11J86 (11J61, 11J13)}
\footnotetext{\textbf{Mots clefs}: formes lin{\'e}aires de logarithmes, cas rationnel, m{\'e}thode de Baker, groupe alg\'ebrique commutatif, taille de sous-sch{\'e}ma formel, lemme d'interpolation $p$-adique, lemme de Siegel absolu.}

\begin{abstract}
Dans ce travail, nous {\'e}tablissons des mesures d'ind{\'e}pendance lin{\'e}aire de logarithmes d'un groupe alg{\'e}brique commutatif dans le cas rationnel. Plus pr{\'e}cis{\'e}ment, soit $k$ un corps de nombres et $v_{0}$ une place quelconque de $k$. Soit $G$ un groupe alg{\'e}brique commutatif d{\'e}fini sur $k$ et $H$ un sous-groupe alg{\'e}brique connexe de $G$, d'alg{\`e}bre de Lie $\lie(H)$. Soit $u\in\lie(G(\mathbf{C}_{v_{0}}))$ un logarithme d'un point $p$ de $G(k)$. Dans le \emph{cas non-p{\'e}riodique} (le point $p$ n'est pas de torsion modulo certains sous-groupes de $G$), nous obtenons des minorations de la distance de $u$ {\`a} $\lie(H)\otimes_{k}\mathbf{C}_{v_{0}}$ qui g{\'e}n{\'e}ralisent en partie les mesures d{\'e}j{\`a} connues dans le cas d'un groupe lin{\'e}aire. Les principales caract{\'e}risques de ces r{\'e}sultats sont d'une part d'am{\'e}liorer la d{\'e}pendance en la hauteur $\log a$ du point $p$, en supprimant une puissance de $\log\log a$, et, d'autre part, d'{\^e}tre valides dans un contexte tr{\`e}s g{\'e}n{\'e}ral. La d{\'e}monstration utilise le formalisme des tailles de sous-sch{\'e}mas formels au sens de Bost en association avec un lemme arithm{\'e}tique de Raynaud. Nous avons {\'e}galement recours {\`a} un lemme de Siegel absolu et, lorsque $v_{0}$ est ultram{\'e}trique, {\`a} un lemme d'interpolation de Roy. \\[.2cm] \textsc{Abstract}. We establish new measures of linear independence of logarithms on commutative algebraic groups in the so-called \emph{rational case}. More precisely, let $k$ be a number field and $v_{0}$ be an arbitrary place of $k$. Let $G$ be a commutative algebraic group defined over $k$ and $H$ be a connected algebraic subgroup of $G$. Denote by $\lie(H)$ its Lie algebra at the origin. Let $u\in\lie(G(\mathbf{C}_{v_{0}}))$ a logarithm of a point $p\in G(k)$. Assuming (essentially) that $p$ is not a torsion point modulo proper connected algebraic subgroups of $G$, we obtain lower bounds for the distance from $u$ to $\lie(H)\otimes_{k}\mathbf{C}_{v_{0}}$. For the most part, they generalize the measures already known when $G$ is a linear group. The main feature of these results is to provide a better dependence in the height $\log a$ of $p$, removing a polynomial term in $\log\log a$. The proof relies on sharp estimates of sizes of formal subschemes associated to $H$ (in the sense of Bost) obtained from a lemma by Raynaud as well as an absolute Siegel lemma and, in the ultrametric case, a recent interpolation lemma by Roy.   
\end{abstract}
\maketitle
\newpage
\section*{Notations et conventions}
{\footnotesize
Soit $g$ un entier naturel $\ge 1$. Pour $\mathbf{t}=(t_{1},\ldots,t_{g})\in\mathbf{N}^{g}$, on d{\'e}signe par $\vert\mathbf{t}\vert$ la longueur $t_{1}+\cdots+t_{g}$ de $\mathbf{t}$, et, si $\mathbf{X}=(X_{1},\ldots,X_{g})$ est un $g$-uplet de variables (ou d'objets math{\'e}matiques plus g{\'e}n{\'e}raux, tels des op{\'e}rateurs diff{\'e}rentiels), on note $\mathbf{X}^{\mathbf{t}}=X_{1}^{t_{1}}\cdots X_{g}^{t_{g}}$. Si $x$ un nombre r{\'e}el, on note $\log^{+}(x)=\log\max{\{1,x\}}$ et $[x]$ la partie enti{\`e}re de $x$.\par  Si $G$ est un sch{\'e}ma en groupes sur un corps commutatif, $t_{G}$ ou $\lie(G)$ d{\'e}signe son espace tangent {\`a} l'origine. Si $E$ est un espace vectoriel, $\mathbf{S}(E)$ (\emph{resp}. $S^{g}(E)$) est l'alg\`ebre sym\'etrique de $E$ (\emph{resp}. la composante de degr\'e $g$ de $\mathbf{S}(E)$) et $\mathbf{P}(E)$ d\'esigne le sch\'ema projectif $\Proj\mathbf{S}(E)$.\par
 Lorsque $k$ est un corps commutatif, on note $\overline{k}$ une cl\^oture alg{\'e}brique de $k$. Soit dor\'enavant $k$ un corps de nombres, d'anneau des entiers $\mathcal{O}_{k}$, et $v$ une place de $k$. \subsubsection*{Normes et valeurs absolues}\begin{enumerate}\item[$\bullet$] Soit $v$ une place ultram\'etrique de $k$, qui correspond \`a un id\'eal premier $\mathfrak{p}$ de $\mathcal{O}_{k}$, et $p$ le nombre premier qui engendre l'id\'eal $\mathfrak{p}\cap\mathbf{Z}$. On note $k_{v}$ (\emph{resp}. $\mathcal{O}_{v}$) le compl\'et\'e $\mathfrak{p}$-adique de $k$ en $\mathfrak{p}$ (\emph{resp}. son anneau de valuation). On munit $k_{v}$ de l'unique valeur absolue $\vert\cdotp\vert_{v}$ qui v\'erifie $\vert p\vert_{v}=p^{-1}$. Cette valeur absolue s'\'etend (de mani\`ere unique) \`a $\overline{k}_{v}$ et, en particulier, aux extensions finies de $k_{v}$. Soit $\mathbf{C}_{v}=\mathbf{C}_{p}$ le compl{\'e}t{\'e} du corps valu{\'e} $(\overline{k}_{v},\vert.\vert_{v})$. Si $E_{v}$ est un $\mathbf{C}_{v}$-espace vectoriel de dimension $g$, alors toute base $(e_{1},\ldots,e_{g})$ de $E_{v}$ d\'efinit une norme $v$-adique $\Vert\cdotp\Vert_{v}$ sur $E_{v}$ par \begin{equation*}\left\Vert\sum_{i=1}^{g}{x_{i}e_{i}}\right\Vert_{v}:=\max_{1\le i\le g}{\{\vert x_{i}\vert_{v}\}}\ .\end{equation*}Ainsi, lorsque $E_{v}=\mathbf{C}_{v}^{g}$ est muni de sa base canonique, on note $\Vert\EuScript{F}\Vert_{v}$ ou $\vert\EuScript{F}\vert_{v}$ la norme d'un vecteur $\EuScript{F}=(f_{1},\ldots,f_{g})\in\mathbf{C}_{v}^{g}$~:\begin{equation*}\Vert\EuScript{F}\Vert_{v}=\vert\EuScript{F}\vert_{v}=\max_{1\le i\le g}{\{\vert f_{i}\vert_{v}\}}\ .\end{equation*}
\item[$\bullet$] Soit $v$ une place archim\'edienne de $k$. On munit $\mathbf{C}_{v}=\mathbf{C}$ de la valeur absolue usuelle. Si $\EuScript{F}=(f_{1},\ldots,f_{g})\in\mathbf{C}^{g}$, on note \begin{equation*}\vert\EuScript{F}\vert_{v}:=\max_{1\le i\le g}{\{\vert f_{i}\vert_{v}\}}\quad\text{et}\quad\Vert\EuScript{F}\Vert_{v}:=\left(\sum_{i=1}^{g}{\vert f_{i}\vert^{2}}\right)^{1/2}\ \cdotp\end{equation*}
\end{enumerate}
Avec ces conventions, la formule du produit s'\'ecrit, $\forall\,\alpha\in k\setminus\{0\}$, $\prod_{v}{\vert\alpha\vert_{v}^{[k_{v}:\mathbf{Q}_{v}]}}=1$ o\`u $v$ parcourt l'ensemble des places de $k$ et $[k_{v}:\mathbf{Q}_{v}]$ est le degr\'e local $1,2$ ou $[k_{v}:\mathbf{Q}_{p}]$ selon le caract\`ere r\'eel, complexe ou $p$-adique de la place $v$.
\subsubsection*{Hauteurs}Soit $\EuScript{F}\in k^{g}\setminus\{0\}$. La hauteur de Weil (logarithmique absolue\footnote{Comme le seront toutes les hauteurs de ce texte.}) de $\EuScript{F}$ est \begin{equation*}h(\EuScript{F})=\sum_{v}{\frac{[k_{v}:\mathbf{Q}_{v}]}{[k:\mathbf{Q}]}\log\vert\EuScript{F}\vert_{v}}\ \cdotp\end{equation*}C'est une hauteur projective ($\forall\,\alpha\in k\setminus\{0\}$, $h(\alpha\EuScript{F})=h(\EuScript{F})$) et elle se prolonge naturellement aux points de $\mathbf{P}^{g-1}(k)$. La hauteur $\mathrm{L}^{2}$ de $\EuScript{F}$ est \begin{equation*}h_{\mathrm{L}^{2}}(\EuScript{F})=\sum_{v}{\frac{[k_{v}:\mathbf{Q}_{v}]}{[k:\mathbf{Q}]}\log\Vert\EuScript{F}\Vert_{v}}\ \cdotp\end{equation*}On a $h(\EuScript{F})\le h_{\mathrm{L}^{2}}(\EuScript{F})\le h(\EuScript{F})+\frac{1}{2}\log(\#\EuScript{F})$ et l'in\'egalit\'e de Liouville $\log\vert\alpha\vert_{v}\ge-[\mathbf{Q}(\alpha):\mathbf{Q}]h(\{1,\alpha\})$ (pour toute place $v$ de $\mathbf{Q}(\alpha)$). Soit $v_{1},\ldots,v_{d}$ des vecteurs lin\'eairement ind\'ependants de $\overline{\mathbf{Q}}^{g}$ (muni de sa base canonique $e_{1},\ldots,e_{g}$). La hauteur de Schmidt de $(v_{1},\ldots,v_{d})$ est la hauteur $\mathrm{L}^{2}$ de l'ensemble des coordonn\'ees de Pl{\"u}cker du produit ext\'erieur $v_{1}\wedge\cdots\wedge v_{d}$ dans la base $e_{i_{1}}\wedge\cdots\wedge e_{i_{d}}$ ($1\le i_{1}<\cdots<i_{d}\le g$). Cette d\'efinition ne d\'epend en r{\'e}alit{\'e} que de l'espace vectoriel $V$ sur $\overline{\mathbf{Q}}$ engendr\'e par les vecteurs $v_{1},\ldots,v_{d}$ et induit la hauteur de Schmidt $h(V)$ de $V$ (et, par convention, $h(\{0\})=0$).}
\section{Introduction}
L'objectif de ce travail est d'{\'e}tablir des r{\'e}sultats g{\'e}n{\'e}raux --- archim{\'e}diens et ultram{\'e}triques --- de la \emph{th{\'e}orie des formes lin{\'e}aires de logarithmes} dans le cas particulier o\`u le lieu des z{\'e}ros des formes lin{\'e}aires est une alg{\`e}bre de Lie \emph{alg{\'e}brique}, c'est-\`a-dire l'alg\`ebre de Lie d'un sous-groupe alg\'ebrique.\par Soit $G$ un groupe alg{\'e}brique commutatif de dimension $g$ d{\'e}fini sur un corps de nombres $k$ et $v_{0}$ une place (quelconque) de $k$. Le groupe de Lie $G(\mathbf{C}_{v_{0}})$ poss{\`e}de une application exponentielle $\exp_{v_{0}}:\mathscr{T}_{v_{0}}\to G(\mathbf{C}_{v_{0}})$ d{\'e}finie sur un voisinage ouvert $\mathscr{T}_{v_{0}}$ de $0$ dans l'alg{\`e}bre de Lie de $G(\mathbf{C}_{v_{0}})$ (not{\'e}e $\lie G(\mathbf{C}_{v_{0}})$ ou $t_{G}(\mathbf{C}_{v_{0}})$ dans la suite). Consid{\'e}rons un {\'e}l{\'e}ment $u\ne 0$ de $\mathscr{T}_{v_{0}}$ d'exponentielle $p:=\exp_{v_{0}}(u)$ $k$-\emph{rationnelle} et donnons-nous par ailleurs une norme $\Vert\cdotp\Vert_{v_{0}}$ sur $\mathscr{T}_{v_{0}}$ (de distance associ{\'e}e $\mathrm{d}_{v_{0}}$) ainsi qu'une fonction hauteur $h\ge 0$ sur $G(\overline{\mathbf{Q}})$ (provenant par exemple d'une hauteur de Weil sur un espace projectif dans lequel $G$ se plonge).\par La th{\'e}orie des formes lin{\'e}aires de logarithmes consiste, dans son aspect quantitatif, \`a fournir des minorations de la distance $\mathrm{d}_{v_{0}}(u,V)$ entre $u$ et une sous-$k$-alg{\`e}bre de Lie $V$ de $\lie(G)$, en fonction des invariants li{\'e}s aux donn{\'e}es introduites (hauteur des quantit{\'e}s alg{\'e}briques, norme de $u$, degr{\'e} de $k$, etc.). Pour qualifier une minoration de $\mathrm{d}_{v_{0}}(u,V)$, on parle aussi de \emph{mesure d'ind{\'e}pendance lin{\'e}aire de logarithmes}. Le donn\'ee principale \`a laquelle on va s'int\'eresser ici est la hauteur du point $p$. \par Le cas $V=\{0\}$ est d\'ej\`a remarquable. Une cons\'equence des r\'esultats pr\'esent\'es au \S~\ref{subsec:resultats} est l'existence d'une fonction $c_{\cst}=c_{\theconstante}(G,k,v_{0},\Vert.\Vert_{v_{0}})\ge 1$, ind{\'e}pendante de $p$, telle que \begin{equation}\label{introaseptequn}\log\Vert u\Vert_{v_{0}}\ge-c_{\theconstante}\max{\{1,h(p)\}}\end{equation}pourvu que le sous-grou\-pe engendr{\'e} par $p$ ne rencontre aucun sous-grou\-pe alg\'ebrique strict de $G(\overline{k})$ sauf en $0$. Il s'agit d'une variante sophistiqu\'ee de l'in\'egalit\'e de Liouville (voir la discussion \`a la suite du corollaire~$1.2$ de~\cite{artepredeux}). Une propri\'et\'e importante de cette minoration est d'{\^e}tre optimale en la hauteur de $p$, comme on peut le voir en se pla{\c c}ant sur le groupe additif $\mathbb{G}_{\mathrm{a}}$. Pour un espace $V$ quelconque, les meilleures mesures connues de $\log\mathrm{d}_{v_{0}}(u,V)$ en fonction de $h(p)$ sont de la forme~: \begin{equation}\label{avril}\log\mathrm{d}_{v_{0}}(u,V)\ge-c_{\cst}\max{\{1,h(p)\}}^{g/t+\epsilon}\end{equation}o\`u $t$ est la codimension de $V$ dans $\lie(G)$, $\epsilon>0$ un nombre r\'eel et $c_{\theconstante}$ une fonction qui ne d\'epend pas de la hauteur de $p$ ({}\footnote{La quantit\'e $h(p)^{\epsilon}$ dans~\eqref{avril} peut {\^e}tre remplac\'ee par une puissance enti\`ere convenable du logarithme de $h(p)$.}). \par Dans cet article, nous montrons que l'on peut supprimer $\epsilon$ dans le minorant~\eqref{avril} dans le cas dit \emph{rationnel} o\`u l'alg\`ebre de Lie $V$ est \emph{alg\'ebrique}. Nous d\'emontrerons le r{\'e}sultat suivant (qui sera rendu plus pr{\'e}cis au \S~\ref{subsec:resultats}).
\begin{theo}\label{enonce2demoi}
Soit $G,k,v_{0},\Vert.\Vert_{v_{0}},p,u,V,t$ les donn{\'e}es g{\'e}n{\'e}rales introduites ci-dessus. Il existe une fonction $c_{\cst}=c_{\theconstante}(G,k,v_{0},\Vert u\Vert_{v_{0}})$ ayant la propri{\'e}t{\'e} suivante. Supposons d'une part que le sous-groupe de $G$ engendr{\'e} par $p$ ne rencontre aucun sous-groupe alg{\'e}brique strict de $G$ (sauf en $0$) et, d'autre part, que $V$ est \emph{l'alg{\`e}bre de Lie d'un sous-groupe alg{\'e}brique} $H$ de $G$. Soit $b$ un nombre r{\'e}el $\ge e$ tel que $\log b$ soit un majorant de la hauteur de $V$. Alors \begin{equation}\label{introaseptequndeuxtroisquatrecinq}\log\mathrm{d}_{v_{0}}(u,V)\ge-c_{\theconstante}(\log b)^{1+\frac{g+1}{t}}\max{\{1,h(p)\}}^{g/t}\ .\end{equation} 
\end{theo}
Jusqu'\`a pr\'esent, seul le cas d'une puissance du groupe multiplicatif $\mathbb{G}_{\mathrm{m}}$ et d'\emph{une} forme lin\'eaire ($t=1$) a vraiment \'et\'e \'etudi\'e (voir~\cite{BaWu93,dong,matveev,durhampph-miw,pph-miw2,miw4,durhamwustholz,Yu2007}\footnote{Vu la grande richesse de la litt{\'e}rature sur ce th{\`e}me, il serait vain d'essayer d'entrer dans tous les d{\'e}tails sans augmenter de mani{\`e}re exponentielle cette introduction\,; le lecteur int{\'e}ress{\'e} pourra se reporter au livre de M.Waldschmidt~\cite{miw4} dont, en particulier, le \S~$10.4$ \og The state of the Art\fg\ ainsi que les pages $545$ {\`a} $547$ qui retracent les principales {\'e}tapes de l'histoire du sujet.}). Bien qu'il soit extr{\^e}mement probable que la minoration~\eqref{introaseptequndeuxtroisquatrecinq} reste vraie en supprimant le terme $(g+1)/t$ dans l'exposant de $\log b$, il s'av{\`e}re que les m{\'e}thodes employ{\'e}es ici pour d{\'e}montrer cette in{\'e}galit{\'e} ne sont pas suffisantes pour obtenir cela.
\tableofcontents
\subsection{Donn\'ees g\'en\'erales}\label{donnesgeneralesasept}
Dans ce paragraphe, nous fixons des notations qui seront utilis{\'e}es tout au long de ce texte. Certains des th{\'e}or{\`e}mes qui vont suivre ne seront valides qu'avec des hypoth{\`e}ses suppl{\'e}mentaires sur les objets introduits ici, hypoth{\`e}ses qui seront alors explicitement mentionn{\'e}es.\par Soit $n$ un entier naturel $\ge 1$, $k$ un corps de nombres de degr{\'e} $D$ et $v_{0}$ une place quelconque (archim{\'e}dienne ou ultram{\'e}trique) de $k$ qui sera privil{\'e}gi\'ee par la suite.\par Soit $G_{1},\ldots,G_{n}$ des groupes alg{\'e}briques (connexes) commutatifs d{\'e}finis sur $k$ et $g_{1},\ldots,g_{n}$ leurs dimensions respectives. Soit $\Phi_{i}:G_{i}\hookrightarrow\mathbf{P}_{k}^{N_{i}}$ un plongement de $G_{i}$ dans l'espace projectif $\mathbf{P}_{k}^{N_{i}}$. En particulier, les degr{\'e}s et fonctions de Hilbert-Samuel g{\'e}om{\'e}triques consid{\'e}r{\'e}s dans la suite sont relatifs aux faisceaux $\mathcal{O}_{G_{i}}(1)$, $1\le i\le n$, induits par ces plongements. Notons $G$ le groupe $G_{1}\times\cdots\times G_{n}$, $g:=g_{1}+\cdots+g_{n}$ sa dimension et $\Phi$ le plongement de $G$ dans le produit des $\mathbf{P}_{k}^{N_{i}}$ induit par les $\Phi_{i}$. Soit $\mathcal{G}_{i}\to\spec\mathcal{O}_{k}[1/m_{i}]$ ($m_{i}\in\mathbf{N}\setminus\{0\}$) un sch\'ema en groupes lisse et dont la fibre g{\'e}n{\'e}rique $\mathcal{G}_{i}\times\spec(k)$ est (isomorphe {\`a}) $G_{i}$. Quitte {\`a} restreindre $\mathcal{G}_{i}$ {\`a} un ouvert plus petit, nous pouvons supposer d'une part que $m_{i}$ est le m{\^e}me entier $m$ pour tous les $i$ et d'autre part que l'anneau $\mathcal{O}_{k}[1/m]$ est principal\footnote{Si bien qu'un $\mathcal{O}_{k}[1/m]$-module projectif (de type fini) est n{\'e}cessairement libre.}. Soit $\mathcal{G}\to\spec\mathcal{O}_{k}[1/m]$ le mod{\`e}le lisse de $G$ induit par les $\mathcal{G}_{i}$. Fixons $v$ une place de $k$ et $i\in\{1,\ldots,n\}$. Consid{\'e}rons $\exp_{i,v}$ une application exponentielle du groupe de Lie $v$-adique $G_{i}(\mathbf{C}_{v})$, d{\'e}finie sur un voisinage ouvert de $0$ dans $t_{G_{i}}(\mathbf{C}_{v})$. Lorsque $v$ est une place archim{\'e}dienne, il est bien connu que cette application se prolonge en un morphisme analytique {\`a} tout l'espace tangent $t_{G_{i}}(\mathbf{C}_{v})$ et d{\'e}finit ainsi une application $\mathbf{C}_{v}$-analytique $\exp_{i,v}:t_{G_{i}}(\mathbf{C}_{v})\to G_{i}(\mathbf{C}_{v})$ surjective. Lorsque $v$ est ultram{\'e}trique, ces propri{\'e}t{\'e}s ne sont plus vraies en g{\'e}n{\'e}ral. Notons alors dans ce cas $\mathscr{T}_{i,v}$ un sous-groupe ouvert de $t_{G_{i}}(\mathbf{C}_{v})$ tel que $\exp_{i,v}$ r{\'e}alise un diff{\'e}omorphisme analytique de $\mathscr{T}_{i,v}$ sur son image $\mathscr{U}_{i,v}$ (ouvert de $G_{i}(\mathbf{C}_{v})$ contenant l'{\'e}l{\'e}ment neutre). Dans la suite, l'exponentielle $\exp_{i,v}$ sera l'application restreinte $\mathscr{T}_{i,v}\to\mathscr{U}_{i,v}$. Afin d'uniformiser les notations, nous {\'e}crirons encore $\mathscr{T}_{i,v}=t_{G_{i}}(\mathbf{C}_{v})$ (\emph{resp}. $\mathscr{U}_{i,v}=G_{i}(\mathbf{C}_{v})$) dans le cas archim{\'e}dien, bien que l'exponentielle \og restreinte\fg\ ne soit plus alors un diff{\'e}omorphisme (en g{\'e}n{\'e}ral). L'espace tangent {\`a} l'origine $t_{\mathcal{G}_{i}}$ de $\mathcal{G}_{i}$ est un $\mathcal{O}_{k}[1/m]$-module libre (car projectif, voir note de bas de page) de rang $g_{i}$ et $t_{\mathcal{G}}$ est {\'e}galement libre de rang $g$. Soit $\mathbf{e}=(e_{1},\ldots,e_{g})$ une base sur $\mathcal{O}_{k}[1/m]$ de $t_{\mathcal{G}}$ obtenue par concat{\'e}nation de bases des $t_{\mathcal{G}_{i}}$, $1\le i\le n$. Quitte {\`a} multiplier chacun des $e_{i}$ par une puissance suffisamment grande de $m$, nous pouvons supposer que, pour toute place ultram{\'e}trique $v$, le disque ouvert \begin{equation*}D(0,r_{p})=\left\{\mathbf{z}=z_{1}e_{1}+\cdots+z_{g}e_{g}\in t_{G}(\mathbf{C}_{v})\,;\ \max_{1\le j\le g}{\vert z_{j}\vert_{v}}<r_{p}\right\},\end{equation*}o\`u \label{pagesept}$r_{p}:=\vert p\vert_{v}^{1/(p-1)}$, est inclus dans $\mathscr{T}_{v}:=\mathscr{T}_{1,v}\times\cdots\times\mathscr{T}_{n,v}$. En effet, il existe un entier $n_{0}\ge 1$, ne d{\'e}pendant que de $(\mathcal{G},m)$, pour lequel, en toute place $v$, le d{\'e}veloppement en s{\'e}rie de l'exponentielle de $G(\mathbf{C}_{v})$ au voisinage de $0$ s'{\'e}crit $$\sum_{\mathbf{n}\in\mathbf{N}^{g}}{\frac{a_{\mathbf{n},v}}{\mathbf{n}!}\mathbf{z}^{\mathbf{n}}}$$avec $a_{\mathbf{n},v}\in\mathcal{O}_{v}[1/m]$, polyn{\^o}me en $1/m$ de degr{\'e} $\le n_{0}\vert\mathbf{n}\vert$. Aux places $v\nmid m$, on a $a_{\mathbf{n},v}\in\mathcal{O}_{v}$ et l'on sait que le disque $D(0,r_{p})$ est contenu dans le domaine de convergence strict de cette s{\'e}rie. Si $v\mid m$, on se ram{\`e}ne au cas pr{\'e}c{\'e}dent en consid{\'e}rant les coordonn{\'e}es de $\mathbf{z}$ dans la base $m^{n_{0}}\mathbf{e}$ de $t_{G}(\mathbf{C}_{v})$.\par L'exponentielle $\exp_{v}:=(\exp_{1,v},\ldots,\exp_{n,v})$ de $G(\mathbf{C}_{v})$ munie de la base $\mathbf{e}$ est appel{\'e}e dans la litt{\'e}rature \emph{exponentielle normalis{\'e}e} (cela fixe un isomorphisme de $\mathscr{U}_{v}:=\mathscr{U}_{1,v}\times\cdots\times\mathscr{U}_{n,v}$ avec un groupe \emph{standard} selon la terminologie de Bourbaki~\cite{bourbagroupedeLie}, III, \S~$7$, n°$3$). La base $\mathbf{e}$ conf{\`e}re {\'e}galement {\`a} $t_{G}(\mathbf{C}_{v})$ une structure d'espace vectoriel norm{\'e}, par transport de la structure hermitienne ($v$ archim{\'e}dienne) ou de la norme du sup ($v$ ultram{\'e}trique) fournie par la base canonique de $\mathbf{C}_{v}^{g}$ (voir \S~Notations et Conventions). Nous notons $\Vert\cdotp\Vert_{v}$ (\emph{resp}. $\mathrm{d}_{v}$) cette norme sur $t_{G}(\mathbf{C}_{v})$ (\emph{resp}. la distance associ{\'e}e {\`a} cette norme).
\par Consid{\'e}rons un point\footnote{La lettre $p$ d{\'e}signe {\`a} la fois ce point et le nombre premier qui divise $v_{0}$. Mais cette maladresse ne devrait pas cr{\'e}er d'ambigu{\"i}t{\'e}.} $p=(p_{1},\ldots,p_{n})$ de $G(k)\cap\mathscr{U}_{v_{0}}$ ainsi qu'un logari\-th\-me $u=(u_{1},\ldots,u_{n})\in\mathscr{T}_{v_{0}}$ de ce point~: \begin{equation*}\exp_{v_{0}}(u)=p\ .\end{equation*}\par Soit $V$ un sous-espace vectoriel de l'espace tangent $t_{G}(k)$, de codimension $t\ge 1$ (ce qui suit est trivial et d{\'e}nu{\'e} d'int{\'e}r{\^e}t lorsque $t=0$, c.{\`a}-d. $V=t_{G}(k)$).\\[0.1cm]\textsc{Hypoth{\`e}se}~: \textbf{Dans tout ce texte, nous supposerons que $V$ est l'alg{\`e}bre de Lie d'un sous-groupe alg{\'e}brique connexe $H$ de $G$.}\\[0.1cm] En d'autres termes, l'espace $V$ est une alg\`ebre de Lie \emph{alg\'ebrique} au sens de~\cite{demazure} (II, \S~$6$, $2.4$, p.~$262$).
\paragraph{\textsc{Importante convention}} Dans toute la suite, le mot \og constante\fg\ qualifie un nombre r{\'e}el $\ge 1$ qui ne d{\'e}pend que de $G,\Phi,\mathcal{G},m,\mathbf{e},(d_{v})_{v},v_{0}$, c'est-{\`a}-dire du groupe alg{\'e}brique $G$ et des donn{\'e}es satellites autour de $G$. Partant, ce nombre r{\'e}el est ind{\'e}pendant de $k$ (il ne d{\'e}pend que d'un corps de d{\'e}finition de $G$), de $p,u,V$ etc. Une telle constante sera d{\'e}sign{\'e}e par la lettre $c$ munie d'un indice.

\subsection{R{\'e}sultats}
\label{subsec:resultats}
Les th{\'e}or{\`e}mes que nous allons {\'e}noncer ici concernent tous le cas rationnel, comme nous l'avons mentionn{\'e} dans l'introduction.\par Fixons auparavant quelques notations suppl\'ementaires. Pour $i\in\{1,\ldots,n\}$, le plongement (en tant qu'espace quasi-projectif) de $G_{i}$ dans $\mathbf{P}_{k}^{N_{i}}$ fournit une hauteur de Weil $h$ sur l'ensemble des points $\overline{\mathbf{Q}}$-rationnels de $G_{i}$ (dont, en particulier, $p_{i}$). Nous noterons encore $h$ la hauteur induite sur $G(\overline{\mathbf{Q}})$. Soit $\rho_{i}$ l'ordre analytique de $G_{i}$ d\'efini comme suit~: $\rho_{i}=1$ si $G_{i}$ est un groupe lin\'eaire et $\rho_{i}=2$ sinon (c'est-\`a-dire lorsque $G_{i}$ a une composante ab\'elienne non triviale). Soit $y\in\{0,1\}$ un param{\`e}tre tel que $y=0$ lorsque $G$ est une vari{\'e}t{\'e} semi-ab{\'e}lienne et $y=1$ sinon.
\par Dans le cas archim\'edien, l'{\'e}nonc{\'e} le plus g{\'e}n{\'e}ral que nous obtenons est le suivant.

\begin{theo}\label{theounasept}
Il existe une constante $c_{\cst}\ge 1$\newcounter{acstsept}\setcounter{acstsept}{\value{constante}} ayant la propri{\'e}t{\'e} suivante. Supposons que $v_{0}$ est une place archim{\'e}dienne. Soit $\mathfrak{e}$ un nombre r{\'e}el $\ge e$ et $\mathfrak{a}$ un entier naturel sup{\'e}rieur ou {\'e}gal {\`a} $D\max{\{1,h(V)\}}/\log\mathfrak{e}$. Notons $U_{0}$ le nombre r{\'e}el \begin{equation}\begin{split} U_{0}:=&(\mathfrak{a}\log\mathfrak{e})\left(\mathfrak{a}^{y}+\frac{D}{\log\mathfrak{e}}\log\left(e+\frac{D}{\log\mathfrak{e}}\right)\right)^{1/t}\\ &\times\prod_{i=1}^{n}{\left(1+\frac{D\underset{0\le s\le c_{\theconstante}\mathfrak{a}}{\mathrm{max}}\{h(sp_{i})\}+(\mathfrak{e}\mathfrak{a}\Vert u_{i}\Vert_{v_{0}})^{\rho_{i}}}{\mathfrak{a}\log\mathfrak{e}}\right)^{g_{i}/t}}\ \cdotp\end{split}\end{equation}Supposons que pour tout entier $s\in\{1,\ldots,c_{\theconstante}\mathfrak{a}\}$ et tout sous-groupe alg{\'e}brique connexe $G'$ de $G$ v{\'e}rifiant $t_{G'}+V\ne t_{G}$ on ait $s(p_{1},\ldots,p_{n})\not\in G'(\overline{k})$. Alors $u\not\in V\otimes_{v_{0}}\mathbf{C}$ et \begin{equation}\label{ineq:principalea7}\log\mathrm{d}_{v_{0}}(u,V)\ge -c_{\theconstante}U_{0}\ \cdotp\end{equation}
\end{theo}
Formellement, ce r{\'e}sultat est tr{\`e}s proche de celui {\'e}nonc{\'e} avec $t=1$ dans le th{\'e}or{\`e}me principal de~\cite{gaudron1}. Hormis la disparition de la hauteur d'une $\mathbf{Q}$-base de $k$, c'est surtout la d{\'e}finition de l'entier $\mathfrak{a}$ qui change radicalement. Dans l'article en question, cet entier {\'e}tait (en substance) $\frac{D}{\log\mathfrak{e}}\log h(p)$ alors qu'ici il d{\'e}pend du sous-espace $V$ mais pas du point $p$. Sachant qu'il existe une constante $c_{\cst}$ pour laquelle $h(sp_{i})\le c_{\theconstante}s^{\rho_{i}}\max{\{1,h(p_{i})\}}$ pour tout entier $s\ge 1$, l'on d{\'e}duit ais{\'e}ment de la d{\'e}finition de $U_{0}$ la d{\'e}pendance standard en $h(p)$ d{\'e}crite dans l'introduction. Si $n=1$, nous pouvons regarder la d{\'e}pendance \emph{minimale} en $h(p)$ de $U_{0}$, c'est-\`a-dire choisir $\mathfrak{e}$ (qui est le seul param{\`e}tre vraiment \og libre\fg\ du th{\'e}or{\`e}me~\ref{theounasept}) de sorte que $U_{0}$, comme fonction uniquement de $h(p)$, soit minimal. Avec des consid{\'e}rations {\'e}l{\'e}mentaires, on s'aper{\c c}oit que, dans cette optique, le meilleur choix pour $\mathfrak{e}$ est $e\sqrt{\max{\{1,h(p)\}}}$, ce qui conduit {\`a} l'estimation \begin{equation}\label{eq:dependanceenhp}\log\mathrm{d}_{v_{0}}(u,V)\ge -c_{\cst}\left(\frac{\max{\{1,h(p)\}}}{\log\max{\{e,h(p)\}}}\right)^{g/t}\end{equation}o\`u $c_{\theconstante}$ est une fonction des donn{\'e}es qui ne d{\'e}pend pas de la hauteur de $p$.  Si nous comparons cela {\`a} la cons{\'e}quence I.3.3 de~\cite{gaudron1} ({\'e}crite avec $t=1$), nous constatons {\`a} nouveau la disparition d'un logarithme de $h(p)$ (qui {\'e}tait au num{\'e}rateur du membre de droite de~\eqref{eq:dependanceenhp}). Autrement dit, nous v{\'e}rifions ainsi que l'am{\'e}lioration en $h(p)$ est bien une caract{\'e}ristique intrins{\`e}que de la mesure~\eqref{ineq:principalea7} et non pas le simple effet d'un choix diff{\'e}rent de param{\`e}tres dans deux {\'e}nonc{\'e}s \og semblables\fg.
\par Nous avons {\'e}galement une version ultram{\'e}trique de l'{\'e}nonc{\'e}~\ref{theounasept}.
\begin{theo}\label{theodeuxasept}
Il existe une constante $c_{\cst}\ge 1$\newcounter{acsthuit}\setcounter{acsthuit}{\value{constante}} ayant la propri{\'e}t{\'e} suivante. Supposons que $v_{0}$ est ultram{\'e}trique et que $\Vert u\Vert_{v_{0}}<r_{p}^{2}$ (o\`u, rappelons-le, $p$ est la caract\'eristique r\'esiduelle de $v_{0}$ et $r_{p}=p^{-1/(p-1)}$). Consid{\'e}rons un nombre r{\'e}el $\mathfrak{r}$ dans l'intervalle ouvert $]1,r_{p}^{2}/\Vert u\Vert_{v_{0}}[$. Soit $\mathfrak{a}$ un entier naturel v{\'e}rifiant
\begin{equation*}\mathfrak{a}\ge\frac{D\max{\{1,h(V)\}}+\log^{+}((\log(\mathfrak{r}))^{-1})}{\log\mathfrak{r}}\ \cdotp\end{equation*}Notons $U_{1}$ le nombre r{\'e}el \begin{equation}(\mathfrak{a}\log(1+\mathfrak{r}))\left(\mathfrak{a}^{y}+\frac{D}{\log\mathfrak{r}}\log\left(e+\frac{D}{\log\mathfrak{r}}\right)\right)^{1/t}\prod_{i=1}^{n}{\left(1+\frac{D\underset{0\le s\le c_{\theconstante}\mathfrak{a}}{\mathrm{max}}\{h(sp_{i})\}}{\mathfrak{a}\log\mathfrak{r}}\right)^{g_{i}/t}}\ \cdotp\end{equation}Supposons que pour tout entier $s\in\{1,\ldots,c_{\theconstante}\mathfrak{a}\}$ et tout sous-groupe alg{\'e}brique connexe $G'$ de $G$ tel que $t_{G'}+V\ne t_{G}$ on ait $s(p_{1},\ldots,p_{n})\not\in G'(\overline{k})$. Alors $u\not\in V\otimes\mathbf{C}_{v_{0}}$ et $\log\mathrm{d}_{v_{0}}(u,V)\ge-c_{\theconstante}U_{1}$.
\end{theo}
\begin{remas}
\noignorespaces
\begin{enumerate}
\item[1)] Les diff{\'e}rences entre les versions archim{\'e}dienne et ultram{\'e}trique r{\'e}sident d'une part dans le changement de $\mathfrak{e}$ par $\mathfrak{r}$ et d'autre part dans une contrainte plus forte sur l'entier $\mathfrak{a}$ dans le cas ultram{\'e}trique. Par ailleurs, toujours dans ce cas, la norme $\Vert u_{j}\Vert_{v_{0}}$ du logarithme $u_{j}$ n'appara{\^\i}t plus.
\item[2)] Soulignons que cette minoration ne d{\'e}voile pas la d{\'e}pendance en la place $v_{0}$, dont une partie est dans la constante $c_{\theconstante}$, non explicite.
\end{enumerate}
\end{remas}
Si la litt{\'e}rature est assez riche et vari{\'e}e en analogues $p$-adiques de mesures d'ind{\'e}pendance lin{\'e}aire de logarithmes lorsque $G$ est un groupe lin{\'e}aire (en particulier gr{\^a}ce aux travaux de Yu~\cite{Yu1998,Yu1999,Yu2007}), elle est en revanche beaucoup plus r{\'e}duite si $G$ a une partie ab{\'e}lienne, voire m{\^e}me inexistante lorsque, comme ici, $G$ est quelconque. Mentionnons la s\'erie d'articles de Bertrand~\cite{bertrand77,bertrand78,bertrandflicker} (la derni{\`e}re r{\'e}f{\'e}rence est un article en commun avec Yu.~Flicker) \`a la fin des ann\'ees~$70$, ainsi que le r\'esultat de R\'emond \& Urfels~\cite{gael-urfels} qui traite le cas du produit de \emph{deux} courbes elliptiques.
\section{Quelques mots sur la d{\'e}monstration des th\'eor\`emes~\ref{theounasept} et~\ref{theodeuxasept}}
Nous n'allons pas expliquer ici le sch\'ema de la d\'emonstration, somme toute assez classique, fond\'e sur la m{\'e}thode de Baker revisit{\'e}e et approfondie par Philippon \& Waldschmidt~\cite{pph-miw}. La d{\'e}marche est rappel{\'e}e au d{\'e}but du \S~\ref{secasept:demonstrationdestheos}. Nous voulons plut{\^o}t d{\'e}gager de fa{\c c}on {\'e}l{\'e}mentaire \emph{la} difficult{\'e} technique sur laquelle achoppaient les preuves dans le cas du tore pour le passage {\`a} un groupe alg{\'e}brique quelconque. Nous en profiterons {\'e}galement pour mettre en lumi{\`e}re certaines modifications techniques de la d\'emonstration, qui la simplifient (dans une certaine mesure), mais au prix, il est vrai, de l'hypoth{\`e}se sur le point $p$ d\'ej\`a rencontr{\'e}e dans l'{\'e}nonc{\'e}~\ref{enonce2demoi}.\par Commen{\c c}ons donc par expliquer l'id{\'e}e fondamentale du cas rationnel usuel dans $\mathbb{G}_{\mathrm{m}}^{n}$ qui conduit \`a de meilleures mesures d'ind{\'e}pendances lin{\'e}aires de logarithmes. Pour cela, simplifions la situation au maximum et ne conservons que $G=\mathbb{G}_{\mathrm{m}}^{2}$, la forme lin{\'e}aire $z_{2}-bz_{1}$ ($b\in\mathbf{Z}$, $z_{1},z_{2}$ coordonn{\'e}es sur $\lie(G)$) et le point $p=(\alpha_{1},\alpha_{2})\in(k\setminus\{0\})^{2}$ de logarithme $(u_{1},u_{2})$. Autrement dit, nous nous int{\'e}ressons {\`a} la forme lin{\'e}aire en deux logarithmes $\Lambda=u_{2}-bu_{1}$. Les preuves \og classiques\fg\ qui m\`enent \`a une minoration de $\vert\Lambda\vert_{v_{0}}$ reposent sur l'{\'e}tude des d{\'e}riv{\'e}es (divis{\'e}es) le long de la droite $z_{2}=bz_{1}$ d'un certain polyn{\^o}me exponentiel $(z_{1},z_{2})\to P(e^{z_{1}},e^{z_{2}})$ ($P\in k[X,Y]$) en les points $s\cdotp(u_{1},u_{2})$, $s\in\mathbf{N}$. Par translation sur le groupe $G$, l'on peut se ramener {\`a} $s=0$ et le r{\'e}sultat clef qui permet d'exploiter l'hypoth{\`e}se $b\in\mathbf{Z}$ est le suivant. \begin{fait}\label{faitnumeroun}Soit $P\in\mathbf{Z}[X,Y]$, $b$ un entier et $\ell$ un entier naturel. Supposons que l'application $z\mapsto P(e^{z},e^{bz})$, analytique au voisinage de $0$, s'annule \`a l'ordre $\ell$ en $0$. Alors le nombre \begin{equation}\label{eqfaitcoeff}\frac{1}{\ell !}\left(\frac{\mathrm{d}}{\mathrm{d}z}\right)^{\ell}P(e^{z},e^{bz})_{\vert z=0}\end{equation}est un entier relatif. \end{fait}Il y a au moins deux preuves assez diff{\'e}rentes de ce r{\'e}sultat. La premi{\`e}re utilise les polyn{\^o}mes binomiaux\begin{equation}\label{binomiaux}\Delta_{0}(X):=1,\quad\Delta_{n}(X):=\frac{X(X+1)\cdots (X+n-1)}{n!}\quad (n\in\mathbf{N}\setminus\{0\})\end{equation}qui prennent des valeurs enti{\`e}res aux points entiers. En consid{\'e}rant un mon{\^o}me $X^{i}Y^{j}$ qui intervient dans $P$ (avec le coefficient $p_{i,j}\in\mathbf{Z}$), la d{\'e}riv{\'e}e $\ell^{\,\text{{\`e}me}}$ de $z\mapsto e^{(i+jb)z}$ en $0$ vaut $(i+jb)^{\ell}/\ell !$. Ce terme est la somme de $\Delta_{\ell}(i+jb)$ et d'une combinaison lin{\'e}aire de $(i+jb)^{h}$ avec $h<\ell$ dont les coefficients ne d{\'e}pendent que de $\ell$. L'hypoth{\`e}se sur $P$ se traduit alors par l'{\'e}galit{\'e} entre le coefficient~\eqref{eqfaitcoeff} et \begin{equation}\label{deuxeqfaitcoeff}\sum_{i,j}{p_{i,j}\Delta_{\ell}(i+jb)}\ ,\end{equation}manifestement un entier, ce qui conclut la preuve. Aussi astucieux soit-il, ce proc{\'e}d{\'e} comporte n{\'e}anmoins une limitation consubstantielle puisque $b$ ne peut {\^e}tre qu'un entier (ou au pire un nombre rationnel), faute de quoi il est difficile d'envisager une g{\'e}n{\'e}ralisation. La seconde preuve du fait~\ref{faitnumeroun} que nous connaissons est bas{\'e}e sur un changement de variables. On pose $T=e^{z}-1$. Comme $b\in\mathbf{Z}$, chacune des fonctions $e^{(i+jb)z}=(1+T)^{i+jb}$ appartient {\`a} l'anneau de s{\'e}ries formelles $\mathbf{Z}[[T]]$. Il en est donc de m{\^e}me pour $P(e^{z},e^{bz})$ et l'hypoth{\`e}se sur $P$ implique que le coefficient~\eqref{eqfaitcoeff} est {\'e}galement le coefficient de $T^{\ell}$ dans ce d{\'e}veloppement. C'est donc un entier. Contrairement {\`a} la d{\'e}monstration pr{\'e}c{\'e}dente, cette m{\'e}thode peut \^etre g{\'e}n{\'e}ralis{\'e}e {\`a} un groupe alg{\'e}brique quelconque. D'ailleurs, rappelons que ce passage de la variable $z$ (sur $\lie(\mathbb{G}_{\mathrm{m}})$) {\`a} la variable $T$ (sur $\mathbb{G}_{\mathrm{m}}$) et ses r\'epercussions arithm\'etiques constituent la cheville ouvri{\`e}re des r{\'e}centes avanc{\'e}es dans le domaine des formes lin{\'e}aires de logarithmes (voir\cite{davidhirata,gaudron1}), mais aussi dans les questions li{\'e}es {\`a} l'alg{\'e}bricit{\'e} de feuilles formelles~\cite{graftieaux2,bost6}. C'est cette observation qui apporte l'essentiel des r{\'e}sultats nouveaux de cet article. \par Il me faut signaler cependant qu'une difficult{\'e} technique {\'e}chappe {\`a} l'analyse du cas d'un tore telle que nous venons de la faire. Dans le cas g{\'e}n{\'e}ral, nous avons besoin de mod{\`e}les lisses des groupes $G$ et $H$ (rappelons que $V=\lie(H)$) sur des anneaux de la forme $\mathcal{O}_{k}[1/m]$, o\`u $m$ est un entier $>0$. Si pour le groupe $G$ cela ne pose aucun probl{\`e}me ($m$ d{\'e}pend de $G$), le groupe $H$, quant {\`a} lui, admet un mod{\`e}le lisse mais sur un anneau localis\'e $\mathcal{O}_{k}[1/mm']$ avec $m'$ un entier qui d{\'e}pend \emph{a priori} de $H$. Par cons{\'e}quent, il est important de contr{\^o}ler l'entier $m'$ fonction du mod{\`e}le de $H$. Cela revient {\`a} avoir des estimations $p$-adiques de nombres alg{\'e}briques plus g{\'e}n{\'e}raux issus de~\eqref{eqfaitcoeff} qui soient les plus pr{\'e}cises possible et qui tiennent compte du mod{\`e}le choisi pour $H$. C'est pourquoi nous emploierons un formalisme particuli{\`e}rement adapt{\'e} {\`a} cette exigence, d{\'e}crit par Bost au \S~$3.1$ de~\cite{bost6}. Le langage g{\'e}om{\'e}trique de ce formalisme, qui s'exprime en termes de \og tailles de sch{\'e}mas formels lisses\fg, {\'e}claire le r{\^o}le exact jou{\'e} par le choix du mod{\`e}le de $H$. Mais la difficult{\'e} technique {\'e}voqu{\'e}e ne dispara{\^\i}t pas pour autant dans ce langage. Un th{\'e}or{\`e}me de Raynaud, donnant un condition pour que l'inclusion entre vari{\'e}t{\'e}s ab{\'e}liennes se prolonge en une \emph{immersion ferm\'ee} pour les mod{\`e}les de N{\'e}ron correspondants, permet alors de contr{\^o}ler tr{\`e}s pr{\'e}cis{\'e}ment l'entier $m'$. Nous d{\'e}taillerons tout cela au \S~\ref{subsec:estultracoefftay}.\par Cet argument arithm{\'e}tique crucial s'accompagne d'une double utilisation d'un lemme de Siegel absolu, {\`a} la fois pour b{\^a}tir le \og classique\fg\ polyn{\^o}me auxiliaire requis par la d{\'e}monstration de transcendance mais aussi pour fixer une $k$-base de $V=\lie(H)$, de \og petite\fg\ hauteur, qui restera la m{\^e}me {\`a} chaque {\'e}tape de la preuve. Ce dernier point {\'e}vite le recours {\`a} certaines bases orthonorm{\'e}es de $\lie(G)\otimes_{v}\mathbf{C}$ et les contr{\^o}les de changements de bases subs{\'e}quents, qui intervenaient auparavant. Quant {\`a} construire le polyn{\^o}me auxiliaire de la sorte, cela procure l'avantage de supprimer la quantit{\'e} $Dh(\xi_{1},\ldots,\xi_{D})$, o\`u $\xi_{1},\ldots,\xi_{D}$ est une $\mathbf{Q}$-base de $k$ (\emph{i.e.}, de mani\`ere \'equivalente, le logarithme du discriminant absolu de $k$, cf.~\cite{RoyThunderRocky}), qui apparaissait dans les mesures de~\cite{gaudron1}. L'emploi d'un lemme de Siegel absolu dans le contexte des formes lin{\'e}aires de logarithmes nous avait {\'e}t{\'e} communiqu{\'e} par S.~David (voir~\cite{crasdeux,davidhirata}). Il remplace le lemme de Thue-Siegel dont on se servait d'ordinaire. Nous le pr{\'e}sentons au \S~\ref{subsec:lemmesiegelabsolu} et nous l'appliquons aux \S~\ref{subsec:constructionbaseW} et~\ref{subsec:constructionasept}. Tous les bienfaits de ce lemme pour la d{\'e}monstration ({\`a} commencer par la clart{\'e} m\^eme de l'argumentation) sont malheureusement un peu ternis par une difficult{\'e} technique que je ne sais pas surmonter sans supposer que le groupe engendr{\'e} par le point $p$ ne rencontre aucun sous-groupe strict de $G(\overline{k})$ sauf en $0$. Nous avions d{\'e}j{\`a} {\'e}t{\'e} contraint d'{\'e}mettre ce type d'hypoth{\`e}ses lorsque nous avions mis en \oe uvre la m{\'e}thode des pentes, qui elle, pourtant, ne requiert aucun lemme de Siegel (voir~\cite{artepredeux}). Le point commun aux deux approches est un certain sous-cas, appel{\'e} \emph{cas p{\'e}riodique}, que je ne sais pas int\'egrer dans les preuves, bien qu'il f{\^u}t d\'ej\`a r{\'e}solu de mani{\`e}re tr{\`e}s astucieuse par Philippon \& Waldschmidt dans leur article~\cite{pph-miw}, gr{\^a}ce {\`a} une extrapolation ({\`a} la mani\`ere de Gel'fond) sur les d{\'e}rivations.

\section{Pr{\'e}paratifs} \label{paragraphepreparatifasept} 
La d{\'e}monstration des th{\'e}or{\`e}mes~\ref{theounasept} et~\ref{theodeuxasept} requiert plusieurs {\'e}nonc{\'e}s d'int{\'e}r{\^e}ts ind{\'e}pendants que nous pr{\'e}sentons dans cette partie.

\subsection{Mise en place de donn{\'e}es suppl{\'e}mentaires}\label{miseenplace}
Consid\'erons le groupe $\mathbb{G}_{\mathrm{a}}\times G$, le point $q=(1,p)\in(\mathbb{G}_{\mathrm{a}}\times G)(k)$ et le sous-espace vectoriel $W$ de $t_{\mathbb{G}_{\mathrm{a}}}\oplus t_{G}$ d{\'e}fini par $t_{\mathbb{G}_{\mathrm{a}}}\oplus V$. L'{\'e}l{\'e}ment $1\oplus u\in t_{\mathbb{G}_{\mathrm{a}}}(\mathbf{C}_{v_{0}})\oplus \mathscr{T}_{v_{0}}$ est un logarithme du point $q$. Pour uniformiser les notations, nous posons $G_{0}:=\mathbb{G}_{\mathrm{a}}$ et $u_{0}:=1$. Par ailleurs, consid{\'e}rons un entier $i$ compris entre $1$ et $n$. Soit $\Phi_{i}$ un plongement de $G_{i}$ dans l'espace projectif $\mathbf{P}^{N_{i}}_{k}$ du type de ceux construits par Serre~\cite{Serre1}. Quitte {\`a} effectuer un changement de coordonn{\'e}es, nous pouvons supposer que l'{\'e}l{\'e}ment neutre de $G_{i}$ est repr{\'e}sent{\'e} par $(1:0:\cdots:0)\in\mathbf{P}^{N_{i}}$. Soit $x\in G_{i}(\mathbf{C}_{v_{0}})$ et $(x_{0}:\cdots:x_{N_{i}})$ les coordonn{\'e}es de $\Phi_{i}(x)$. On note \begin{equation}\label{systemeaddition} A_{x}^{(i)}=(A_{x,0}^{(i)}(\mathbf{X},\mathbf{Y}):\cdots:A_{x,N_{i}}^{(i)}(\mathbf{X},\mathbf{Y}))\end{equation}une famille de polyn{\^o}mes ({\`a} coefficients dans $\mathcal{O}_{v_{0}(k)}$) qui exprime la loi d'addition de $G_{i}$ au voisinage de $x$. Dans cette formule, $\mathbf{X}$ et $\mathbf{Y}$ sont des $(N_{i}+1)$-uplets de variables et chacun des $A_{x,j}^{(i)}(\mathbf{X},\mathbf{Y})$, $0\le j\le N_{i}$, est homog{\`e}ne de m{\^e}me degr{\'e} sur chacune des variables $\mathbf{X},\mathbf{Y}$, inf{\'e}rieur {\`a} une constante $c_{\cst}$\newcounter{degdespolys}\setcounter{degdespolys}{\value{constante}}, qui peut \^etre choisie uniforme en $x$ (quasi-compacit{\'e} de $G_{i}$) et en $i$ (nombre fini). Cette constante ne d{\'e}pend que de $(G,\Phi)$. Dans la suite, pour ne pas alourdir excessivement les notations, nous omettrons souvent la r{\'e}f{\'e}rence {\`a} $x$ en indice et nous {\'e}crirons $A_{j}^{(i)}$ au lieu de $A_{x,j}^{(i)}$. Soit $v$ une place quelconque de $K$. Il est possible {\'e}galement de repr{\'e}senter l'exponentielle $v$-adique de $G_{i}(\mathbf{C}_{v})$ par des fonctions $(\theta_{v,i,j})_{0\le j\le N_{i}}$, analytiques et sans z{\'e}ros communs dans $\mathscr{T}_{i,v}$, telles que $(\theta_{v,i,0}(0),\ldots,\theta_{v,i,N_{i}}(0))=(1,0,\ldots,0)$~: \begin{equation*}\exp_{i,v}(z)=(\theta_{v,i,0}(z):\cdots:\theta_{v,i,N_{i}}(z)),\quad z\in\mathscr{T}_{i,v}\ .\end{equation*}Nous noterons \begin{equation*}\Theta_{v,i}:z\in\mathscr{T}_{i,v}\mapsto(\theta_{v,i,0}(z),\ldots,\theta_{v,i,N_{i}}(z))\in\mathbf{C}_{v}^{N_{i}+1}\end{equation*}et, si $j$ est un entier naturel inf\'erieur \`a $N_{i}$, \begin{equation*}\Psi_{v,i,j}:z\in\mathscr{T}_{i,v}\setminus\theta_{v,i,j}^{-1}(\{0\})\mapsto\left(\frac{\theta_{v,i,0}}{\theta_{v,i,j}}(z),\ldots,\frac{\theta_{v,i,N_{i}}}{\theta_{v,i,j}}(z)\right)\in\mathbf{C}_{v}^{N_{i}+1}\ \cdotp\end{equation*}Bien que ce soit un abus de notations, nous nous permettrons d'\'ecrire ces formules pour $z\in\mathscr{T}_{v}$ au lieu de la $i^{\text{\`eme}}$ composante de $z$ sur $\mathscr{T}_{i,v}$. Comme nous l'avons vu au \S~\ref{donnesgeneralesasept}, lorsque $v$ est ultram{\'e}trique, le choix de la base $\mathbf{e}$ permet d'{\'e}crire chacune des coordonn{\'e}es $\theta_{v,i,j}(z)$ sous la forme d'une s{\'e}rie $\sum_{\mathbf{n}}{\frac{a_{\mathbf{n},v,i,j}}{\mathbf{n}!}\mathbf{z}^{\mathbf{n}}}$ o\`u $\mathbf{z}=(z_{1},\ldots,z_{g_{i}})$ sont les coordonn{\'e}es de $z$ dans la base $\mathbf{e}$ et $a_{\mathbf{n},v,i,j}\in\mathcal{O}_{v}$. De plus, lorsque $v$ est archim\'edienne, les fonctions $\theta_{v,i,j}$, $0\le j\le N_{i}$, sont d'ordre analytique $\le\rho_{i}$ et il existe une constante $c_{\cst}\ge 1$\newcounter{prime}\setcounter{prime}{\value{constante}} telle que, pour tout entier $i\in\{1,\ldots,n\}$, pour tout vecteur $z$ de $\mathscr{T}_{i,v}$, on ait \begin{equation}\label{asepteq:inegalitetreize}-c_{\theprime}(1+\Vert z\Vert_{v})^{\rho_{i}}\le\log\max_{0\le j\le N_{i}}\vert\theta_{v,i,j}(z)\vert\le c_{\theprime}(1+\Vert z\Vert_{v})^{\rho_{i}}\ .\end{equation}La $k$-structure de l'espace tangent $t_{G_{i}}$ entra{\^\i}ne une stabilit{\'e} par d{\'e}rivation (selon un vecteur de $t_{G_{i}}(k)$) de l'anneau $k[(\theta_{v,i,j}/\theta_{v,i,0})_{0\le j\le N_{i}}]$. Ces propri{\'e}t{\'e}s seront utilis{\'e}es aux paragraphes~\ref{subsec:estultracoefftay} et~\ref{subsec:estarchicoefftay}.\par
Dans la suite nous noterons $\mathbf{P}$ l'espace multiprojectif $\mathbf{P}^{1}\times\mathbf{P}^{N_{1}}\times\cdots\times\mathbf{P}^{N_{n}}$ (le corps de base \'etant $k$, $\overline{k}$ ou $\mathbf{C}_{v_{0}}$ selon le contexte). Il est naturellement muni du faisceau canonique $\mathcal{O}_{\mathbf{P}}(1,\ldots,1)$ et, si $\mathsf{V}$ est une sous-vari\'et\'e (ferm\'ee) de $\mathbf{P}$, l'entier $\deg\mathsf{V}$ (\emph{resp}. le polyn\^ome $H(\mathsf{V};X_{0},\ldots,X_{n})$) d\'esigne le degr\'e (\emph{resp}. le polyn\^ome de Hilbert-Samuel) de $\mathsf{V}$ relatif \`a ce faisceau. On consid\`ere \'egalement le multidegr\'e\begin{equation*}\mathscr{H}(\mathsf{V};X_{0},\ldots,X_{n}):=(\dim\mathsf{V})!\lim_{\alpha\to+\infty}\frac{H(\mathsf{V};\alpha .X_{0},\ldots,\alpha. X_{n})}{\alpha^{\dim\mathsf{V}}}\end{equation*}et on \'etend cette d\'efinition aux sous-sch\'emas int\`egres de $\mathbf{P}$ en prenant l'adh{\'e}rence de Zariski dans $\mathbf{P}$. Le plongement $G_{0}\times G\hookrightarrow\mathbf{P}$ permet alors de d\'efinir le polyn\^ome $\mathscr{H}(G';X_{0},\ldots,X_{n})$ pour tout sous-sch{\'e}ma en groupes $G'$ de $G_{0}\times G$. Rappelons que ses coefficients sont des entiers naturels de somme \'egale \` a $\deg G'$.
\subsection{Param{\`e}tres et choix d'un sous-groupe}
\label{parametresetchoixdunsousgroupe} Soit $x,\,\widetilde{D}_{0},\,\widetilde{D}_{1},\ldots,\widetilde{D}_{n},\,\widetilde{T},\,C_{0}$ des nombres r{\'e}els strictement positifs et $0<S_{0}\le S$ des entiers. Posons, pour chaque entier $i\in\{0,\ldots,n\}$, $\widetilde{D}_{i}^{\#}:=x\widetilde{D}_{i}$ et $D_{i}:=[\widetilde{D}_{i}^{\#}]$, ainsi que $T:=[\widetilde{T}]$. Nous supposerons que l'entier $T$ est non nul. En guise de support {\`a} l'intuition, mentionnons que $x$ est une variable \og d'ajustement\fg, les $D_{i}$ des degr{\'e}s de polyn{\^o}mes, $T$ un ordre de d{\'e}rivation, $S$ un nombre de points (tous multiples de $q$) et $C_{0}$ une constante positive (que l'on peut prendre enti{\`e}re) plus grande que toutes celles qui interviendront dans ce texte.\par
Lorsque $G'$ est un sous-groupe alg{\'e}brique de $G_{0}\times G$, on note $\lambda':=\codim_{W}W\cap t_{G'}$ et $r':=\codim_{G}G'$.
\begin{defi}
Soit $G'$ un sous-groupe alg{\'e}brique connexe de $G_{0}\times G$ tel que $t_{G'}+W\ne t_{G_{0}\times G}$. On d{\'e}finit \begin{equation*} A(G'):=\left(\frac{\widetilde{T}^{\lambda'}\card\left(\frac{\Sigma_{q}(S)+G'(\overline{k})}{G'(\overline{k})}\right)\mathscr{H}(G'\,;\,\widetilde{D}_{0},\ldots,\widetilde{D}_{n})}{C_{0}\mathscr{H}(G_{0}\times G\,;\,\widetilde{D}_{0},\ldots,\widetilde{D}_{n})}\right)^{\frac{1}{r'-\lambda'}}\end{equation*}et $B(G'):=A(G')^{\frac{r'-\lambda'}{r'}}\max{\{1,A(G')\}}^{\frac{\lambda'}{r'}}$. 
\end{defi} 
Remarquons alors que l'ensemble $$\{B(G')\,;\ B(G')\le B(\{0\})\ \text{et}\ t_{G'}+W\ne t_{G_{0}\times G}\}$$est fini. En effet, si $B(G')\le B(\{0\})$ alors $A(G')\le B(\{0\})$ donc $\mathscr{H}(G'\,;\,\widetilde{D}_{0},\ldots,\widetilde{D}_{n})$ et par cons{\'e}quent $\deg G'$ sont born{\'e}s. Le degr{\'e} de $G'$ {\'e}tant un entier, il n'y a qu'un nombre fini de valeurs possibles pour lui et comme les coefficients du polyn{\^o}me $\mathscr{H}(G'\,;\,X_{0},\ldots,X_{n})$ sont des entiers compris entre $0$ et $\deg G'$, il n'y en a {\'e}galement qu'un nombre fini. Il est alors clair que $B(G')$ ne prend qu'un nombre fini de valeurs lorsque $G'$ varie (parmi les sous-groupes tels que $B(G')\le B(\{0\})$). Cela justifie la d\'efinition suivante.
\begin{defi}\label{defiaseptdex}
On d{\'e}finit le nombre r{\'e}el strictement positif \begin{equation*}x:=\underset{t_{G'}+W\ne t_{G_{0}\times G}}{\mathrm{min}}\{B(G')\}\end{equation*}o\`u $G'$ varie parmi les sous-groupes alg{\'e}briques connexes de $G_{0}\times G$ tel que $t_{G'}+W$ est strictement inclus dans l'espace tangent $t_{G_{0}\times G}$. On note {\'e}galement $\widetilde{G}$ un sous-groupe parmi les $G'$ en question tel que $x=B(\widetilde{G})$.
\end{defi}
\begin{lemm}\label{lemme24deasept}
Supposons qu'il existe un sous-groupe alg{\'e}brique connexe $H_{1}\subseteq G_{0}\times G$ tel que $t_{H_{1}}+W\ne t_{G_{0}\times G}$ et $A(H_{1})\le 1$. Alors $x\le 1$ et pour tout sous-groupe alg{\'e}brique connexe $G'$ de $G_{0}\times G$ v{\'e}rifiant $t_{G'}+W\ne t_{G_{0}\times G}$, on a \begin{equation}\label{ineq:multiplicite} \widetilde{T}^{\lambda'}\card\left(\frac{\Sigma_{q}(S)+G'(\overline{k})}{G'(\overline{k})}\right)\mathscr{H}(G'\,;\,D^{\#}_{0},\ldots,D^{\#}_{n})\ge C_{0}\mathscr{H}(G_{0}\times G\,;\,D_{0}^{\#},\ldots,D_{n}^{\#})\ \cdotp\end{equation}De plus, cette in{\'e}galit{\'e} est une \'egalit\'e pour $G'=\widetilde{G}$.
\end{lemm}
 \begin{proof}
De l'existence de $H_{1}$, on d{\'e}duit imm{\'e}diatement que $B(H_{1})\le 1$ et donc $x\le 1$. Par ailleurs, pour un sch{\'e}ma en groupes $G'$ qui v{\'e}rifie les hypoth{\`e}ses de l'{\'e}nonc{\'e}, consid{\'e}rons le nombre r{\'e}el \begin{equation*}\mho_{G'}:=\frac{\widetilde{T}^{\lambda'}\card\left(\frac{\Sigma_{q}(S)+G'(\overline{k})}{G'(\overline{k})}\right)\mathscr{H}(G'\,;\,D^{\#}_{0},\ldots,D^{\#}_{n})}{C_{0}\mathscr{H}(G_{0}\times G\,;\,D_{0}^{\#},\ldots,D_{n}^{\#})}\ \cdotp\end{equation*}On a $x^{r'}\mho_{G'}=A(G')^{r'-\lambda'}$ par homog{\'e}n{\'e}it{\'e} de $\mathscr{H}$. \begin{itemize}
\item[$\bullet$] Si $A(G')\ge 1$ on a $\mho_{G'}\ge 1/x^{r'}\ge 1$.\item[$\bullet$] Si $A(G')\le 1$ on a $B(G')=A(G')^{\frac{r'-\lambda'}{r'}}\ge x$ donc encore $\mho_{G'}\ge 1$. \end{itemize}Cela d{\'e}montre l'in{\'e}galit{\'e}~\eqref{ineq:multiplicite}. De plus comme $A(\widetilde{G})\le 1$ (sinon $x=B(\widetilde{G})$ serait $>1$), on a $x=A(\widetilde{G})^{\frac{\widetilde{r}-\widetilde{\lambda}}{\widetilde{r}}}$ puis $\mho_{\widetilde{G}}=1$ et il y a donc bien {\'e}galit{\'e} dans~\eqref{ineq:multiplicite} pour $\widetilde{G}$.
\end{proof}
\subsection{Rang d'un syst{\`e}me d'{\'e}quations lin{\'e}aires}\label{subsec:rangsyslineaire}
Soit $(P_{0}=1,P_{1},\ldots,P_{D_{0}})$ une base de $k[X]_{\le D_{0}}$ et $\mathbf{D}=(D_{0},\ldots,D_{n})$. Consid{\'e}rons l'espace vectoriel $k[\mathbf{P}]_{\mathbf{D}}$ des polyn{\^o}mes multihomog{\`e}nes en les variables $\mathbf{X}_{i}=(X_{0}^{(i)},\ldots,X_{N_{i}}^{(i)})$, $0\le i\le n$ (en posant $N_{0}:=1$), de multidegr{\'e}s $\mathbf{D}$. Lorsque $v$ est une place de $k$ et \begin{equation*} P=\sum_{\vert\lambda_{i}\vert=D_{i}}{q_{\boldsymbol{\lambda}}\prod_{i=0}^{n}{\mathbf{X}_{i}^{\lambda_{i}}}}\in k[\mathbf{P}]_{\mathbf{D}},\end{equation*}($\boldsymbol{\lambda}=(\lambda_{0},\ldots,\lambda_{n})\in\prod_{i=0}^{n}{\mathbf{N}^{N_{i}+1}}$, $p_{\boldsymbol{\lambda}}\in k$), nous noterons $F_{P,v}$, ou plus simplement $F$ s'il n'y a pas d'ambiguït{\'e}, l'application \begin{equation} \label{defideFPv} F_{P,v}:z\mapsto\sum_{\boldsymbol{\lambda}}{q_{\boldsymbol{\lambda}}z_{0}^{\lambda_{0}}\prod_{i=1}^{n}{\Theta_{v,i}^{\lambda_{i}}(z_{i})}}\end{equation}d{\'e}finie pour $z=(z_{0},z_{1},\ldots,z_{n})\in\mathbf{C}_{v}\times\mathscr{T}_{v}$ et {\`a} valeurs dans $\mathbf{C}_{v}$. Dans cette expression, nous avons identifi{\'e} l'{\'e}l{\'e}ment $\lambda_{0}\in\mathbf{N}^{2}$ de longueur $D_{0}$ avec sa projection sur $\{0\}\times\mathbf{N}$, not{\'e}e aussi $\lambda_{0}$ (nous commettrons souvent cet abus de notation). Nous noterons $(p_{\boldsymbol{\lambda}})$ les coefficients de $P$ dans la base induite par $P_{\lambda_{0}}$, ce qui se traduit pour $F_{P,v}$ par \begin{equation*}F_{P,v}(z)=\sum_{\boldsymbol{\lambda}}{p_{\boldsymbol{\lambda}}P_{\lambda_{0}}(z_{0})\prod_{i=1}^{n}{\Theta_{v,i}^{\lambda_{i}}(z_{i})}}\ .\end{equation*}Soit $E$ la composante de degr{\'e} $\mathbf{D}$ de l'espace vectoriel $k[\mathbf{P}]$ quotient{\'e} par l'id\'eal des polyn\^omes identiquement nuls sur $G_{0}\times G$ (c'est-\`a-dire tels que $F_{P}$ soit identiquement nul). Soit $\mathbf{w}=(w_{0},\ldots,w_{g-t})$ une base de $W$ et $\cD_{w_{i}}$ l'op{\'e}rateur diff{\'e}rentiel associ{\'e} {\`a} $w_{i}$. Soit $S_{1},T_{1}$ des entiers naturels non nuls. Consid{\'e}rons le syst{\`e}me \begin{equation}\label{eqasept:huit}\forall\,(s,\boldsymbol{\tau})\in\mathbf{N}\times\mathbf{N}^{\dim W}\ ;\ 0\le s\le S_{1},\ \vert\boldsymbol{\tau}\vert\le T_{1},\quad \cD_{\mathbf{w}}^{\boldsymbol{\tau}}F_{P,v_{0}}(s,su)=0\end{equation}en les variables $q_{\boldsymbol{\lambda}}$ ou, de mani{\`e}re {\'e}quivalente, en les variables $p_{\boldsymbol{\lambda}}$. En consid{\'e}rant un sous-groupe alg{\'e}brique $G'$ de $G_{0}\times G$, nous allons majorer le rang $\rho$ du syst{\`e}me~\eqref{eqasept:huit} en fonction de $G'$. Pour cela, nous adoptons la m{\^e}me d{\'e}marche que celle du lemme~$6.7$ de~\cite{pph-miw} (voir aussi la preuve du lemme~$6.1$ de~\cite{sinnou2}). Soit $\mathbf{w}'$ une base d'un suppl{\'e}mentaire de $W\cap t_{G'}$ dans $W$ et $\lambda'=\dim W-\dim(W\cap t_{G'})$. Avec ces donn{\'e}es, on montre que $\rho$ est inf{\'e}rieur au rang du syst{\`e}me \begin{equation*}\forall\,(s,\boldsymbol{\tau})\in\mathbf{N}\times\mathbf{N}^{\lambda'}\ ;\ 0\le s\le S_{1},\ \vert\boldsymbol{\tau}\vert\le T_{1},\quad \cD_{\mathbf{w}'}^{\boldsymbol{\tau}}P(sq+G')=0,\end{equation*}et donc aussi plus petit que \begin{equation*}\card\{\boldsymbol{\tau}\in\mathbf{N}^{\lambda'};\ \vert\boldsymbol{\tau}\vert\le T_{1}\}\times\card\left(\frac{\Sigma_{q}(S_{1})+G'(\overline{k})}{G'(\overline{k})}\right)\dim\left(\mathbf{C}_{v_{0}}[\mathbf{P}]/I(G')\right)_{2\mathbf{D}}\ .\end{equation*}En vertu d'un th{\'e}or{\`e}me de Nesterenko~\cite{nes1985}, le dernier terme est lui-m{\^e}me contr{\^o}l{\'e} en fonction de $\mathscr{H}(G')$~:\begin{equation*}\dim\left(\mathbf{C}_{v_{0}}[\mathbf{P}]/I(G')\right)_{2\mathbf{D}}\le c_{\cst}\mathscr{H}(G';D_{0}',\ldots,D_{n}')\end{equation*}o\`u $D_{i}':=\max{\{1,D_{i}\}}$ et $c_{\theconstante}$ est une constante (explicite). Nous obtenons ainsi l'existence d'une constante $c_{\cst}$ telle que \begin{equation}\label{eqasept:dix}\rho\le c_{\theconstante}T_{1}^{\lambda'}\card\left(\frac{\Sigma_{q}(S_{1})+G'(\overline{k})}{G'(\overline{k})}\right)\mathscr{H}(G';D_{0}',\ldots,D_{n}')\ .\end{equation}En appliquant ce r{\'e}sultat au sous-groupe $G'=\widetilde{G}$ introduit dans la d{\'e}finition~\ref{defiaseptdex}, nous avons la 
\begin{prop} Supposons que $sq\not\in\widetilde{G}(\overline{k})$ pour tout $s\in\{1,\ldots,S\}$. Alors le rang $\rho$ du syst{\`e}me~\eqref{eqasept:huit} avec $S_{1}:=S_{0}$ et $T_{1}:=2(g+1)T$ v{\'e}rifie $$\rho\le C_{0}^{3/2}\frac{S_{0}}{S}\mathscr{H}(G;D_{0}',\ldots,D_{n}')\ .$$
\end{prop}
\begin{proof}
En effet la derni{\`e}re assertion du lemme~\ref{lemme24deasept} permet de simplifier la majoration~\eqref{eqasept:dix} en \begin{equation*}\rho\le c_{\theconstante}C_{0}\frac{\card\left(\frac{\Sigma_{q}(S_{0})+\widetilde{G}(\overline{k})}{\widetilde{G}(\overline{k})}\right)}{\card\left(\frac{\Sigma_{q}(S)+\widetilde{G}(\overline{k})}{\widetilde{G}(\overline{k})}\right)}\mathscr{H}(G;D_{0}',\ldots,D_{n}')\end{equation*}car chacune des applications partielles $$x_{i}\mapsto\mathscr{H}(\widetilde{G};x_{0},\ldots,x_{n})/\mathscr{H}(G;x_{0},\ldots,x_{n})$$est d{\'e}croissante (voir propri{\'e}t{\'e}~$4.4$, p.~$414$, de~\cite{hirata2}) et $D_{i}^{\#}\le 2D_{i}'$. Pour conclure, nous utilisons l'hypoth{\`e}se faite sur $q$.
\end{proof}
\begin{rema}
L'hypoth{\`e}se sur $q$ correspond {\`a} celle faite sur $p$ dans les th{\'e}or{\`e}mes~\ref{theounasept} et~\ref{theodeuxasept} et elle n'intervient dans toute la preuve qu'{\`a} cet endroit pr{\'e}cis.
\end{rema}

\subsection{Remarque auxiliaire (non-nullit{\'e} du $\max{\{D_{i}\}}$)}
\begin{lemm}\label{lemmenonzeroasept}
Supposons que $x\le 1$ et que $sq\not\in\widetilde{G}(\overline{k})$ pour tout $s\in\{1,\ldots,C_{0}\deg(G_{0}\times G)\}$. Si les deux conditions suivantes sont satisfaites \begin{enumerate}
\item $S\ge C_{0}\deg(G_{0}\times G)$,\item $\widetilde{T}(S+1)\underset{1\le i\le n}{\mathrm{max}}{\{\widetilde{D}_{i}\}}\ge C_{0}\deg(G_{0}\times G)\widetilde{D}_{0}$,\end{enumerate}alors les entiers $D_{i}$, $1\le i\le n$, ne sont pas tous nuls.
\end{lemm}
\begin{proof}
Par construction, il s'agit de montrer que $\underset{1\le i\le n}{\mathrm{max}}\{x\widetilde{D}_{i}\}\ge 1$. Comme $x=B(\widetilde{G})=A(\widetilde{G})^{\frac{\widetilde{r}-\widetilde{\lambda}}{\widetilde{r}}}$ (car pr{\'e}cis{\'e}ment $x\le 1$ par hypoth{\`e}se), on a \begin{equation}\label{forsansnomasept}\left(x\underset{1\le i\le n}{\mathrm{max}}\{\widetilde{D}_{i}\}\right)^{\widetilde{r}}=\frac{\widetilde{T}^{\widetilde{\lambda}}\card\left(\frac{\Sigma_{q}(S)+\widetilde{G}(\overline{k})}{\widetilde{G}(\overline{k})}\right)\mathscr{H}(\widetilde{G}\,;\,\widetilde{D}_{0},\ldots,\widetilde{D}_{n})\underset{1\le i\le n}{\mathrm{max}}{\{\widetilde{D}_{i}\}}^{\widetilde{r}}}{C_{0}\mathscr{H}(G_{0}\times G\,;\,\widetilde{D}_{0},\ldots,\widetilde{D}_{n})}\ \cdotp\end{equation}Soit $\pi_{0}:G_{0}\times G\to G_{0}$ la projection canonique sur $G_{0}$ et $\pi:G_{0}\times G\to G$ celle sur $G$. Si $\pi_{0}(\widetilde{G})=\{0\}$ alors $\widetilde{G}=\{0\}\times\pi(\widetilde{G})$ donc \begin{enumerate}\item[$\bullet$]$\card\left(\frac{\Sigma_{q}(S)+\widetilde{G}(\overline{k})}{\widetilde{G}(\overline{k})}\right)=S+1$ (car la premi{\`e}re composante de $q$ est $1$),\item[$\bullet$] $\widetilde{\lambda}\ge 1$ (car sinon $W\subseteq t_{\widetilde{G}}$ et on aurait $t_{G_{0}}\subseteq t_{\pi_{0}(\widetilde{G})}=\{0\}$), \item[$\bullet$] $\mathscr{H}(\widetilde{G}\,;\,\widetilde{D}_{0},\ldots,\widetilde{D}_{n})\ge\underset{\genfrac{}{}{0pt}{}{i_{1}+\cdots+i_{n}=\dim\widetilde{G}}{0\le i_{j}\le g_{j}}}{\mathrm{min}}\{\widetilde{D}_{1}^{i_{1}}\cdots\widetilde{D}_{n}^{i_{n}}\}$\end{enumerate}et de la formule~\eqref{forsansnomasept} on d{\'e}duit \begin{equation}\left(x\underset{1\le i\le n}{\mathrm{max}}\{\widetilde{D}_{i}\}\right)^{\widetilde{r}}\ge\frac{\widetilde{T}(S+1)\underset{1\le i\le n}{\mathrm{max}}\{\widetilde{D}_{i}\}}{C_{0}\deg(G_{0}\times G)\widetilde{D}_{0}}\end{equation}et cette derni{\`e}re quantit{\'e} est sup{\'e}rieure {\`a} $1$ par hypoth{\`e}se. \par Si $\pi_{0}(\widetilde{G})=G_{0}$ on a alors \begin{equation}\mathscr{H}(\widetilde{G}\,;\,\widetilde{D}_{0},\ldots,\widetilde{D}_{n})\ge\widetilde{D}_{0}\underset{\genfrac{}{}{0pt}{}{i_{1}+\cdots+i_{n}=\dim\widetilde{G}-1}{0\le i_{j}\le g_{j}}}{\mathrm{min}}\{\widetilde{D}_{1}^{i_{1}}\cdots\widetilde{D}_{n}^{i_{n}}\}\end{equation}et\begin{equation}\left(x\underset{1\le i\le n}{\mathrm{max}}\{\widetilde{D}_{i}\}\right)^{\widetilde{r}}\ge\frac{\widetilde{T}^{\widetilde{\lambda}}\card\left(\frac{\Sigma_{q}(S)+\widetilde{G}(\overline{k})}{\widetilde{G}(\overline{k})}\right)}{C_{0}\deg(G_{0}\times G)}\ge 1\end{equation}par hypoth{\`e}se.
\end{proof}
\subsection{Lemme de multiplicit{\'e}s}
\label{asept:lemmedemuliplicite}
Rappelons que si $\mathbf{w}=(w_{0},\ldots,w_{g-t})$ d{\'e}signe une base de $W$, on note $\cD_{w_{i}}$ l'op{\'e}rateur diff{\'e}rentiel associ{\'e} {\`a} $w_{i}$. L'objet de ce paragraphe est de montrer que, sous des hypoth{\`e}ses \og minimales\fg\ (en particulier sans la n{\'e}cessit{\'e} d'un choix tr{\`e}s pr{\'e}cis des param{\`e}tres {\`a} cette {\'e}tape) et gr{\^a}ce au lemme de multiplicit{\'e}s de Philippon~\cite{philippon2}, il est possible d'affirmer que le syst{\`e}me \begin{equation}\forall\,(m,\mathbf{t})\in\mathbf{N}\times\mathbf{N}^{\dim W},\ 0\le m\le (g+1)S,\ \vert\mathbf{t}\vert \le(g+1)T,\ \cD_{\mathbf{w}}^{\mathbf{t}}F_{P,v_{0}}(m,mu)=0\end{equation}n'admet pas de solution polynomiale $P$ non nulle (les notations sont celles du \S~\ref{subsec:rangsyslineaire}).
\begin{lemm}\label{lemme:multiplicite}Supposons que les entiers $D_{1},\ldots,D_{n}$ ne sont pas tous nuls et que $x\le 1$. Supposons {\'e}galement que, pour tout sous-groupe alg{\'e}brique $G'$ de $G$ tel que $t_{G'}+V\ne t_{G}$, pour tout entier $s\in\{1,\ldots,C_{0}\}$, on ait $sp\not\in G'(\overline{k})$ et supposons enfin que\begin{equation}\label{hypounasept}\widetilde{T}\ge c_{\cst}\max{\left\{\frac{\widetilde{D}_{0}}{(S+1)^{1-y}},\widetilde{D}_{1},\ldots,\widetilde{D}_{n},1\right\}}\end{equation}avec $c_{\theconstante}=8^{g}(\deg G_{0}\times G)\prod_{i=1}^{n}{\deg G_{i}}$. Alors il n'existe pas de polyn{\^o}me $P\in E_{v_{0}}\setminus\{0\}$ tel que la d{\'e}riv{\'e}e $\cD_{\mathbf{w}}^{\mathbf{t}}F_{P,v_{0}}(m,mu)$ soit nulle pour tout $(m,\mathbf{t})\in\mathbf{N}\times\mathbf{N}^{\dim W},\ 0\le m\le (g+1)S$ et $\vert\mathbf{t}\vert \le(g+1)T$.
\end{lemm}
\begin{proof}
Supposons qu'un tel polyn{\^o}me $P\ne 0$ existe. Alors, d'apr{\`e}s le lemme de multiplicit{\'e}s~\cite{philippon2}, il existe un sous-groupe alg{\'e}brique connexe et strict $G'$ de $G_{0}\times G$ tel que \begin{equation}\begin{split}\label{ineqlemmezerosasept}&T^{\lambda'}\card\left(\frac{\Sigma_{q}(S)+G'(\overline{k})}{G'(\overline{k})}\right)\mathscr{H}(G'\,;\,D_{0}',\ldots,D_{n}')\\ &\le 2^{g}\mathscr{H}(G_{0}\times G\,;\,D_{0}',\ldots,D_{n}')\end{split}\end{equation}o\`u, rappelons-le, $\lambda'=\codim_{W}W\cap t_{G'}$ et $D_{i}'=\max{\{1,D_{i}\}}$. Nous allons examiner s{\'e}par{\'e}ment les cas $t_{G'}+W=t_{G_{0}\times G}$ et $t_{G'}+W\ne t_{G_{0}\times G}$ afin de conclure que $G'$ ne peut pas exister. \par Si $t_{G'}+W=t_{G_{0}\times G}$ alors $\lambda'=r'=\codim_{G_{0}\times G}G'$ et l'in{\'e}galit{\'e}~\eqref{ineqlemmezerosasept} entra{\^\i}ne \begin{equation*} T^{r'}\le 2^{g}(\deg(G_{0}\times G))\max{\{D_{0},\ldots,D_{n}\}}^{r'},\end{equation*}ce qui contredit l'hypoth{\`e}se~\eqref{hypounasept} lorsque $y$ vaut $1$ (cas g\'en\'eral). Dans le cas semi-ab{\'e}lien ($y=0$), le groupe $G'$ s'{\'e}crit $G_{0}'\times A$ avec $G_{0}'\subseteq G_{0}$ et $A\subseteq G$. Si $G_{0}'=G_{0}$ alors l'in{\'e}galit{\'e}~\eqref{ineqlemmezerosasept} implique \begin{equation*} T^{r'}\mathscr{H}(A;D_{1}',\ldots,D_{n}')\le 2^{g}(\deg(G_{0}\times G))(D_{1}')^{g_{1}}\cdots(D_{n}')^{g_{n}}\end{equation*}ce qui est incompatible avec $\widetilde{T}>4^{g}(\deg(G_{0}\times G))\max{\{\widetilde{D}_{1},\ldots,\widetilde{D}_{n}\}}$. Si $G_{0}'=\{0\}$ alors \begin{equation*}\card\left(\frac{\Sigma_{q}(S)+G'(\overline{k})}{G'(\overline{k})}\right)=S+1\end{equation*}et~\eqref{ineqlemmezerosasept} devient \begin{equation*} T^{r'}(S+1)\mathscr{H}(A;D_{1}',\ldots,D_{n}')\le 2^{g}(\deg(G_{0}\times G))D_{0}'(D_{1}')^{g_{1}}\cdots(D_{n}')^{g_{n}}\end{equation*}ce qui implique $T^{r'}(S+1)\le 2^{g}(\deg(G_{0}\times G))D_{0}'\max{\{D_{1},\ldots,D_{n}\}}^{r'-1}$, ce qui est encore impossible d'apr{\`e}s l'hypoth{\`e}se~\eqref{hypounasept}. Ce cas ne peut donc pas se produire et on a n{\'e}cessairement $t_{G'}+W\ne t_{G_{0}\times G}$. Soit alors $1\le\kappa_{1}<\cdots<\kappa_{h}\le n$ les entiers pour lesquels $D_{\kappa_{i}}\ne 0$, $1\le i\le h$ et notons $\pi_{\kappa}:G_{0}\times G\to\prod_{i=1}^{h}{G_{\kappa_{i}}}$ la projection canonique. \par Si $D_{0}'=1$, on a \begin{equation}\label{ineqimpo}\binom{\dim G'}{\dim\pi_{\kappa}(G')}\mathscr{H}(\pi_{\kappa}(G');D_{\kappa_{1}},\ldots,D_{\kappa_{h}})\le\mathscr{H}(G';D_{0}',\ldots,D_{n}')\end{equation}ce qui entra{\^\i}ne \textit{via}~\eqref{ineqlemmezerosasept} \begin{equation*} T^{\lambda'}\card\left(\frac{\Sigma_{q}(S)+G'(\overline{k})}{G'(\overline{k})}\right)\le c_{\cst}\frac{\mathscr{H}(\pi_{\kappa}(G);D_{\kappa_{1}},\ldots,D_{\kappa_{h}})}{\mathscr{H}(\pi_{\kappa}(G');D_{\kappa_{1}},\ldots,D_{\kappa_{h}})}\end{equation*}\newcounter{predeux}\setcounter{predeux}{\value{constante}}avec \begin{equation*} c_{\thepredeux}:=\frac{2^{g}(g+1)!(\dim\pi_{\kappa}(G'))!(\dim G'-\dim\pi_{\kappa}(G'))!}{(\dim\pi_{\kappa}(G))!(\dim G')!}\prod_{\genfrac{}{}{0pt}{}{m\ne\kappa_{i}}{1\le i\le h}}{\frac{\deg G_{m}}{g_{m}!}}\ \cdotp\end{equation*}Comme chacune des applications partielles \begin{equation*} x_{i}\mapsto\frac{\mathscr{H}(\pi_{\kappa}(G);D_{\kappa_{1}},\ldots,x_{i},\ldots,D_{\kappa_{h}})}{\mathscr{H}(\pi_{\kappa}(G');D_{\kappa_{1}},\ldots,x_{i},\ldots,D_{\kappa_{h}})},\quad 1\le i\le h,\end{equation*}est croissante sur $]0,+\infty[$, on d{\'e}duit de la majoration $D_{i}\le x\widetilde{D}_{i}$ l'in{\'e}galit{\'e}
\begin{equation}\begin{split}\label{ineqzeroasept}&T^{\lambda'}\card\left(\frac{\Sigma_{q}(S)+G'(\overline{k})}{G'(\overline{k})}\right)\\ &\qquad\le c_{\thepredeux}\frac{\mathscr{H}(\pi_{\kappa}(G)\,;\,\widetilde{D}_{\kappa_{1}},\ldots,\widetilde{D}_{\kappa_{h}})}{\mathscr{H}(\pi_{\kappa}(G')\,;\,\widetilde{D}_{\kappa_{1}},\ldots,\widetilde{D}_{\kappa_{h}})}\,x^{\dim\pi_{\kappa}(G)-\dim\pi_{\kappa}(G')}\ .\end{split}\end{equation}Posons \begin{equation*} G'':=\pi_{\kappa}(G')\times\prod_{m\not\in\{\kappa_{1},\ldots,\kappa_{h}\}}{G_{m}}\end{equation*}(vu, apr{\`e}s permutation {\'e}ventuelle des facteurs, comme un sous-groupe de $G_{0}\times G$). Ce groupe alg{\'e}brique est diff{\'e}rent de $G$ sinon $\pi_{\kappa}(G')=\pi_{\kappa}(G)$ et l'in{\'e}galit{\'e}~\eqref{ineqzeroasept} entra{\^\i ne} $T^{\lambda'}\card\left(\frac{\Sigma_{q}(S)+G'(\overline{k})}{G'(\overline{k})}\right)\le c_{\thepredeux}$. L'hypoth{\`e}se sur le point $p$ (appliqu{\'e}e {\`a} la projection de $G'$ sur $G$) entra{\^\i}ne alors $C_{0}\le c_{\thepredeux}$ ce qui est impossible si $C_{0}$ est assez grand. Nous allons maintenant obtenir {\`a} partir de~\eqref{ineqzeroasept} une in{\'e}galit{\'e} pour $G''$ analogue {\`a}~\eqref{ineqlemmezerosasept}. En effet, observons d'une part que \begin{equation*}\frac{\mathscr{H}(G_{0}\times G\,;\,\widetilde{D}_{0},\ldots,\widetilde{D}_{n})}{\mathscr{H}(G''\,;\,\widetilde{D}_{0},\ldots,\widetilde{D}_{n})}=c_{\cst}\frac{\mathscr{H}(\pi_{\kappa}(G)\,;\,\widetilde{D}_{\kappa_{1}},\ldots,\widetilde{D}_{\kappa_{h}})}{\mathscr{H}(\pi_{\kappa}(G')\,;\,\widetilde{D}_{\kappa_{1}},\ldots,\widetilde{D}_{\kappa_{h}})}\end{equation*}\newcounter{pretrois}\setcounter{pretrois}{\value{constante}}avec $c_{\thepretrois}:=\frac{(g+1)!\left(\dim\pi_{\kappa}(G')\right)!}{(\dim G'')!\left(\dim\pi_{\kappa}(G)\right)!}$ et, d'autre part, on a\begin{enumerate}\item $\dim\pi_{\kappa}(G)-\dim\pi_{\kappa}(G')=\codim_{G_{0}\times G}G''=:r''$,\item $\lambda'':=\codim_{W}W\cap t_{G''}\le\codim_{W}W\cap t_{G'}=\lambda'$ (car $G'\subseteq G''$), \item $\card\left(\frac{\Sigma_{q}(S)+G''(\overline{k})}{G''(\overline{k})}\right)\le\card\left(\frac{\Sigma_{q}(S)+G'(\overline{k})}{G'(\overline{k})}\right)$.\end{enumerate}L'in{\'e}galit{\'e}~\eqref{ineqzeroasept} devient\begin{equation}\label{ineqaseptdeux}\frac{T^{\lambda''}\card\left(\frac{\Sigma_{q}(S)+G''(\overline{k})}{G''(\overline{k})}\right)\mathscr{H}(G''\,;\,\widetilde{D}_{0},\ldots,\widetilde{D}_{n})}{\mathscr{H}(G_{0}\times G\,;\,\widetilde{D}_{0},\ldots,\widetilde{D}_{n})}\le \frac{c_{\thepredeux}}{c_{\thepretrois}}x^{r''}\ \cdotp\end{equation}En reprenant alors exactement les m{\^e}mes arguments que dans la premi{\`e}re partie de la preuve et en rempla{\c c}ant la constante $2^{g}$ par $c_{\thepredeux}/c_{\thepretrois}$ (il faut observer que ce quotient est strictement plus petit que $4^{g}\prod_{i=1}^{n}{\deg G_{i}}$), on d{\'e}montre que $t_{G''}+W\ne t_{G_{0}\times G}$. L'in{\'e}galit{\'e}~\eqref{ineqaseptdeux} se lit alors en fonction de $A(G'')$ (en minorant $T$ par $\widetilde{T}/2$)~:\begin{equation}\label{estifinasept} A(G'')^{r''-\lambda''}\le\left(\frac{2^{\lambda''}c_{\thepredeux}}{c_{\thepretrois}C_{0}}\right)x^{r''}\ \cdotp\end{equation}Comme $C_{0}>2^{g}c_{\thepredeux}/c_{\thepretrois}$ et $x\le 1$, cette in{\'e}galit{\'e} implique $A(G'')\le 1$ donc $x\le B(G'')=A(G'')^{\frac{r''-\lambda''}{r''}}$, ce qui contredit~\eqref{estifinasept}. On vient donc de montrer que n{\'e}cessairement $D_{0}'=D_{0}$, c.{\`a}-d. $D_{0}\ge 1$. \`A quelques variantes pr{\`e}s, la m{\^e}me preuve conduit encore {\`a} une contradiction. En effet, consid{\'e}rons $\pi_{0,\kappa}$ la projection $G_{0}\times G\to G_{0}\times\pi_{\kappa}(G)$ et posons \begin{equation} G^{\star}=\pi_{0,\kappa}(G')\times\prod_{m\not\in\{\kappa_{1},\ldots,\kappa_{h}\}}{G_{m}},\end{equation}vu comme sous-groupe de $G_{0}\times G$. Comme \begin{equation}\binom{\dim G'}{\dim\pi_{0,\kappa}(G')}\mathscr{H}(\pi_{0,\kappa}(G')\,;\,D_{0},D_{\kappa_{1}},\ldots,D_{\kappa_{h}})\le\mathscr{H}(G'\,;\,D_{0}',\ldots,D_{n}'),\end{equation}un raisonnement similaire au pr{\'e}c{\'e}dent ({\`a} partir de l'in{\'e}galit{\'e}~\eqref{ineqimpo}) conduit encore {\`a} une impossibilit{\'e} {\`a} condition de remplacer $\pi_{\kappa}(G')$ par $\pi_{0,\kappa}(G')$ et $G''$ par $G^{\star}$ (dans les constantes $c_{\thepredeux}$ et $c_{\thepretrois}$ en particulier).\par \emph{Conclusion}~: Dans tous les cas, l'existence de $G'$ aboutit {\`a} une contradiction, ce qui d{\'e}montre ainsi le lemme~\ref{lemme:multiplicite}.
\end{proof}

\subsection{Poids de la droite affine}
\label{subsec:poidsdeladroiteaffine}
Rappelons que $(P_{\lambda_{0}})_{\lambda_{0}\in\{0,\ldots,D_{0}\}}$ d\'esigne une famille libre de polyn{\^o}mes en une variable (voir \S~\ref{subsec:rangsyslineaire}). Dans ce paragraphe, nous introduisons une quantit{\'e} mi-arithm{\'e}tique mi-analytique qui mesure l'influence du choix de cette famille sur les param\`etres $U_{0}$ et $U_{1}$ des th{\'e}or{\`e}mes~\ref{theounasept} et~\ref{theodeuxasept}. C'est ce que nous voulons {\'e}voquer par la terminologie \og poids de $\mathbb{G}_{\mathrm{a}}$\fg.
\begin{defi}Si $v_{0}$ est archim{\'e}dienne, nous appelons \emph{poids} de $\mathbb{G}_{\mathrm{a}}$ relatif {\`a} la famille $(P_{\lambda_{0}})$, aux param{\`e}tres $(T,S,\mathfrak{e})$ et {\`a} la place $v_{0}$, la quantit{\'e} \begin{equation*}\aleph((P_{\lambda_{0}})):=h\left(\left\{\frac{1}{t_{0}!}P_{\lambda_{0}}^{(t_{0})}(s),\ \genfrac{}{}{0pt}{}{\genfrac{}{}{0pt}{}{}{0\le\lambda_{0}\le D_{0}}}{\genfrac{}{}{0pt}{}{0\le s\le S}{0\le t_{0}\le T}}\right\}\right)+\frac{1}{D}\log\max_{\genfrac{}{}{0pt}{}{0\le t_{0}\le T}{\vert z\vert\le\mathfrak{e} S}}\left\vert\frac{1}{t_{0}!}P_{\lambda_{0}}^{(t_{0})}(z)\right\vert_{v_{0}}\ \cdotp\end{equation*}Lorsque $v_{0}$ est ultram{\'e}trique, le poids de $\mathbb{G}_{\mathrm{a}}$ (relatif {\`a} $(T,S,\mathfrak{r})$) est la quantit{\'e} obtenue en rempla{\c c}ant ci-dessus $\mathfrak{e} S$ (qui est en indice dans le dernier terme) par $\mathfrak{r}$.
\end{defi}Le poids de $\mathbb{G}_{\mathrm{a}}$ est le terme r{\'e}siduel qui provient de l'introduction m{\^e}me du groupe $\mathbb{G}_{\mathrm{a}}$ dans la d\'emonstration des th{\'e}or{\`e}mes~\ref{theounasept} et~\ref{theodeuxasept}. Cet inconv{\'e}nient s'av{\`e}re largement compens{\'e} par au moins deux avantages que procure $\mathbb{G}_{\mathrm{a}}$. D'une part, il s'accompagne d'un param{\`e}tre $D_{0}$ qui facilite la construction du polyn{\^o}me auxiliaire en rendant la condition de Siegel plus simple {\`a} satisfaire. D'autre part, il permet de modifier le point $p$ en un point $q=(1,p)$ moins vuln{\'e}rable aux ph{\'e}nom{\`e}nes de torsion modulo des sous-groupes particuliers de $G_{0}\times G$ ({\`a} commencer par le sous-groupe nul, $q$ n'{\'e}tant alors jamais de torsion). Cela est particuli{\`e}rement important pour le lemme de multiplicit{\'e} qui fait intervenir le cardinal du quotient $(\Sigma_{q}(S)+G'(\overline{k}))/G'(\overline{k})$, $G'$ {\'e}tant un sous-groupe alg{\'e}brique de $\mathbb{G}_{\mathrm{a}}\times G$ (voir {\`a} cet {\'e}gard le \S~\ref{asept:lemmedemuliplicite}).\par Pour minimiser ce poids, nous allons utiliser la famille des polyn{\^o}mes de Matveev, d\'efinie de la mani\`ere suivante. Soit $(\Delta_{n})_{n\in\mathbf{N}}$ la famille des polyn{\^o}mes binomiaux d\'efinie par~\eqref{binomiaux}.
\begin{defi}\'Etant donn{\'e} $\lambda_{0},D_{0}^{\flat}\in\mathbf{N}$, le polyn{\^o}me de Matveev $\delta_{D_{0}^{\flat}}(X;\lambda_{0})$ est $\Delta_{D_{0}^{\flat}}(X)^{q}\Delta_{r}(X)$ o\`u les entiers $q$ et $r$ sont respectivement les quotient et reste de la division euclidienne de $\lambda_{0}$ par $D_{0}^{\flat}$.
\end{defi} Le degr{\'e} de $\delta_{D_{0}^{\flat}}(X;\lambda_{0})$ est $\lambda_{0}$. Par cons{\'e}quent, lorsque $D_{0}^{\flat}$ est fix{\'e}, la famille $\delta_{D_{0}^{\flat}}(X;\lambda_{0})$, $\lambda_{0}=0,\ldots,D_{0}$, forme une base de $k[X]_{\le D_{0}}$.
\begin{lemm}\label{lemme:estimationdespoids}
Il existe une constante absolue $c_{\cst}\ge 1$\newcounter{aseptarchicst}\setcounter{aseptarchicst}{\value{constante}} pour laquelle nous disposons des estimations suivantes.
\begin{enumerate}
\item[$\bullet$] Si $v_{0}$ est archim{\'e}dienne alors \begin{equation}\begin{split}\label{ineq:estimationdupoidsdesfeldman}&\aleph\left(\left(\delta_{D_{0}^{\flat}}(X;\lambda_{0})\right)_{\lambda_{0}}\right)\\ &\quad\le c_{\theconstante}\left(D_{0}\log\left(e+\frac{S}{D_{0}^{\flat}}\right)+\min{\left(D_{0},T\right)}D_{0}^{\flat}+\frac{D_{0}}{D}\log\left(1+\frac{\mathfrak{e} S}{D_{0}^{\flat}}\right)\right)\ \cdotp
\end{split}\end{equation}
\item[$\bullet$] Si $v_{0}$ est ultram{\'e}trique alors l'in{\'e}galit{\'e} ci-dessus reste vraie en rempla{\c c}ant le dernier terme par $\frac{D_{0}}{D}\log\left(\mathfrak{r}\right)$.
\end{enumerate}
\end{lemm}
Dans le cas archim{\'e}dien la d{\'e}monstration de cette majoration d{\'e}coule des estimations sur les polyn{\^o}mes de Matveev qui sont donn{\'e}es dans le livre de Waldschmidt~\cite{miw4}, p.~$269$ et suivantes. Dans le cas ultram{\'e}trique, l'{\'e}valuation de la d{\'e}riv{\'e}e $\delta_{D_{0}^{\flat}}(X;\lambda_{0},t_{0})$ pour $t_{0}\le T$ et $\vert z\vert\le\mathfrak{r}$ repose sur la formule de Leibniz et l'existence d'une constante absolue $c_{\cst}>0$ telle que $\vert\lambda_{0}!\vert_{v_{0}}\ge c_{\theconstante}^{\lambda_{0}}$. Notons par ailleurs qu'en choisissant $D_{0}^{\flat}=1$ l'on retrouve une estimation du poids de la famille $(X^{\lambda_{0}})_{0\le\lambda_{0}\le D_{0}}$:\begin{equation*}\aleph((X^{\lambda_{0}}))\le c_{\theconstante}D_{0}\left(\log S+\frac{\log\mathfrak{e}}{D}\right)\end{equation*}(et $\mathfrak{r}$ {\`a} la place de $\mathfrak{e}$ dans le cas ultram{\'e}trique). 

\begin{rema}Bien que tous les th{\'e}or{\`e}mes {\'e}nonc{\'e}s reposent sur le m{\^e}me choix de la base $(P_{i})_{i}$ (base des polyn{\^o}mes de Matveev en l'occurrence), il nous semble pr{\'e}f{\'e}rable de conserver une base indiff{\'e}renci{\'e}e jusqu'{\`a} la toute fin de la d{\'e}monstration (\S~\ref{subsec:conclusionasept}). Outre une justification \textit{a posteriori} du choix des param{\`e}tres, cela permet {\'e}galement d'obtenir d'une part des variantes d'{\'e}nonc{\'e}s {\`a} moindre frais et d'autre part une meilleure compr{\'e}hension du r{\^o}le jou{\'e} par ce facteur $\mathbb{G}_{\mathrm{a}}$ suppl{\'e}mentaire au cours de la preuve, comme nous venons de le voir. Ce proc{\'e}d{\'e} a d{\'e}j{\`a} {\'e}t{\'e} mis en {\oe}uvre (sous une forme tr{\`e}s l{\'e}g{\`e}rement diff{\'e}rente) par Waldschmidt (voir~\cite{miw4}, pp.~$477-480$). 
\end{rema}

\subsection{Lemme de Siegel absolu}\label{subsec:lemmesiegelabsolu}Nous pr\'esentons un raffinement du lemme de Siegel absolu de Roy \& Thunder~\cite{roy-thunder}, signal\'e par David et Philippon dans~\cite{david-pph}, et qui repose sur une in\'egalit\'e de Zhang relative aux minima successifs d'une vari\'et\'e arithm\'etique.
\begin{lemm}\label{lemmedesiegelabsolu}
Soit $m$ un entier naturel $\ge 1$ et $\mathsf{V}$ un sous-espace vectoriel de $\overline{\mathbf{Q}}^{m+1}$, de dimension $d\ge 1$. Il existe une base $(v_{1},\ldots,v_{d})$ de $\mathsf{V}$ telle que \begin{equation} \sum_{i=1}^{d}{h_{\mathrm{L}^{2}}(v_{i})}\le h(\mathsf{V})+d\log(d)\end{equation}o\`u $h(\mathsf{V})$ est la hauteur (logarithmique absolue) de Schmidt de $\mathsf{V}$. 
\end{lemm}
Le r\'esultat pr\'ec\'edemment cit\'e de Roy \& Thunder conduit \`a une majoration de ce type mais avec un terme lin\'eaire en $d^{2}$ en lieu et place du $d\log(d)$. En fait l'argument de~\cite{david-pph}, remarque du \S~$4.2$, p.~523-524, fournit un r\'esultat un peu plus fort~: pour tout nombre r\'eel $\epsilon>0$, il existe une base $(v_{1},\ldots,v_{d})$ de $\mathsf{V}$ (qui d\'epend de cet $\epsilon$) telle que \begin{equation}\label{inequalitySiegel} \sum_{i=1}^{d}{h_{\mathrm{L}^{2}}(v_{i})}\le h(\mathsf{V})+\sum_{j=1}^{d-1}{\sum_{i=1}^{j}{\frac{1}{2i}}}+\epsilon\ \cdotp\end{equation}

\section{D{\'e}monstrations des th{\'e}or{\`e}mes~\ref{theounasept} et~\ref{theodeuxasept}}
\label{secasept:demonstrationdestheos}
Nous fixons une fois pour toutes un corps de nombres $K$, qui contient $k$, de sorte que tous les nombres alg{\'e}briques consid{\'e}r{\'e}s au cours de la d\'emonstration et qui, bien s{\^u}r, sont en nombre fini, sont inclus dans $K$. Il est commode d'introduire un tel corps pour les estimations locales de ces nombres alg{\'e}briques.

\subsubsection*{Description de la preuve}
La d{\'e}marche suivie est assez classique et elle est commune aux deux th{\'e}or{\`e}mes {\`a} d{\'e}montrer. Il s'agit de construire un {\'e}l{\'e}ment $\alpha$ de $K\setminus\{0\}$, de \og petite\fg\ hauteur, et dont toutes les valeurs absolues $v$-adiques, pour $v$ une place de $K$ au-dessus de $v_{0}$, sont major{\'e}es par un terme lin{\'e}aire en la distance $\mathrm{d}_{v_{0}}(u,V)$. De sorte que de la formule du produit appliqu{\'e}e {\`a} $\alpha$ se d{\'e}duit une minoration de cette distance, ce qui est l'assertion des th{\'e}or{\`e}mes~\ref{theounasept} et~\ref{theodeuxasept}. Comme cela est fr{\'e}quent en transcendance, l'{\'e}l{\'e}ment $\alpha$ en question provient d'un coefficient de Taylor d'une fonction de la forme $P\circ\exp_{(G_{0}\times G)(\mathbf{C}_{v_{0}})}$ restreinte {\`a} l'espace $W$. Dans cette expression, $P$ est un polyn{\^o}me construit au moyen du lemme de Siegel absolu {\'e}nonc{\'e} dans le paragraphe pr{\'e}c{\'e}dent et $\exp_{(G_{0}\times G)(\mathbf{C}_{v_{0}})}:=(\exp_{G_{0}(\mathbf{C}_{v_{0}})},\exp_{v_{0}})$ d{\'e}signe l'exponentielle ({\`a} valeurs dans l'espace multiprojectif $\mathbf{P}$) du groupe de Lie $(G_{0}\times G)(\mathbf{C}_{v_{0}})$.

\subsection{Choix des param{\`e}tres}
\label{choixunasept}
 Dans ce paragraphe, nous pr{\'e}cisons toutes ces donn{\'e}es pour la d{\'e}monstration proprement dite des th{\'e}or{\`e}mes~\ref{theounasept} et~\ref{theodeuxasept}. Cependant les choix que nous d{\'e}voilons ici ne seront v{\'e}ritablement utilis{\'e}s qu'{\`a} l'{\'e}tape \og Extrapolation\fg\ (\S~\ref{aspetsection:extrapolation}) \textit{via} les lemmes~\ref{lemmearchi:choixdesparametres} et~\ref{lemmeultra:choixdesparametres} qui vont suivre. Posons\begin{equation*}\Upsilon:=\left\{(s,\mathbf{t})\in\mathbf{N}\times\mathbf{N}^{\dim W}\,;\ 0\le s\le S_{0}\ \text{et}\ \vert\mathbf{t}\vert\le 2(g+1)T\right\}\cdotp\end{equation*}Dans un souci de clart{\'e}, nous distinguons les choix selon que $v_{0}$ est une place archim{\'e}dienne ou ultram{\'e}trique.
\subsubsection{$v_{0}$ archim{\'e}dienne}Rappelons que $\mathfrak{a}$ d{\'e}signe un entier sup\'erieur ou \'egal \`a $D\max{\{1,h(V)\}}/\log\mathfrak{e}$. Le param{\`e}tre $C_{0}$ est une constante suffisamment grande par rapport {\`a} toutes les constantes $c_{i}$ qui interviendront dans la suite. Posons alors $S_{0}:=C_{0}^{3}\mathfrak{a}$, $S:=C_{0}^{6}\mathfrak{a}$, \begin{equation*}\begin{split}U:=&C_{0}^{25g}(\mathfrak{a}\log\mathfrak{e})\left(\mathfrak{a}^{y}+\frac{D}{\log\mathfrak{e}}\log\left(e+\frac{D}{\log\mathfrak{e}}\right)\right)^{1/t}\\ &\times\prod_{i=1}^{n}{\left(1+\frac{D\underset{0\le s\le(g+1)C_{0}^{6}\mathfrak{a}}{\mathrm{max}}{\{h(sp_{i})\}}+(\mathfrak{e}\mathfrak{a}\Vert u_{i}\Vert_{v_{0}})^{\rho_{i}}}{\mathfrak{a}\log\mathfrak{e}}\right)^{g_{i}/t}}\end{split}\end{equation*}(il s'agit essentiellement de $C_{0}^{25g}U_{0}$ mais o\`u la constante \og ind{\'e}finie\fg\ $c_{\theacstsept}$ qui est en indice du $\max{\{h(sp_{i})\}}$ dans $U_{0}$ est remplac{\'e}e par $(g+1)C_{0}^{6}$), \begin{equation*}\widetilde{T}:=\frac{C_{0}U}{S_{0}\log\mathfrak{e}},\qquad T:=[\widetilde{T}],\end{equation*}\begin{equation*}\widetilde{D}_{i}:=\frac{U}{C_{0}^{2}\left(D\max_{0\le s\le (g+1)S}{\{h(sp_{i})\}}+(\mathfrak{e} S\Vert u_{i}\Vert_{v_{0}})^{\rho_{i}}+S_{0}\log\mathfrak{e}\right)}\end{equation*}si $1\le i\le n$. Soit {\'e}galement\begin{equation*}D_{0}^{\flat}:=\left[\frac{S_{0}\log\mathfrak{e}}{DC_{0}^{3}}\right]\quad\text{et}\quad\widetilde{D_{0}}:=\frac{U}{C_{0}^{4}\left(D\log\left(e+\frac{D}{\log\mathfrak{e}}\right)+S_{0}^{y}\log\mathfrak{e}\right)}\ \cdotp\end{equation*}La d{\'e}finition~\ref{defiaseptdex} introduit un nombre r{\'e}el $x>0$ et l'on note $D_{i}:=[x\widetilde{D}_{i}]$, pour tout $i\in\{1,\ldots,n\}$. Consid{\'e}rons alors pour $(P_{\lambda_{0}})$ la famille des polyn{\^o}mes de Matveev $$(\delta_{D_{0}^{\flat}}(X;\lambda_{0}))_{0\le\lambda_{0}\le D_{0}}$$d{\'e}finie au \S~\ref{subsec:poidsdeladroiteaffine}.\par Voici r{\'e}sum{\'e}es en quelques lignes les principales conditions que satisfont ces param{\`e}tres.
\begin{lemm}\label{lemmearchi:choixdesparametres}
On a 
\begin{enumerate}
\item[\ding{192}] $\widetilde{T}\ge C_{0}^{2}$, $S/S_{0}\ge C_{0}^{2}$, $S_{0}\ge C_{0}^{2}$,
\item[\ding{193}] $\widetilde{T}\ge C_{0}\max{\left\{\widetilde{D}_{0}/(S+1)^{1-y},\widetilde{D}_{1},\ldots,\widetilde{D}_{n}\right\}}$,
\item[\ding{194}] $\widetilde{T}^{g-t+1}(S+1)\le C_{0}\widetilde{D}_{0}\widetilde{D}_{1}^{g_{1}}\cdots\widetilde{D}_{n}^{g_{n}}$,
\item[\ding{195}] $U\ge C_{0}^{3/2}D\aleph((P_{\lambda_{0}}))$,
\item[\ding{196}] $S_{0}\log\mathfrak{e}\ge C_{0}^{3}D\max{\{1,h(V)\}}$,
\item[\ding{197}] $U\ge C_{0}^{2}D_{i}\left(D\max_{0\le s\le(g+1)S}{\{h(sp_{i})\}}+(\mathfrak{e} S\Vert u_{i}\Vert_{v_{0}})^{\rho_{i}}\right)$ pour tout $1\le i\le n$.
\end{enumerate}
\end{lemm}
\begin{proof}
Les points~\ding{192} et~\ding{193} d{\'e}coulent imm{\'e}diatement de la d{\'e}finition des param{\`e}tres. L'in{\'e}galit{\'e}~\ding{194} correspond pr{\'e}cis{\'e}ment {\`a} la d{\'e}finition de $U_{0}$ {\`a} une constante pr{\`e}s. La condition~\ding{195} d{\'e}coule du lem\-me~\ref{lemme:estimationdespoids} et les points~\ding{196} et~\ding{197} sont {\'e}vidents {\`a} partir des d{\'e}finitions de $\mathfrak{a}$ et des $\widetilde{D}_{i}$. 
\end{proof}
Accessoirement, on pourra aussi noter que $U\ge C_{0}^{2}D\log(D_{0}S)$. La condition~\ding{194} et le lemme~\ref{lemme24deasept} avec $H_{1}=\{0\}$ impliquent $x\le 1$.
\subsubsection{$v_{0}$ ultram{\'e}trique}
Reprenons les notations $\mathfrak{r},\mathfrak{a}$ du th{\'e}or{\`e}me~\ref{theodeuxasept}. Posons $S_{0}:=C_{0}^{3}\mathfrak{a}$, $S:=C_{0}^{6}\mathfrak{a}$, 
\begin{equation*}\begin{split}U:=& C_{0}^{25g}(\mathfrak{a}\log\mathfrak{r})\left(\mathfrak{a}^{y}+\frac{D}{\log\mathfrak{r}}\log\left(e+\frac{D}{\log\mathfrak{r}}\right)\right)^{1/t}\\ &\quad\prod_{i=1}^{n}{\left(1+\frac{D\underset{0\le s\le(g+1)C_{0}^{6}\mathfrak{a}}{\mathrm{max}}\{h(sp_{i})\}}{\mathfrak{a}\log\mathfrak{r}}\right)^{g_{i}/t}},\end{split}\end{equation*}puis $\widetilde{T}:=C_{0}U/(S_{0}\log\mathfrak{r})$, $T:=[\widetilde{T}]$ et \begin{equation*}\widetilde{D}_{i}:=\frac{U}{C_{0}^{2}\left(D\max_{s\le(g+1)S}{\{h(sp_{i})\}}+S_{0}\log\mathfrak{r}\right)}\quad\text{pour $1\le i\le n$.}\end{equation*}Soit \begin{equation*}D_{0}^{\flat}:=\left[\frac{S_{0}\log\mathfrak{r}}{C_{0}^{3}D}\right]\quad\text{et}\quad\widetilde{D}_{0}:=\frac{U}{C_{0}^{4}\left(D\log\left(e+\frac{D}{\log\mathfrak{r}}\right)+S_{0}^{y}\log\mathfrak{r}\right)}\ \cdotp\end{equation*}Ensuite, comme dans le cas archim{\'e}dien, nous prenons le nombre r{\'e}el $x>0$ de la d{\'e}finition~\ref{defiaseptdex} (p.~\pageref{defiaseptdex}), nous formons les entiers $D_{i}:=[x\widetilde{D}_{i}]$, $i\in\{1,\ldots,n\}$, et la famille $(P_{\lambda_{0}})$ est la  m{\^e}me que celle du cas archim{\'e}dien.\par Les conditions remplies par ces param{\`e}tres sont les conditions~\ding{192}, \ding{193}, \ding{194} du lemme~\ref{lemmearchi:choixdesparametres} (en particulier, on a $x\le 1$) et celles apport{\'e}es par le lemme suivant.
\begin{lemm}\label{lemmeultra:choixdesparametres}
On a
\begin{enumerate}
\item[\ding{205}] $U\ge C_{0}^{3/2}D\aleph((P_{\lambda_{0}}))$,
\item[\ding{206}] $S_{0}\log\mathfrak{r}\ge C_{0}\left(D\max{\{1,h(V)\}}+\log S_{0}\right)$,
\item[\ding{207}] $U\ge C_{0}^{2}DD_{i}\max_{s\le (g+1)S}{\{h(sp_{i})\}}$.
\end{enumerate}
\end{lemm}
\begin{proof}
Ces trois in{\'e}galit{\'e}s sont faciles {\`a} v{\'e}rifier {\`a} partir du choix des param{\`e}tres, la premi{\`e}re, par exemple, {\'e}tant une cons{\'e}quence de la majoration $$\aleph((P_{\lambda_{0}}))\le c_{\theaseptarchicst}\left(D_{0}\log\left(e+\frac{S}{D_{0}^{\flat}}\right)+TD_{0}^{\flat}+\frac{D_{0}\log\mathfrak{r}}{D}\right)$$induite par le lemme~\ref{lemme:estimationdespoids}.
\end{proof}
\begin{remas}
\noignorespaces\begin{enumerate}\item
La pr{\'e}sence du logarithme de $S_{0}$ dans la condition~\ding{206}, pr{\'e}sence qui sera requise lors de l'extrapolation $p$-adique\footnote{Voir le lemme~\ref{lemmedinterpolationpadique} et en particulier le r{\'e}el $\kappa$.}, explique la modification du param{\`e}tre $\mathfrak{a}$ par rapport au cas archim{\'e}dien, avec l'ajout du terme $$\frac{\log^{+}\left((\log(\mathfrak{r}))^{-1}\right)}{\log(\mathfrak{r})}\ .$$
\item Il est facile de v{\'e}rifier qu'avec ces choix les hypoth{\`e}ses du lemme~\ref{lemmenonzeroasept} sont satisfaites et, par cons{\'e}quent, qu'au moins un des entiers $D_{1},\ldots,D_{n}$ est non nul. En revanche, bien que $\widetilde{D}_{0}$ soit clairement sup{\'e}rieur {\`a} $1$, il se pourrait que $D_{0}$ soit nul\footnote{Je ne sais pas si cette {\'e}ventualit\'e peut se produire.}. Cela n'a (paradoxalement) aucune cons{\'e}quence dans la suite de la d\'emonstration.
\item Lorsque dans la preuve on choisit de mettre sur la partie $\mathbb{G}_{\mathrm{a}}$ la base des mon{\^o}mes usuels, cela remplace le terme $\mathfrak{a}^{y}+\frac{D}{\log\mathfrak{e}}\log(e+D/\log\mathfrak{e})$ qui est dans $U$ par $\mathfrak{a}^{y}+\frac{D\log\mathfrak{a}}{\log\mathfrak{e}}$. Cela rajoute donc une d{\'e}pendance suppl{\'e}mentaire en le logarithme de la hauteur du sous-espace $V$.
\end{enumerate}
\end{remas}
\subsubsection{Choix d'une base de $W$}
\label{subsec:constructionbaseW}
Soit $w_{0}$ la base canonique de $t_{\mathbb{G}_{\mathrm{a}}}$ et $(w_{1},\ldots,w_{g-t})$ une base de $V\otimes\overline{\mathbf{Q}}$ fournie par le lemme~\ref{lemmedesiegelabsolu} (l'espace $V\otimes\overline{\mathbf{Q}}$ {\'e}tant identifi{\'e} {\`a} un sous-espace de $\overline{\mathbf{Q}}^{g}$ \textit{via} la base $\mathbf{e}$ de $t_{G}$ introduite au \S~\ref{donnesgeneralesasept}). Par d{\'e}finition, nous avons $\Vert w_{0}\Vert_{v}=1$ pour toute place $v$ de $K$ et $h_{\mathrm{L}^{2}}(w_{1})+\cdots+h_{\mathrm{L}^{2}}(w_{g-t})\le h(V)+g\log(g)$. Nous supposerons que ces vecteurs sont d{\'e}finis sur $K$ (voir pr{\'e}ambule) et que chacune des normes $\Vert w_{j}\Vert_{v}$ avec $j\in\{1,\ldots,g-t\}$ et $v$ une place de $K$ au-dessus de $v_{0}$ est sup{\'e}rieure ou {\'e}gale {\`a} $1$. Cela est toujours possible quitte {\`a} multiplier $w_{i}$ par un nombre rationnel convenable et {\`a} utiliser l'invariance par homoth{\'e}tie des hauteurs $\mathrm{L}^{2}$. De cette mani{\`e}re, nous fixons une base $\mathbf{w}:=(w_{0},\ldots,w_{g-t})$ de $W$. Tous les {\'e}nonc{\'e}s qui vont suivre jusqu'{\`a} la fin du \S~\ref{subsec:constructionasept} restent vrais avec une base quelconque de $W$.
\subsection{Estimations ultram{\'e}triques d'un coefficient de Taylor}
\label{subsec:estultracoefftay}
Dans tout ce paragraphe, $v$ d{\'e}signe une place ultram{\'e}trique du corps de nombres $K$ et $p$ la caract{\'e}ristique r{\'e}siduelle de $v$.
\subsubsection{}Nous rappelons la notion de \emph{taille d'un sous-sch{\'e}ma formel lisse} telle qu'elle a {\'e}t{\'e} d{\'e}finie par Bost au \S~$3.1$ de~\cite{bost6}. Rappelons qu'au \S~\ref{donnesgeneralesasept}, nous avons introduit un mod{\`e}le lisse $\mathcal{G}\to\spec\mathcal{O}_{k}[\frac{1}{m}]$ de $G$. De la sorte, si $v$ ne divise pas $m$, nous pouvons consid{\'e}rer le compl{\'e}t{\'e} formel $\widehat{\mathcal{G}}_{v}$ de $\mathcal{G}\times\spec(\mathcal{O}_{v})$ \`a l'origine ($\mathcal{O}_{v}$ {\'e}tant l'anneau de valuation du compl{\'e}t{\'e} $K_{v}$). C'est un groupe formel lisse sur $\spec(\mathcal{O}_{v})$ et le choix de coordonn{\'e}es locales {\'e}tales au voisinage de l'origine fournit un isomorphisme de sch\'emas formels (sur $\mathcal{O}_{v}$) \begin{equation*}\widehat{\mathcal{G}}_{v}\simeq\widehat{\mathbf{A}}_{\mathcal{O}_{v}}^{g}:=\specf\mathcal{O}_{v}[[X_{1},\ldots,X_{g}]]\ .\end{equation*}Si $\mathfrak{X}$ est un sous-sch{\'e}ma formel lisse de $\widehat{G}_{K_{v}}\simeq\widehat{\mathcal{G}}_{v}\widehat{\otimes}\spec K_{v}$, on dispose d'un nombre r{\'e}el $R_{\mathcal{G},v}(\mathfrak{X})\in [0,1]$, appel{\'e} \emph{taille de $\mathfrak{X}$ relativement au mod{\`e}le $\widehat{\mathcal{G}}_{v}$ de $\widehat{G}_{K_{v}}$}, d{\'e}fini de la mani{\`e}re suivante. Consid{\'e}rons l'image de $\mathfrak{X}$ (not{\'e}e encore $\mathfrak{X}$) dans $\widehat{\mathbf{A}}_{K_{v}}^{g}$ \textit{via} le choix de coordonn{\'e}es pr{\'e}c{\'e}dent. Le groupe $\mathrm{Aut}(\widehat{\mathbf{A}}_{K_{v}}^{g})$ des automorphismes de $\widehat{\mathbf{A}}_{K_{v}}^{g}$ s'identifie {\`a} l'ensemble des $g$-uplets de s{\'e}ries formelles $f=(f_{1},\ldots,f_{g})\in K_{v}[[X_{1},\ldots,X_{g}]]^{g}$ tels que $f(0)=0$ et la matrice jacobienne $\mathrm{D}_{0}f=(\partial f_{i}/\partial x_{j}(0))_{i,j}$ soit inversible. Pour $\varphi=\sum_{\mathbf{n}\in\mathbf{N}^{g}}{a_{\mathbf{n}}\mathbf{X}^{\mathbf{n}}}\in K_{v}[[\mathbf{X}]]$ et $r>0$, on note $$\Vert\varphi\Vert_{r}:=\sup_{\mathbf{n}\in\mathbf{N}^{g}}{\vert a_{\mathbf{n}}\vert_{v}r^{\vert\mathbf{n}\vert}}\in[0,+\infty]$$et \begin{equation*}G_{\omega}(r):=\left\{f\in\mathrm{Aut}(\widehat{\mathbf{A}}_{K_{v}}^{g})\,;\,\mathrm{D}_{0}f\in\mathrm{GL}_{g}(\mathcal{O}_{v})\ \text{et}\ \Vert f\Vert_{r}:=\max_{1\le i\le g}{\Vert f_{i}\Vert_{r}}\le r\right\}\ \cdotp\end{equation*}Alors, par d{\'e}finition, la taille de $\mathfrak{X}$ est \begin{equation*}R_{\mathcal{G},v}(\mathfrak{X}):=\sup{\left\{r\in\,]0,1]\,;\,\exists\,f\in G_{\omega}(r)\,;\,f^{*}(\mathfrak{X})=\widehat{\mathbf{A}}_{K_{v}}^{d}\times\{0\}\right\}}\end{equation*}o\`u $d:=\dim\mathfrak{X}$. La borne sup\'erieure est prise dans $[0,1]$, ainsi $R_{\mathcal{G},v}(\mathfrak{X})=0$ si l'ensemble pr{\'e}c{\'e}dent est vide. Dans cette {\'e}criture, $f^{*}(\mathfrak{X})$ d{\'e}signe l'image inverse de $\mathfrak{X}$ par $f$. Le nombre $R_{\mathcal{G},v}(\mathfrak{X})$ est strictement positif lorsque $\mathfrak{X}$ est analytique. Observons {\'e}galement que si $\mathfrak{X}$ provient d'un sous-sch{\'e}ma de $\mathcal{G}\times\spec\mathcal{O}_{v}$, lisse le long de l'origine, alors $R_{\mathcal{G},v}(\mathfrak{X})=1$. Lorsque cette derni{\`e}re condition de lissit{\'e} n'est pas remplie, nous disposons n{\'e}anmoins d'une estimation un peu meilleure que seulement $R_{\mathcal{G},v}(\mathfrak{X})\ge 0$. Supposons que $\mathfrak{X}$ est le compl{\'e}t{\'e} formel le long de l'origine d'un sous-groupe alg{\'e}brique de $G_{K_{v}}$. Rappelons que le groupe formel $\widehat{\mathcal{G}}_{v}$ poss{\`e}de une exponentielle formelle qui, en termes des coordonn{\'e}es $X_{i}$, s'{\'e}crit comme un $g$-uplet $\mathbf{E}=(E_{1}(\mathbf{X}),\ldots,E_{g}(\mathbf{X}))$ de s{\'e}ries formelles de $K_{v}[[\mathbf{X}]]$, telles que les coefficients de $\mathbf{X}^{\mathbf{n}}$ ($\mathbf{n}\in\mathbf{N}^{g}$) de $E_{1},\ldots,E_{g}$ sont de la forme $\alpha_{\mathbf{n}}/\mathbf{n}!$ avec $\alpha_{\mathbf{n}}\in\mathcal{O}_{v}$. On peut normaliser cette exponentielle de sorte que la diff{\'e}rentielle {\`a} l'origine $\mathrm{D}_{0}\mathbf{E}$ soit l'identit{\'e}. Soit $r_{p}=\vert p\vert_{v}^{1/(p-1)}$ (d{\'e}j{\`a} introduit p.~\pageref{pagesept}). Le $g$-uplet $\mathbf{E}$ est un {\'e}l{\'e}ment de $G_{\omega}(r_{p})$ puisque $\vert\mathbf{n}!\vert_{v}\ge r_{p}^{\vert\mathbf{n}\vert-1}$. Au moyen de cette application exponentielle, il est alors ais{\'e} de construire un automorphisme $f\in G_{\omega}(r_{p})$ tel que $f^{*}(\mathfrak{X})=\widehat{\mathbf{A}}^{d}_{K_{v}}\times\{0\}$. Ce qui entra{\^\i}ne la minoration $R_{\mathcal{G},v}(\mathfrak{X})\ge r_{p}$ (pour toute place $v\nmid m$). Toutefois cette estimation n'est vraiment utile qu'en un nombre fini de places $v$ puisque le produit infini $\prod_{p}{r_{p}}$ diverge.\par Revenons au cas d'un sous-sch{\'e}ma formel lisse $\mathfrak{X}$ quelconque de $\widehat{G}_{K_{v}}$. L'espace tangent {\`a} l'origine $t_{\mathfrak{X}}$ de $\mathfrak{X}$ est muni d'une structure enti\`ere en consid{\'e}rant le module \begin{equation*} t_{\mathfrak{X}}\cap t_{\mathcal{G}_{v}}=\left\{z\in t_{\mathfrak{X}}\subseteq t_{G}(K_{v})\,;\ \Vert z\Vert_{v}\le 1\right\}\end{equation*}sur l'anneau $\mathcal{O}_{v}$, ce qui conf{\`e}re {\`a} l'espace dual $t_{\mathfrak{X}}^{\mathsf{v}}$ puis {\`a} l'espace sym{\'e}trique $S^{\ell}(t_{\mathfrak{X}}^{\mathsf{v}})$ de degr{\'e} $\ell\in\mathbf{N}$ une norme not{\'e}e $\Vert.\Vert_{S^{\ell}(t_{\mathfrak{X}}^{\mathsf{v}})}$.\par 
Ces d{\'e}finitions conduisent alors au lemme suivant. \begin{lemm}[lemme~$3.3$ de~\cite{bost6}]\label{lemmecrucial} Soit $\ell\in\mathbf{N}$, $\Omega$ un sous-sch{\'e}ma ouvert de $\mathcal{G}_{v}$ contenant la section nulle et $s$ une fonction r{\'e}guli{\`e}re sur $\Omega$ telle que $s_{K_{v}}$ s'annule ainsi que ses d\'eriv\'ees d'ordre $<\ell$ le long de $\mathfrak{X}$ en l'\'el\'ement neutre de $G_{K_{v}}$. Alors le jet $\mathrm{j}_{\mathfrak{X}}^{\ell}s$ d'ordre $\ell$ le long de $\mathfrak{X}$ en $0$ --- vu comme {\'e}l{\'e}ment de $S^{\ell}(t_{\mathfrak{X}})^{\mathsf{v}}$ \emph{(\footnote{C'est-\`a-dire, si l'on {\'e}crit $s(x)=F(z).s_{0}(x)$ o\`u $x=\exp_{v}(z)$ et $F$ analytique dans un voisinage de $0$ de $t_{G}(\mathbf{C}_{v})$, on a $\mathrm{j}_{\mathfrak{X}}^{\ell}s=\mathrm{jet}_{t_{\mathfrak{X}}}^{\ell}F(0).s_{0}$.})} --- v{\'e}rifie \begin{equation}\Vert\mathrm{j}_{\mathfrak{X}}^{\ell}s\Vert_{S^{\ell}(t_{\mathfrak{X}})^{\mathsf{v}}}\le R_{\mathcal{G},v}(\mathfrak{X})^{-\ell}\ \cdotp\end{equation}
\end{lemm}
Soit $\widehat{H}_{v}$ le compl{\'e}t{\'e} formel {\`a} l'origine du sch{\'e}ma $H\times\spec K_{v}$. C'est un sous-sch{\'e}ma formel lisse de $\widehat{G}_{K_{v}}$. Soit $P=\sum_{\boldsymbol{\lambda}}{p_{\boldsymbol{\lambda}}\mathbf{X}^{\boldsymbol{\lambda}}}$ un polyn{\^o}me de $K_{v}[\mathbf{P}^{N_{1}}\times\cdots\times\mathbf{P}^{N_{n}}]$. Soit $p_{\widetilde{\boldsymbol{\lambda}}}$ un coefficient de $P$ pour lequel $\vert p_{\widetilde{\boldsymbol{\lambda}}}\vert_{v}=\max_{\boldsymbol{\lambda}}{\{\vert p_{\boldsymbol{\lambda}}\vert_{v}\}}$. En appliquant le lemme {\`a} la section de $\mathcal{O}_{\mathcal{G}_{v}}$ d{\'e}finie par $(P/p_{\widetilde{\boldsymbol{\lambda}}})\circ\exp_{v}$ au voisinage de $0$ et $\mathfrak{X}=\widehat{H}_{v}$, nous obtenons le 
\begin{coro}\label{coro:taylorpolyenzero}Soit $\ell$ un entier non nul et $P$ le polyn{\^o}me ci-dessus. Supposons que les d{\'e}riv{\'e}es $\cD_{\mathbf{w}}^{\boldsymbol{\tau}}(P\circ\exp_{v})(0)$ soient toutes nulles lorsque $\boldsymbol{\tau}\in\mathbf{N}^{\dim V}$ v{\'e}rifie $\vert\boldsymbol{\tau}\vert<\ell$. Supposons {\'e}galement que $v$ ne divise pas $m$. Alors la valeur absolue $v$-adique du coefficient de Taylor
\begin{equation}\label{Taylorunasept}\frac{\cD_{w_{1}}^{\tau_{1}}\cdots\cD_{w_{g-t}}^{\tau_{g-t}}}{\tau_{1}!\cdots\tau_{g-t}!}\left(P\circ\exp_{v}(z)\right)_{\vert z=0},\end{equation}pour $(\tau_{1},\ldots,\tau_{g-t})\in\mathbf{N}^{\dim V}$ de longueur $\ell$, est major{\'e}e par \begin{equation*}R_{\mathcal{G},v}(\widehat{H}_{v})^{-\ell}\max_{\boldsymbol{\lambda}}{\{\vert p_{\boldsymbol{\lambda}}\vert_{v}\}}\prod_{i=1}^{g-t}{\Vert w_{i}\Vert_{v}^{\tau_{i}}}\ .\end{equation*}
\end{coro} 
 D{\'e}signons par $\chi_{H}$ la somme finie \begin{equation}\label{defidechidusubsec}\chi_{H}:=\frac{1}{[K:\mathbf{Q}]}\sum_{\genfrac{}{}{0pt}{}{v\nmid m}{v\ \text{ultram{\'e}trique}}}{[K_{v}:\mathbf{Q}_{p}]\log R_{\mathcal{G},v}(\widehat{H}_{v})^{-1}}\in\mathbf{R}^{+}\ \cdotp\end{equation}Nous disposons d'une estimation \emph{uniforme} de $\chi_{H}$ qui nous a {\'e}t{\'e} signal{\'e}e par J.-B.~Bost.
\begin{lemm}\label{lemme:Raynaud}
Il existe une constante $c_{\cst}\ge 1$ ne d{\'e}pendant que de $G$ (en particulier ind{\'e}pendante du corps $K$ et du groupe $H$) telle que $\chi_{H}\le c_{\theconstante}$. 
\end{lemm}
Ce r{\'e}sultat est une cons{\'e}quence du th{\'e}or{\`e}me suivant de Raynaud.
\begin{theo}[corollaire du th{\'e}or{\`e}me~$4$, \S~$7.5$, de~\cite{bosch}, p.~$187$]\label{thm:Raynaud}
Soit $k_{0}$ un corps de nombres et $A_{1}\subseteq A_{2}$ deux vari{\'e}t{\'e}s ab{\'e}liennes sur $k_{0}$. Soit $\mathcal{A}_{1},\mathcal{A}_{2}$ leurs mod{\`e}les de N{\'e}ron respectifs sur $\mathcal{O}_{k_{0}}$. Soit $v$ une place finie de $k_{0}$, de caract{\'e}ristique r{\'e}siduelle $p$ et d'indice de ramification $e((k_{0})_{v}\mid\mathbf{Q}_{p})$ \emph{(\footnote{De la sorte, si $\varpi$ est une uniformisante de l'anneau de valuation $\mathcal{O}_{v}\subseteq (k_{0})_{v}$, l'id{\'e}al engendr{\'e} par $p$ est $(\varpi)^{e((k_{0})_{v}\mid\mathbf{Q}_{p})}$.})}. Si $e((k_{0})_{v}\mid\mathbf{Q}_{p})<p-1$ et si $\mathcal{A}_{2}$ a bonne r{\'e}duction en $v$ alors $\mathcal{A}_{1}$ a bonne r{\'e}duction en $v$ et l'unique morphisme $\mathcal{A}_{1}\to\mathcal{A}_{2}$ issu de la propri{\'e}t{\'e} de mod{\`e}le de N{\'e}ron de $\mathcal{A}_{2}$ est une immersion ferm{\'e}e.   
\end{theo}
\begin{lemm}\label{lemme:modelelissetore}
Soit $G_{1}$ un sous-groupe alg\'ebrique connexe de $\mathbb{G}_{\mathrm{a}}^{d_{0}}\times\mathbb{G}_{\mathrm{m}}^{d_{1}}$, d\'efini sur un corps de nombres $k_{0}$. Il existe un sous-sch\'ema en groupes $\mathcal{G}_{1}\hookrightarrow\mathbb{G}_{\mathrm{a},\mathcal{O}_{k_{0}}}^{d_{0}}\times\mathbb{G}_{\mathrm{m},\mathcal{O}_{k_{0}}}^{d_{1}}$, lisse sur $\spec\mathcal{O}_{k_{0}}$ et de fibre g\'en\'erique $G_{1}$.
\end{lemm}
\begin{proof}
Par d\'ecomposition, il suffit de traiter les cas $G_{1}\subseteq\mathbb{G}_{\mathrm{a}}^{d_{0}}$ et $G_{1}\subseteq\mathbb{G}_{\mathrm{m}}^{d_{1}}$. Dans le premier cas, on choisit le fibr\'e vectoriel $\mathbf{V}\left((t_{G_{1}}(k_{0})\cap\mathcal{O}_{k_{0}}^{d_{0}})^{\mathsf{v}}\right)$. Dans le second cas, il existe un sous-groupe facteur direct $\Phi$ de $\mathbf{Z}^{d_{1}}$ tel que, si $M:=\mathbf{Z}^{d_{1}}/\Phi$ et si $\mathbf{Z}[M]$ est la $\mathbf{Z}$-alg\`ebre engendr\'ee par $M$, on a $G_{1}=\spec\mathbf{Z}[M]\times_{\spec\mathbf{Z}}\spec k_{0}$. On choisit $\mathcal{G}_{1}:=\spec\mathbf{Z}[M]\times_{\spec\mathbf{Z}}\spec\mathcal{O}_{k_{0}}$, qui est lisse sur $\spec\mathcal{O}_{k_{0}}$ car $M$ est sans torsion. 
\end{proof}
\begin{proof}[D{\'e}monstration du lemme~\ref{lemme:Raynaud}]
Le groupe alg{\'e}brique $G$ est connexe par hypoth{\`e}se. D'apr{\`e}s le th{\'e}or{\`e}me de d{\'e}composition de Chevalley, et apr{\`e}s une {\'e}ventuelle extension finie du corps de nombres de d{\'e}finition $k_{0}$ de $G$, il existe des entiers naturels $d_{0}$ et $d_{1}$ et une vari{\'e}t{\'e} ab{\'e}lienne $A$, d\'efinie sur $k_{0}$, tels que $G$ soit une extension du groupe lin{\'e}aire $G_{0}:=\mathbb{G}_{\mathrm{a},k_{0}}^{d_{0}}\times\mathbb{G}_{\mathrm{m},k_{0}}^{d_{1}}$ par $A$. Le sous-groupe $H$ de $G$ est alors une extension d'un sous-groupe alg\'ebrique (connexe) $G_{1}$ de $G_{0}$ par une sous-vari\'et\'e ab\'elienne $B$ de $A$. Soit $\mathcal{G}_{0}$ (\emph{resp}. $\mathcal{A}$) le sch\'ema en groupes $\mathbb{G}_{\mathrm{a},\mathcal{O}_{k_{0}}}^{d_{0}}\times\mathbb{G}_{\mathrm{m},\mathcal{O}_{k_{0}}}^{d_{1}}$ (\emph{resp}. le mod\`ele de N\'eron de $A$ sur $\spec\mathcal{O}_{k_{0}}$). De m{\^e}me, soit $\mathcal{G}_{1}\hookrightarrow\mathcal{G}_{0}$ le mod\`ele lisse de $G_{1}$ donn\'e par le lemme~\ref{lemme:modelelissetore} et soit $\mathcal{B}$ le mod\`ele de N\'eron de $B$. Il existe un ensemble fini $F$ de places ultram\'etriques de $k_{0}$ (qui ne d\'epend pas de $H$) tel que, si $v\not\in F$, les propri\'et\'es suivantes sont v\'erifi\'ees~:
\begin{enumerate}
\item[1)] La suite $0\to\mathcal{G}_{0,v}\to\mathcal{G}_{v}\to\mathcal{A}_{v}\to 0$ est exacte (l'indice $v$ signifie que nous avons consid\'er\'e le produit fibr\'e avec $\spec\mathcal{O}_{v}$),
\item[2)] Le sch\'ema $\mathcal{B}_{v}$ est un sous-sch\'ema ab\'elien de $\mathcal{A}_{v}$.\end{enumerate}Cette seconde assertion d\'ecoule du th\'eor\`eme~\ref{thm:Raynaud}~: elle n'est pas satisfaite seulement si $v$ est une place de mauvaise r\'eduction pour $\mathcal{A}$ ou si la caract\'eristique r\'esiduelle $p$ de $v$ est plus petite que $e((k_{0})_{v}\mid\mathbf{Q}_{p})+1$, quantit\'e elle-m\^eme inf\'erieure \`a $[k_{0}:\mathbf{Q}]+1$ (nombre fini de telles places).\par Soit $v\not\in F$. Posons $d:=d_{0}+d_{1}$ et $h:=\dim A$. On choisit des coordonn\'ees locales $x_{1},\ldots,x_{d}$ (\emph{resp}. $y_{1},\ldots,y_{h}$) sur $\mathcal{G}_{0,v}$ (\emph{resp}. $\mathcal{A}_{v}$) \'etales en l'origine, telles que $x_{1},\ldots,x_{\dim G_{1}}$ (\emph{resp}. $y_{1},\ldots,y_{\dim B}$) soient des coordonn\'ees locales de $\mathcal{G}_{1,v}$ (\emph{resp}. $\mathcal{B}_{v}$) (crit\`ere de Jacobi, voir la proposition du chapitre~2 de~\cite{bosch}). On obtient ainsi un syst\`eme de coordonn\'ees locales \'etales qui param\'etrisent $H_{v}:=H\times\spec (k_{0})_{v}$ au voisinage de l'origine, et qui, relativement \`a $(x_{1},\ldots,y_{h})$, est de taille $1$. Ainsi, pour $v\not\in F$, on a $R_{\mathcal{G},v}(\widehat{H_{v}})=1$. En les autres places, on utilise la minoration $R_{\mathcal{G},v}(\widehat{H_{v}})\ge\vert p\vert_{v}^{1/(p-1)}$ {\'e}voqu{\'e}e un peu plus haut. Pour conclure il suffit d'observer que $\chi_{H}$ est inf\'erieur \`a une expression du m\^eme type o\`u le corps $K$ est remplac\'e par $k_{0}$ et $v$ parcourant les places ultram\'etriques de $k_{0}$ qui ne divisent pas $m$. L'existence de la constante $c_{\theconstante}$ s'ensuit. 
\end{proof}
\begin{rema}
En r{\'e}alit{\'e}, il n'est pas n{\'e}cessaire d'avoir une estimation aussi fine de $\chi_{H}$ pour d{\'e}montrer les th{\'e}or{\`e}mes~\ref{theounasept} et~\ref{theodeuxasept}. Lorsque nous laissons $\chi_{H}$ \og ind{\'e}termin{\'e}\fg\ (ce que nous ferons dans la suite), nous constatons qu'il appara{\^\i}t dans le facteur $\max{\{1,h(V)\}}$ (intervenant dans le param{\`e}tre $\mathfrak{a}$), qui devient alors $\max{\{1,h(V),\chi_{H}\}}$. Autrement dit, il suffirait de montrer l'existence d'une constante $c_{\cst}$ telle que $\chi_{H}\le c_{\theconstante}\max{\{1,h(V)\}}$ pour obtenir \emph{exactement} les m{\^e}mes {\'e}nonc{\'e}s que les th{\'e}or{\`e}mes~\ref{theounasept} et~\ref{theodeuxasept}. Et cela est possible de mani{\`e}re {\'e}l{\'e}mentaire sans avoir recours au r{\'e}sultat de Raynaud. La d{\'e}marche consiste {\`a} se ramener comme ci-dessus au cas ab{\'e}lien puis {\`a} comparer la structure enti{\`e}re sur $t_{B}$ donn{\'e}e par le mod{\`e}le de N{\'e}ron $\mathcal{B}$ de $B$ et celle donn{\'e}e par le module satur{\'e} $t_{B}\cap t_{\mathcal{A}}$, {\`a} partir duquel se calcule la norme des jets le long de $t_{B}$. Le calcul du quotient $(t_{B}\cap t_{\mathcal{A}})/t_{\mathcal{B}}$ (voir, par exemple, p.~$33$ de~\cite{bost2}) fournit alors une constante $c_{\cst}$ telle que $\chi_{H}\le c_{\theconstante}\max{\{1,\log\deg B\}}$ o\`u le degr{\'e} est relatif {\`a} un plongement (quelconque) de $A$ dans un espace projectif. Il ne reste plus qu'{\`a} observer que $\log\deg B$ est du m{\^e}me ordre de grandeur que la hauteur de $t_{B}$ (voir, par exemple, le lemme~$4.8$ de~\cite{artepredeux}), elle-m{\^e}me major{\'e}e par $c_{\cst}\max{\{1,h(V)\}}$ pour une certaine constante $c_{\theconstante}$.  
\end{rema}
\subsubsection{} Avant de passer aux estimations archim{\'e}diennes, nous aurons besoin d'une cons{\'e}quence (qui en est aussi une g{\'e}n{\'e}ralisation en quelque sorte) du corollaire~\ref{coro:taylorpolyenzero}. Soit maintenant $P$ un polyn{\^o}me de l'espace $E$, consid{\'e}r{\'e} au \S~\ref{subsec:rangsyslineaire} et $F:=P\circ\exp_{(G_{0}\times G)(\mathbf{C}_{v_{0}})}$ l'application associ{\'e}e {\`a} $P$ en la place $v_{0}$, d{\'e}finie par~\eqref{defideFPv}. Soit $\ell\in\mathbf{N}$, $z$ un vecteur de $\mathscr{T}_{v_{0}}$ d'exponentielle $K$-rationnelle et $z_{0}$ un {\'e}l{\'e}ment de $t_{G_{0}}(K)\simeq K$ (le corps $K$ est plong\'e dans $K_{v_{0}}$). Dans l'{\'e}nonc{\'e} qui va suivre, nous supposerons que $F$ s'annule {\`a} l'ordre $\ell$ le long de $W$ au point $(z_{0},z)$, ou, en d'autres termes:\begin{equation}\label{diversesannulations}\forall\,\mathbf{i}=(i_{0},\ldots,i_{g-t})\in\mathbf{N}^{\dim W}\;\ \vert\mathbf{i}\vert\le\ell-1,\quad \cD_{\mathbf{w}}^{\mathbf{i}}F(z_{0},z)=0\end{equation}(si $\ell=0$ cette condition est vide). Pour chaque entier $i\in\{1,\ldots,n\}$, nous fixons un entier $\varepsilon_{i}$ de l'ensemble $\{0,\ldots,N_{i}\}$ pour lequel $\theta_{v_{0},i,\varepsilon_{i}}(z)\ne 0$.
\begin{lemm}\label{sdgfkjsgvjq}
Dans ces conditions, \'etant donn\'e un multiplet $(\tau_{0},\ldots,\tau_{g-t})$ de longueur $\ell$, le coefficient de Taylor tordu \begin{equation}\label{coeffdetaylorultra}\frac{1}{\prod_{i=1}^{n}{\theta_{v_{0},i,\varepsilon_{i}}(z)^{D_{i}}}}\frac{\mathcal{D}_{w_{0}}^{\tau_{0}}\cdots\mathcal{D}_{w_{g-t}}^{\tau_{g-t}}}{\tau_{0}!\cdots\tau_{g-t}!}F(z_{0},z),\end{equation}appartient {\`a} $K$.
\end{lemm} 
\begin{proof}
Nous allons utiliser les formules de translations sur le groupe alg{\'e}brique $G(K_{v_{0}})$ pour nous ramener en $0$\footnote{Il s'agit de \og l'astuce d'Anderson-Baker-Coates\fg\ (voir, par exemple, la d\'emonstration du lemme~$13$ de~\cite{gaudron1}).}. Soit $\widetilde{Q}$ le polyn{\^o}me d{\'e}fini par \begin{equation}\widetilde{Q}(\mathbf{X}_{1},\ldots,\mathbf{X}_{n})=\sum_{\boldsymbol{\lambda}}{p_{\boldsymbol{\lambda}}\frac{P_{\lambda_{0}}^{(\tau_{0})}(z_{0})}{\tau_{0}!}\prod_{i=1}^{n}{A^{(i)}(\Psi_{v_{0},i,\varepsilon_{i}}(z),\mathbf{X}_{i})^{\lambda_{i}}}}\end{equation}o\`u $(A^{(i)})_{i=1\cdots n}$ sont les multiplets de polyn{\^o}mes repr{\'e}sentant les formules d'addition sur $G_{K_{v_{0}}}$ au voisinage de $x=\exp_{v_{0}}(z)$ (voir \S~\ref{miseenplace}). La formule de Leibniz et l'hypoth{\`e}se~\eqref{diversesannulations} montrent que le coefficient~\eqref{coeffdetaylorultra} {\'e}gale \begin{equation}\label{coeffdeuxaspet}\frac{\cD_{w_{0}}^{\tau_{0}}\cdots\cD_{w_{g-t}}^{\tau_{g-t}}}{\tau_{0}!\cdots\tau_{g-t}!}\left(\frac{F(z_{0}+z_{0}',z+z')}{\prod_{i=1}^{n}{\theta_{v_{0},i,\varepsilon_{i}}(z+z')^{D_{i}}}}\right)_{\vert (z_{0}',z')=(0,0)}\ .\end{equation}Or, par d{\'e}finition, \begin{equation*}\frac{\theta_{v_{0},i,j}(z+z')}{\theta_{v_{0},i,\varepsilon_{i}}(z+z')}=\frac{A_{j}^{(i)}(\Psi_{v_{0},i,\varepsilon_{i}}(z),\Theta_{v_{0},i}(z'))}{A_{\varepsilon_{i}}^{(i)}(\Psi_{v_{0},i,\varepsilon_{i}}(z),\Theta_{v_{0},i}(z'))}\end{equation*}pour $z'$ proche de $0$, donc \begin{equation*}\begin{split}\frac{F(z_{0}+z_{0}',z+z')}{\prod_{i=1}^{n}{\theta_{v_{0},i,\varepsilon_{i}}(z+z')^{D_{i}}}}=&\left(\sum_{\boldsymbol{\lambda}}{p_{\boldsymbol{\lambda}}P_{\lambda_{0}}(z_{0}+z_{0}')\prod_{i=1}^{n}{A^{(i)}(\Psi_{v_{0},i,\varepsilon_{i}}(z),\Theta_{v_{0},i}(z'))^{\lambda_{i}}}}\right)\\ &\times\frac{1}{\prod_{i=1}^{n}{A_{\varepsilon_{i}}^{(i)}(\Psi_{v_{0},i,\varepsilon_{i}}(z),\Theta_{v_{0},i}(z'))^{D_{i}}}}\ \cdotp\end{split}\end{equation*}En d{\'e}rivant par rapport {\`a} $z_{0}'$, on a \begin{equation*}\frac{\cD_{w_{0}}^{\tau_{0}}}{\tau_{0}!}\left(\frac{F(z_{0}+z_{0}',z+z')}{\prod_{i=1}^{n}{\theta_{v_{0},i,\varepsilon_{i}}(z+z')^{D_{i}}}}\right)_{\vert z_{0}'=0}=\frac{F_{\widetilde{Q},v_{0}}(z')}{\prod_{i=1}^{n}{A_{\varepsilon_{i}}^{(i)}(\Psi_{v_{0},i,\varepsilon_{i}}(z),\Theta_{v_{0},i}(z'))^{D_{i}}}}\ \cdotp\end{equation*}Cette {\'e}galit{\'e} implique en particulier que les d{\'e}riv{\'e}es $\cD_{\mathbf{w}}^{\mathbf{i}}F_{\widetilde{Q},v_{0}}(0)$ pour $\vert\mathbf{i}\vert<\ell-\tau_{0}$ sont toutes nulles. En appliquant alors l'op{\'e}rateur $$\cD_{w_{1}}^{\tau_{1}}\cdots\cD_{w_{g-t}}^{\tau_{g-t}}/(\tau_{1}!\cdots\tau_{g-t}!)$$aux deux membres puis, {\`a} nouveau, la formule de Leibniz pour le second, on d{\'e}duit l'{\'e}galit{\'e} \begin{equation}\begin{split}\label{equationronde}&\frac{1}{\prod_{i=1}^{n}{\theta_{v_{0},i,\varepsilon_{i}}(z)^{D_{i}}}}\frac{\mathcal{D}_{w_{0}}^{\tau_{0}}\cdots\mathcal{D}_{w_{g-t}}^{\tau_{g-t}}}{\tau_{0}!\cdots\tau_{g-t}!}F_{P,v_{0}}(z_{0},z)\\ &\quad=\frac{1}{\prod_{i=1}^{n}{(A_{\varepsilon_{i}}^{(i)}(\Psi_{v_{0},i,\varepsilon_{i}}(z),(1:0:\cdots:0)))^{D_{i}}}}\frac{\mathcal{D}_{w_{1}}^{\tau_{1}}\cdots\mathcal{D}_{w_{g-t}}^{\tau_{g-t}}}{\tau_{1}!\cdots\tau_{g-t}!}F_{\widetilde{Q},v_{0}}(0)\ \cdotp\end{split}\end{equation}Maintenant, par choix de $z$, chacune des coordonn{\'e}es de $\Psi_{v_{0},i,\varepsilon_{i}}(z)$ appartient \`a $K$, ainsi bien s{\^u}r que les coefficients des polyn{\^o}mes $A_{\varepsilon_{i}}^{(i)}$ et des polyn{\^o}mes des formules diff\'erentielles v{\'e}rifi{\'e}es par les $\theta_{v_{0},i,j}$ en $0$ (rappelons que $\Theta_{v_{0},i}(0)=(1,0,\ldots,0)$). Par suite, le membre de droite de~\eqref{equationronde} est clairement un {\'e}l{\'e}ment de $K$, ce qui d{\'e}montre le lemme.
\end{proof}
En chemin, nous avons montr{\'e} que le polyn{\^o}me $\widetilde{Q}$ satisfaisait aux hypoth{\`e}ses du corollaire~\ref{coro:taylorpolyenzero}. De l'{\'e}galit{\'e}~\eqref{equationronde} et d'une estimation imm{\'e}diate des coefficients de $\widetilde{Q}$ d{\'e}coule alors la 
\begin{prop}\label{prop:estultrametriqueasept}
Avec les notations et hypoth{\`e}ses ci-dessus et si $v\nmid m$, la valeur absolue $v$-adique du coefficient~\eqref{coeffdetaylorultra} est major{\'e}e par \begin{equation}\begin{split}\label{majorantpadique}&\max{\{\vert p_{\boldsymbol{\lambda}}\vert_{v}\}}\prod_{i=1}^{n}{\underset{0\le j\le N_{i}}{\mathrm{max}}\left\{\left\vert\frac{\theta_{v_{0},i,j}}{\theta_{v_{0},i,\varepsilon_{i}}}(z)\right\vert_{v}\right\}^{c_{\thedegdespolys}D_{i}}}\underset{\lambda_{0}}{\mathrm{max}}\left\{\left\vert\frac{1}{\tau_{0}!}P_{\lambda_{0}}^{(\tau_{0})}(z_{0})\right\vert_{v}\right\}\\ &\times\prod_{i=1}^{g-t}{\Vert w_{i}\Vert_{v}^{\tau_{i}}}\times R_{\mathcal{G},v}(\widehat{H}_{v})^{-\ell}\times\prod_{i=1}^{n}{\frac{1}{\vert A_{\varepsilon_{i}}^{(i)}(\Psi_{v_{0},i,\varepsilon_{i}}(z),(1,0,\ldots,0))\vert_{v}^{D_{i}}}}\ \cdotp\end{split}\end{equation}Si $v\mid m$ la majoration ci-dessus reste vraie en rempla{\c c}ant $R_{\mathcal{G},v}(\widehat{H}_{v})^{-\ell}$ par $c_{\cst}^{\ell+D_{1}+\cdots+D_{n}}$ o\`u $c_{\theconstante}$ est une constante $\ge 1$.
\end{prop}
 Des estimations du m{\^e}me type aux places archim{\'e}diennes sont plus rudimentaires et font l'objet du paragraphe suivant.
\subsection{Estimations archim{\'e}diennes de coefficients de Taylor}
\label{subsec:estarchicoefftay}
Dans tout ce paragraphe, $v$ est une place archim{\'e}dienne de $K$. Nous voulons obtenir ici une majoration du coefficient tordu~\eqref{coeffdetaylorultra} en la place $v$. Commen{\c c}ons tout d'abord par noter le 
\begin{lemm}\label{lemmee:majarchitaylor}
Il existe une constante $c_{\cst}$ telle que soit v{\'e}rifi{\'e}e la propri{\'e}t{\'e} suivante. Soit $\boldsymbol{\lambda}=(\lambda_{1},\ldots,\lambda_{n})\in\mathbf{N}^{N_{1}+1}\times\cdots\times\mathbf{N}^{N_{n}+1}$, $z=(z_{1},\ldots,z_{n})$ un vecteur de $t_{G}(\mathbf{C}_{v})$, $\mathbf{f}=(f_{1},\ldots,f_{g})$ une base de $t_{G}(\mathbf{C}_{v})$ et $\boldsymbol{\tau}=(\tau_{1},\ldots,\tau_{g})\in\mathbf{N}^{g}$. Pour chaque entier $i\in\{1,\ldots,n\}$, consid{\'e}rons $\varepsilon_{i}\in\{0,\ldots,N_{i}\}$ un entier pour lequel $\theta_{v,i,\varepsilon_{i}}(z)\ne 0$. Alors la valeur absolue du coefficient de Taylor \begin{equation*}\frac{\cD_{\mathbf{f}}^{\boldsymbol{\tau}}}{\boldsymbol{\tau}!}\left(\prod_{i=1}^{n}{\Psi_{v,i,\varepsilon_{i}}^{\lambda_{i}}}\right)(z)\end{equation*}est major{\'e}e par \begin{equation*}c_{\theconstante}^{\vert\boldsymbol{\tau}\vert+\vert\boldsymbol{\lambda}\vert}\left(\prod_{j=1}^{g}{\Vert f_{j}\Vert_{v}^{\tau_{j}}}\right)\prod_{i=1}^{n}{\max_{0\le j\le N_{i}}{\left\{\left\vert\frac{\theta_{v,i,j}}{\theta_{v,i,\varepsilon_{i}}}(z)\right\vert_{v}\right\}}^{c_{\theconstante}\vert\lambda_{i}\vert}}\ \cdotp\end{equation*}
\end{lemm}
Ce r{\'e}sultat quoique technique ne pose aucune difficult{\'e} particuli{\`e}re lorsque l'on se rappelle que pour chaque $i\in\{1,\ldots,n\}$ l'anneau $$k[\Psi_{v,i,\varepsilon_{i}}]:=k[(\theta_{v,i,j}/\theta_{v,i,\varepsilon_{i}})_{0\le j\le N_{i}}]$$est stable par d{\'e}rivation le long d'un vecteur de $t_{G}(k)$. Il faut cependant observer que la constante $c_{\theconstante}$ est bien ind{\'e}pendante de la base $\mathbf{f}$ choisie. Mais cela se voit imm{\'e}diatement en {\'e}crivant les vecteurs $f_{j}/\Vert f_{j}\Vert_{v}$ dans une base orthonorm{\'e}e quelconque de $t_{G}(\mathbf{C}_{v})$ (qui, elle, ne d{\'e}pend que $G$ et de la norme sur $t_{G}(\mathbf{C}_{v})$) et en majorant leurs composantes par $1$.\par Maintenant, consid{\'e}rons une situation semblable {\`a} celle du dernier {\'e}nonc{\'e} du paragraphe pr{\'e}c{\'e}dent. \'Etant donn{\'e} un polyn{\^o}me $P\in E$, un entier $\ell\in\mathbf{N}$ et un vecteur $(z_{0},z)\in t_{G_{0}\times G}(\mathbf{C}_{v_{0}})$ d'exponentielle $K$-rationnelle, il s'agit, sous l'hypoth{\`e}se d'annulation~\eqref{diversesannulations} de donner une borne du coefficient de Taylor tordu~\eqref{coeffdetaylorultra}. La r{\'e}ponse se trouve aussit{\^o}t dans l'{\'e}galit{\'e}~\eqref{equationronde} (qui repose sur les formules d'addition \og explicites\fg\ de $G$) et le lemme~\ref{lemmee:majarchitaylor}, ce qui conduit {\`a} la 
\begin{prop}\label{aseptprop:estiarchimediennes}
Il existe une constante $c_{\cst}\ge 1$\newcounter{archicst}\setcounter{archicst}{\value{constante}} ayant la propri\'et\'e suivante. Pour toute place archim\'edienne $v$ et avec les notations et hypoth{\`e}ses ci-dessus, la valeur absolue $v$-adique du nombre~\eqref{coeffdetaylorultra} ($\vert\boldsymbol{\tau}\vert=\ell$) est major\'ee par \begin{equation*}\begin{split}&c_{\thearchicst}^{T+\log D_{0}+D_{1}+\cdots+D_{n}}\left(\prod_{i=1}^{g-t}{\Vert w_{i}\Vert_{v}^{\tau_{i}}}\right)\max_{\boldsymbol{\lambda}}{\{\vert p_{\boldsymbol{\lambda}}\vert_{v}\}}\max_{\lambda_{0}}{\left\{\left\vert\frac{1}{\tau_{0}!}P_{\lambda_{0}}^{(\tau_{0})}(z_{0})\right\vert_{v}\right\}}\\ &\times\prod_{i=1}^{n}{\frac{1}{\left\vert A_{\varepsilon_{i}}^{(i)}(\Psi_{v_{0},i,\varepsilon_{i}}(z),(1:0:\cdots:0))\right\vert_{v}^{D_{i}}}}\prod_{i=1}^{n}{\max_{0\le j\le N_{i}}{\left\{\left\vert\frac{\theta_{v_{0},i,j}}{\theta_{v_{0},i,\varepsilon_{i}}}(z)\right\vert_{v}\right\}^{c_{\thearchicst}D_{i}}}}\ \cdotp\end{split}\end{equation*}
\end{prop}   
\begin{rema}Le logarithme de $D_{0}$ qui appara{\^\i}t en exposant de $c_{\theconstante}$ provient simplement de la dimension de $E$ (qui vaut $H(G_{0}\times G;D_{0},\ldots,D_{n})$).\end{rema}

\subsection{Construction du polyn{\^o}me auxiliaire}\label{subsec:constructionasept}
Soit $\Upsilon$ l'ensemble d\'efini au d\'ebut du \S~\ref{choixunasept}. Soit $E$ l'espace vectoriel quo\-tient $\left(k[\mathbf{P}]/I(G_{0}\times G)\right)_{\mathbf{D}}$ introduit au \S~\ref{subsec:rangsyslineaire} et $F$ le sous-espace de $\overline{\mathbf{Q}}^{\dim E}$ d\'efini par \begin{equation}\label{defideF}F:=\left\{(p_{\boldsymbol{\lambda}})_{\boldsymbol{\lambda}}\in\overline{\mathbf{Q}}^{\dim E}\,;\ \forall\,(s,\mathbf{t})\in\Upsilon,\ \sum_{\boldsymbol{\lambda}}{p_{\boldsymbol{\lambda}}\mathcal{D}_{\mathbf{w}}^{\mathbf{t}}(\Theta^{\boldsymbol{\lambda}})(s(1,u))}=0\right\}\end{equation}o\`u $\Theta^{\boldsymbol{\lambda}}$ est une notation abr\'eg\'ee pour \begin{equation}\label{asepteq:pagequinze}P_{\lambda_{0}}\prod_{i=1}^{n}{\Theta_{v_{0},i}^{\lambda_{i}}}=P_{\lambda_{0}}\prod_{i=1}^{n}{\prod_{j=0}^{N_{i}}{\theta_{v_{0},i,j}^{\lambda_{i,j}}}}\end{equation}($\lambda_{i}=(\lambda_{i,0},\ldots,\lambda_{i,N_{i}})$ est un multiplet de longueur $D_{i}$). En d'autres termes, cet espace s'identifie \`a l'ensemble des polyn\^omes $P$ de $E\otimes_{v_{0}}\overline{\mathbf{Q}}$ tel que, pour tout $(s,\mathbf{t})\in\Upsilon$, la d{\'e}riv\'ee $\mathcal{D}_{\mathbf{w}}^{\mathbf{t}}F_{P,v_{0}}(s,su)$ est nulle. Nous savons que $F$ n'est pas r\'eduit \`a $\{0\}$. En effet, le syst\`eme lin\'eaire d'\'equations d\'efinissant $F$ a \'et\'e \'etudi\'e au \S~\ref{subsec:rangsyslineaire} et nous avons vu que son rang $\rho$ \'etait major\'e par $C_{0}^{3/2}\frac{S_{0}}{S}\mathscr{H}(G_{0}\times G;D_{0}',\ldots,D_{n}')$. Or, d'une part, le rapport $S_{0}/S$ est inf\'erieur \`a $1/C_{0}^{2}$ (choix des param\`etres) et, d'autre part, comme les param{\`e}tres $D_{i}$, $0\le i\le n$, ne sont pas tous nuls, il existe une constante $c_{\cst}$ telle que $\mathscr{H}(G_{0}\times G;D_{0}',\ldots,D_{n}')\le c_{\theconstante}\dim E$. Nous en d\'eduisons que $\rho\le(\dim E)/2$ (pourvu que $C_{0}$ soit assez grand) et donc $\dim F\ge (\dim E)/2$. Dans ces conditions, le lemme de Siegel absolu \'enonc\'e au \S~\ref{subsec:lemmesiegelabsolu} fournit naturellement un \'el\'ement de $F\setminus\{0\}$ de \og petite\fg\ hauteur, et, plus pr\'ecis\'ement, on a la 
\begin{prop}\label{prop:construction}Il existe une constante $c_{\cst}\ge 1$ et une famille $(p_{\boldsymbol{\lambda}})_{\boldsymbol{\lambda}}\in F\setminus\{0\}$ de hauteur (logarithmique absolue) $\mathrm{L}^{2}$ major\'ee par \begin{equation*} \begin{split}&c_{\theconstante}\left(\log D_{0}+\max_{1\le i\le n}{\{D_{i}\}}+\aleph((P_{\lambda_{0}}))+\sum_{i=1}^{n}{D_{i}\max_{0\le s\le (g+1)S}{\{h(sp_{i})\}}}\right.\\ &\qquad\Bigg.+T(1+\chi_{H}+h_{\mathrm{L}^{2}}(w_{1})+\cdots+h_{\mathrm{L}^{2}}(w_{g-t}))\Bigg)\ \cdotp\end{split}\end{equation*}
\end{prop}
L'\'el\'ement $(p_{\boldsymbol{\lambda}})_{\boldsymbol{\lambda}}$ repr\'esente les coefficients du polyn\^ome auxiliaire $P$ que l'on cherchait \`a construire (ainsi les notations $(p_{\boldsymbol{\lambda}})_{\boldsymbol{\lambda}}$ et $P$ sont-elles d\'esormais fix\'ees jusqu'\`a la fin de la preuve des th\'eor\`emes~\ref{theounasept} et~\ref{theodeuxasept}). Par la suite nous supposerons qu'un des coefficients $p_{\boldsymbol{\lambda}}$ vaut $1$ (ce qui est loisible puisque la hauteur $\mathrm{L}^{2}$ est projective), si bien que chacun des termes locaux intervenant dans $h_{\mathrm{L}^{2}}((p_{\boldsymbol{\lambda}}))$ est positif.
\begin{proof}[D\'emonstration de la proposition~\ref{prop:construction}] Comme nous l'avons mentionn\'e ci-dessus, c'est le lem\-me~\ref{lemmedesiegelabsolu} qui fournit l'\'el\'ement $(p_{\boldsymbol{\lambda}})_{\boldsymbol{\lambda}}$ de $F\setminus\{0\}$ recherch\'e. Sa hauteur est major\'ee par $h(F)/(\dim F)+\log(\dim F)$ (en prenant le vecteur de hauteur minimale parmi ceux de la base apport\'ee par ce lemme). La principale difficult\'e est d'\'evaluer soigneusement la hauteur de $F$ pour ne pas faire appara{\^\i}tre un terme en $T\log(T)$. \par Pour un multi-indice $\boldsymbol{\lambda}$ comme ci-dessus et $(s,\mathbf{t})\in\Upsilon$, consid{\'e}rons un entier $\varepsilon_{i,s}\in\{0,\ldots,N_{i}\}$ pour lequel $\theta_{v_{0},i,\varepsilon_{i,s}}(su)\ne 0$. Posons $\varepsilon_{s}:=(\varepsilon_{1,s},\ldots,\varepsilon_{n,s})$ et d{\'e}signons par $a_{\boldsymbol{\lambda},(s,\mathbf{t})}$ le nombre alg{\'e}brique ({\'e}l{\'e}ment de $K\subseteq K_{v_{0}}$) \begin{equation*}a_{\boldsymbol{\lambda},(s,\mathbf{t})}:=\frac{1}{\mathbf{t}!}\cD_{\mathbf{w}}^{\mathbf{t}}\left(\Psi_{v_{0},\varepsilon_{s}}^{\boldsymbol{\lambda}}\right)(s,su)\end{equation*}o\`u $\Psi_{v_{0},\varepsilon_{s}}^{\boldsymbol{\lambda}}$ est une notation condens{\'e}e pour $P_{\lambda_{0}}\prod_{i=1}^{n}{\Psi_{v_{0},i,\varepsilon_{i,s}}^{\lambda_{i}}}$. L'{\'e}galit{\'e} \begin{equation}\label{eq:eq}\Theta_{v_{0}}^{\boldsymbol{\lambda}}=\left(\prod_{i=1}^{n}{\theta_{v_{0},i,\varepsilon_{i,s}}^{D_{i}}}\right)\times\Psi_{v_{0},\varepsilon_{s}}^{\boldsymbol{\lambda}},\end{equation}la formule de Leibniz et le fait que $\Upsilon$ contienne \`a $(s,\mathbf{t})$ fix{\'e} tous les $(s,\mathbf{t}')$ avec $\vert\mathbf{t}'\vert<\vert\mathbf{t}\vert$ permettent de d{\'e}crire $F$ avec les {\'e}quations \begin{equation*}\forall\,(s,\mathbf{t})\in\Upsilon,\quad\sum_{\boldsymbol{\lambda}}{p_{\boldsymbol{\lambda}}a_{\boldsymbol{\lambda},(s,\mathbf{t})}}=0\ .\end{equation*}Parmi ces {\'e}quations, choisissons-en $\codim F=\dim E-\dim F$ lin{\'e}airement ind{\'e}pendantes et formons alors la matrice $\mathsf{F}$ \`a $\dim E$ lignes et $\codim F$ colonnes constitu\'ee des $a_{\boldsymbol{\lambda},(s,\mathbf{t})}$ correspondants. De la sorte, $\mathsf{F}$ est de rang maximal et les vecteurs colonnes $C_{1},\ldots,C_{\codim F}$ de $\mathsf{F}$ forment une base de l'orthogonal $F^{\perp}$ de $F$ (pour le produit scalaire usuel). La formule de dualit\'e de Schmidt\footnote{Voir~\cite{schmidt1967}.} $h(F)=h(F^{\perp})$ permet alors de calculer la hauteur de $F$ au moyen de ces vecteurs colonnes. Les coordonn\'ees de Pl{\"u}cker (dans la base canonique de $\overline{\mathbf{Q}}^{\dim E}$) des produits ext\'erieurs $C_{1}\wedge\cdots\wedge C_{\codim F}$ conduisent \`a l'\'egalit\'e \begin{equation*}\begin{split}h(F)=&\frac{1}{[K:\mathbf{Q}]}\sum_{\sigma:K\hookrightarrow\mathbf{C}}{[K_{\sigma}:\mathbf{R}]\log\left(\sum_{\mathsf{F}_{0}}{\vert\det\mathsf{F}_{0}\vert_{\sigma}^{2}}\right)^{1/2}}\\ &\quad+\frac{1}{[K:\mathbf{Q}]}\sum_{v\nmid\infty}{[K_{v}:\mathbf{Q}_{p}]\log\max_{\mathsf{F}_{0}}{\{\vert\det\mathsf{F}_{0}\vert_{v}\}}}\end{split}\end{equation*}o\`u $\mathsf{F}_{0}$ parcourt les mineurs maximaux ($\codim F\times\codim F$) de la matrice $\mathsf{F}$. Nous appellerons $h_{\infty}(F)$ (\emph{resp}. $h_{f}(F)$) la premi\`ere (\emph{resp}. la seconde) somme du membre de droite de cette \'egalit\'e.\par\'Etant donn\'e une place archim\'edienne $\sigma$ de $K$ et un tel mineur $\mathsf{F}_{0}$, l'in{\'e}galit{\'e} de Hadamard entra{\^\i}ne \begin{equation*}\vert\det\mathsf{F}_{0}\vert_{\sigma}\le\prod_{i=1}^{\codim F}{\Vert C_{i}\Vert_{\sigma}}\end{equation*}(les colonnes de $\mathsf{F}_{0}$ {\'e}tant en norme plus petites que celles de $\mathsf{F}$). Comme il y a $\binom{\dim E}{\dim F}$ mineurs maximaux possibles pour $\mathsf{F}$, on obtient \begin{equation}\label{majdelahauteurinfinie}h_{\infty}(F)\le\sum_{(s,\mathbf{t})}{h_{\infty,\mathrm{L}^{2}}((a_{\boldsymbol{\lambda},(s,\mathbf{t})})_{\boldsymbol{\lambda}})}+\frac{1}{2}\log\binom{\dim E}{\dim F}\end{equation}(la notation $h_{\infty,\mathrm{L}^{2}}$ signifie que nous n'avons pris que la somme portant sur les places archim{\'e}diennes de $K$ des normes $\mathrm{L}^{2}$ des vecteurs $(a_{\boldsymbol{\lambda},(s,\mathbf{t})})_{\boldsymbol{\lambda}}$). Le lemme~\ref{lemmee:majarchitaylor} fournit un majorant de chacun des $\vert a_{\boldsymbol{\lambda},(s,\mathbf{t})}\vert_{\sigma}$ et donc une estimation de la somme ci-dessus.\par En ce qui concerne l'estimation du d{\'e}terminant $\det\mathsf{F}_{0}$ en une place ultram\'etrique $v$, nous constatons qu'il n'est pas possible d'utiliser en l'\'etat la proposition~\ref{prop:estultrametriqueasept} puisque les coefficients de la matrice $\mathsf{F}_{0}$ ne proviennent pas (\emph{a priori}) de polyn\^omes qui s'annulent aux ordres de d\'erivations pr\'ec\'edents. Nous allons donc faire appara{\^\i}tre de tels polyn\^omes en proc\'edant de la mani\`ere suivante. Consid\'erons une colonne $(s,\mathbf{t})$ de $\mathsf{F}_{0}$ telle que la longueur de $\mathbf{t}$ soit maximale et notons $L_{1},\ldots,L_{\codim F}$ les lignes du mineur $\mathsf{F}_{0}$, qui, elles-m\^emes, correspondent (respectivement) aux lignes $\lambda_{i_{1}},\ldots,\lambda_{i_{\codim F}}$ de $\mathsf{F}$. Si nous supprimons la colonne $(s,\mathbf{t})$ de $\mathsf{F}_{0}$, les lignes $\widetilde{L}_{j}$ de la matrice restante ($j$ est un entier compris entre $1$ et $\codim F$, et $\widetilde{L}_{j}:=(a_{\lambda_{i_{j}},(s,\mathbf{t}')})_{\mathbf{t}'\ne\mathbf{t}}$) sont li\'ees sur $\overline{\mathbf{Q}}$ et il est possible de trouver une relation de d\'ependance lin\'eaire entre ces lignes sur $K\cap\mathcal{O}_{v}$ avec au moins un des coefficients \'egal \`a $1$. Pour cela, il suffit de consid{\'e}rer une relation sur $\mathcal{O}_{K}$ puis de diviser par le coefficient dont la valeur absolue $v$-adique est minimale (non nulle). Autrement dit, il existe $j_{0}\in\{1,\ldots,\codim F\}$ et $(\alpha_{j})_{j\ne j_{0}}\in(K\cap\mathcal{O}_{v})^{\codim F-1}$ tels que $$\widetilde{L}_{j_{0}}=\sum_{\genfrac{}{}{0pt}{}{j=1}{j\ne j_{0}}}^{\codim F}{\alpha_{j}\widetilde{L}_{j}}\ .$$Ainsi, en soustrayant $\sum_{j\ne j_{0}}{\alpha_{j}L_{j}}$ \`a la $j_{0}^{\text{\`eme}}$ ligne de $\mathsf{F}_{0}$, nous obtenons une matrice de m\^eme d\'eterminant que $\mathsf{F}_{0}$ et dont la $j_{0}^{\text{\`eme}}$ ligne est compos\'ee de z\'eros sauf \`a la position $(s,\mathbf{t})$ o\`u le coefficient est $$\xi_{(s,\mathbf{t})}:=a_{\lambda_{i_{j_{0}}},(s,\mathbf{t})}-\sum_{j\ne j_{0}}{\alpha_{j}a_{\lambda_{i_{j}},(s,\mathbf{t})}}\ \cdotp$$Par cons\'equent, la valeur absolue $v$-adique de $\det\mathsf{F}_{0}$ est \'egale \`a $\vert\xi_{(s,\mathbf{t})}\vert_{v}$ multipli\'e par la valeur absolue du d\'eterminant d'un mineur $\widetilde{F}_{0}$ de taille $\codim F-1$ de la matrice $\mathsf{F}$ \`a laquelle on a retir\'e la colonne $(s,\mathbf{t})$ et la ligne $j_{0}$. Soit $(q_{\boldsymbol{\lambda}})_{\boldsymbol{\lambda}}$ l'\'el\'ement de $K^{\dim E}$ d\'efini par $q_{\lambda_{i_{j_{0}}}}=1$, $q_{\lambda_{i_{j}}}=-\alpha_{j}$ si $j\in\{1,\ldots,\codim F\}\setminus\{j_{0}\}$ et $q_{\boldsymbol{\lambda}}=0$ si $L_{\boldsymbol{\lambda}}$ n'est pas une ligne de $\mathsf{F}_{0}$. Le polyn\^ome $Q$ correspondant \`a ces coordonn\'ees v\'erifie $$\cD_{\mathbf{w}}^{\mathbf{t}'}\left(\sum_{\boldsymbol{\lambda}}{q_{\boldsymbol{\lambda}}\Psi_{v_{0},\varepsilon_{s}}^{\boldsymbol{\lambda}}}\right)(s,su)=0$$pour tout $\mathbf{t}'\in\mathbf{N}^{\dim W}$ de longueur strictement inf{\'e}rieure {\`a} $\vert\mathbf{t}\vert$ (puisque $\vert\mathbf{t}\vert$ a {\'e}t{\'e} choisi maximal). En particulier l'{\'e}galit{\'e}~\eqref{eq:eq} et la formule de Leibniz entra{\^\i}nent, pour tout $\vert\mathbf{t}'\vert<\vert\mathbf{t}\vert$, $\cD_{\mathbf{w}}^{\mathbf{t}'}F_{Q,v_{0}}(s,su)=0$ et \begin{equation*}\xi_{(s,\mathbf{t})}=\frac{1}{\prod_{i=1}^{n}{\theta_{v_{0},i,\varepsilon_{i,s}}(su)^{D_{i}}}}\frac{\cD_{\mathbf{w}}^{\mathbf{t}}}{\mathbf{t}!}F_{Q,v_{0}}(s,su)\ \cdotp\end{equation*}La proposition~\ref{prop:estultrametriqueasept} peut donc s'appliquer \`a $\xi_{(s,\mathbf{t})}$, ce qui fournit une majoration de $\vert\xi_{(s,\mathbf{t})}\vert_{v}$ (ind\'ependante de la taille $v$-adique des coefficients de $Q$). En op{\'e}rant de la m\^eme mani\`ere pour $\widetilde{\mathsf{F}}_{0}$, nous d\'eduisons par r{\'e}currence imm{\'e}diate que $\vert\det\mathsf{F}_{0}\vert_{v}$ s'\'ecrit comme un produit $\prod_{i=1}^{\codim F}{\vert\xi_{(s_{i},\mathbf{t}_{i})}\vert_{v}}$ o\`u chacun des $\xi_{(s_{i},\mathbf{t}_{i})}$ est born\'e comme dans la formule~\eqref{majorantpadique} de la proposition~\ref{prop:estultrametriqueasept}. La majoration du logarithme de $\max_{\mathsf{F}_{0}}\{\vert\det\mathsf{F}_{0}\vert_{v}\}$ qui en d{\'e}coule (valable pour toute place ultram\'etrique $v$) et l'in\'egalit\'e archim\'edienne~\eqref{majdelahauteurinfinie} entra{\^\i}nent alors la proposition~\ref{prop:construction}.

\end{proof}
\subsection{Extrapolation}
\label{aspetsection:extrapolation}
\'Etant donn\'e $(s,\boldsymbol{\tau})\in\mathbf{N}\times\mathbf{N}^{\dim W}$ et $P$ le polyn{\^o}me construit au paragraphe pr{\'e}c{\'e}dent, nous disposons d'un coefficient de Taylor tordu \begin{equation}\label{asepteq:taylortordu}\frac{1}{\prod_{i=1}^{n}{\theta_{v_{0},i,\varepsilon_{i}}(su)^{D_{i}}}}\frac{\cD_{\mathbf{w}}^{\boldsymbol{\tau}}}{\boldsymbol{\tau}!}\left(P\circ\exp_{(G_{0}\times G)(\mathbf{C}_{v_{0}})}\right)(s,su)\quad\text{(voir~\eqref{coeffdetaylorultra}).}\end{equation}Pour une raison technique qui appara{\^i}tra \`a la fin de la preuve de la proposition~\ref{asept:propextrapolation} (p.~\pageref{asept:propextrapolation}), nous supposerons que $\varepsilon_{i}\in\{0,\ldots,N_{i}\}$ est choisi de telle sorte que \begin{equation*}\vert\theta_{v_{0},i,\varepsilon_{i}}(su)\vert_{v_{0}}=\max_{0\le j\le N_{i}}{\{\vert\theta_{v_{0},i,j}(su)\vert_{v_{0}}\}}.\end{equation*}En particulier $\theta_{v_{0},i,\varepsilon_{i}}(su)\ne 0$. D'apr\`es le lemme~\ref{lemme:multiplicite} et le choix des param\`etres (voir l'in{\'e}galit{\'e}~\ding{193} du lemme~\ref{lemmearchi:choixdesparametres}), il existe un couple $(s,\boldsymbol{\tau})$ avec $0\le s\le(g+1)S$ et $\vert\boldsymbol{\tau}\vert\le(g+1)T$ pour lequel le nombre~\eqref{asepteq:taylortordu} est non nul. Parmi ces couples, choisissons-en un tel que $(s,\vert\boldsymbol{\tau}\vert)$ soit minimal pour l'ordre lexicographique dans $\mathbf{N}^{2}$, et notons $\alpha$ le terme~\eqref{asepteq:taylortordu} correspondant. On notera que par construction de $P$ on a n\'ecessairement $s\ge S_{0}+1$ et $\alpha\in K$ (lemme~\ref{sdgfkjsgvjq}). Les propositions~\ref{prop:estultrametriqueasept} et~\ref{aseptprop:estiarchimediennes} apportent des estimations de $\vert\alpha\vert_{v}$ en toutes les places $v$ de $K$ et, par suite, de la hauteur de Weil de $\alpha$.\par Consid{\'e}rons une place \emph{quelconque} de $K$ au-dessus de $v_{0}$, place que nous noterons encore $v_{0}$. Nous allons donner ici une majoration de $\vert\alpha\vert_{v_{0}}$ qui d{\'e}pend de la distance $\mathrm{d}_{v_{0}}(u,V)$ de sorte que si celle-ci est \og trop petite\fg, il y aura une contradiction avec la formule du produit. \par Invent\'ee par A.~Baker, la d{\'e}marche consiste {\`a} d{\'e}placer la question sur une droite de $W$. Plus pr{\'e}cis{\'e}ment, si $\widetilde{u}$ est un {\'e}l{\'e}ment de $V\otimes\mathbf{C}_{v_{0}}$ tel que $\mathrm{d}_{v_{0}}(u,V)=\Vert u-\widetilde{u}\Vert_{v_{0}}$, il revient au m{\^e}me --- modulo un terme d'erreur lin{\'e}aire en $\mathrm{d}_{v_{0}}(u,V)$ --- de majorer~\eqref{asepteq:taylortordu} avec $(s,s\widetilde{u})$ au lieu de $(s,su)$ (sous r{\'e}serve, dans le cas $p$-adique, d'avoir v{\'e}rifi{\'e} que $\widetilde{u}$ appartenait bien {\`a} $\mathscr{T}_{v_{0}}$). Le point remarquable est alors le suivant. La fonction analytique d{\'e}finie dans le disque unit{\'e} de $\mathbf{C}_{v_{0}}$ par \begin{equation*}f_{\boldsymbol{\tau}}(z)=\frac{\cD_{\mathbf{w}}^{\boldsymbol{\tau}}}{\boldsymbol{\tau}!}\left(P\circ\exp_{(G_{0}\times G)(\mathbf{C}_{v_{0}})}\right)(z,z\widetilde{u})\end{equation*}admet des d{\'e}riv{\'e}es (divis{\'e}es) qui sont elles-m\^emes des combinaisons lin{\'e}aires tr{\`e}s simples des $f_{\boldsymbol{\tau}'}$, avec $\vert\boldsymbol{\tau}'\vert=\vert\boldsymbol{\tau}\vert+$ordre de d{\'e}rivation, car $(1,\widetilde{u})\in W\otimes\mathbf{C}_{v_{0}}$. Un lemme d'interpolation permet alors de majorer ais{\'e}ment $\vert f_{\boldsymbol{\tau}}(s)\vert_{v_{0}}$ en fonction des valeurs $\vert f_{\boldsymbol{\tau}'}(s_{0})\vert_{v_{0}}$ pour $0\le s_{0}\le S_{0}$ et $\vert\boldsymbol{\tau}'\vert\le 2(g+1)T$. \par Ce sch{\'e}ma de d{\'e}monstration ne d{\'e}pend pas vraiment de la nature, archim{\'e}dienne ou $p$-adique, de la place $v_{0}$. N{\'e}anmoins, il me semble pr{\'e}f{\'e}rable dans un souci de clart{\'e} pour la pr{\'e}sentation de distinguer ces deux cas.
\subsubsection{$v_{0}$ archim{\'e}dienne}
Avant d'{\'e}noncer la proposition principale, commen{\c c}ons par quelques r{\'e}sultats pr{\'e}liminaires usuels dans ce contexte.
\begin{lemm}\label{aseptlemme:schwarz}
Il existe une constante $c_{\cst}\ge 1$ ayant la propri{\'e}t{\'e} suivante. Soit $\mathbf{z}=(z_{0},\ldots,z_{n})\in t_{G_{0}}(\mathbf{C}_{v_{0}})\oplus t_{G}(\mathbf{C}_{v_{0}})$, $\boldsymbol{\lambda}=(\lambda_{0},\lambda_{1},\ldots,\lambda_{n})\in\mathbf{N}\times\prod_{i=1}^{n}{\mathbf{N}^{N_{i}+1}}$ et $\mathbf{h}=(h_{0},\ldots,h_{g-t})\in\mathbf{N}^{\dim W}$. Alors la valeur absolue du coefficient $\frac{1}{\mathbf{h}!}\cD_{\mathbf{w}}^{\mathbf{h}}\Theta^{\boldsymbol{\lambda}}(\mathbf{z})$ est major{\'e}e par \begin{equation*}c_{\theconstante}^{\vert\mathbf{h}\vert+\vert\boldsymbol{\lambda}\vert}\left(\prod_{j=1}^{g-t}{\Vert w_{j}\Vert^{h_{j}}_{v_{0}}}\right)\max_{\lambda_{0}}{\left\{\left\vert\frac{1}{h_{0}!}P_{\lambda_{0}}^{(h_{0})}(z_{0})\right\vert\right\}}\exp{\left\{c_{\theconstante}\sum_{i=1}^{n}{\vert\lambda_{i}\vert(1+\Vert z_{i}\Vert_{v_{0}})^{\rho_{i}}}\right\}}\end{equation*}(voir~\eqref{asepteq:pagequinze} pour la d{\'e}finition de $\Theta^{\boldsymbol{\lambda}}$). 
\end{lemm}
La preuve de ce lemme est essentiellement la m{\^e}me que celle du lemme~\ref{lemmee:majarchitaylor} en tenant compte de l'in{\'e}galit{\'e}~\eqref{asepteq:inegalitetreize}.\par Soit $m\in\mathbf{N}$, $\mathbf{h}\in\mathbf{N}^{\dim W}$ et $F=F_{P,v_{0}}$ la fonction associ\'ee \`a $P$. Soit {\'e}galement (comme dans l'introduction) $\widetilde{u}$ un vecteur de $V\otimes_{v_{0}}\mathbf{C}$ pour lequel $\mathrm{d}_{v_{0}}(u,V)=\Vert u-\widetilde{u}\Vert_{v_{0}}$. Le th{\'e}or{\`e}me des accroissements finis appliqu{\'e} {\`a} la fonction \begin{equation*}x\longmapsto\frac{1}{\mathbf{h}!}\cD_{\mathbf{w}}^{\mathbf{h}}F(m,mu+xm(\widetilde{u}-u))\end{equation*}de la variable r{\'e}elle $x\in[0,1]$ entra{\^\i}ne alors imm{\'e}diatement le 
\begin{lemm}(Comparer avec le lemme~$8$ de~\cite{gaudron1}).\label{aseptlemme:difference}
Il existe une constante $c_{\cst}\ge 1$ ayant la propri{\'e}t{\'e} suivante. Avec les notations ci-dessus, la valeur absolue de la diff{\'e}rence $$\frac{1}{\mathbf{h}!}\cD_{\mathbf{w}}^{\mathbf{h}}F(m,mu)-\frac{1}{\mathbf{h}!}\cD_{\mathbf{w}}^{\mathbf{h}}F(m,m\widetilde{u})$$est major{\'e}e par \begin{equation*}\begin{split}&c_{\theconstante}^{T+\log(D_{0})}\max_{\boldsymbol{\lambda}}{\{\vert p_{\boldsymbol{\lambda}}\vert_{v_{0}}\}}\left(\prod_{j=1}^{g-t}{\Vert w_{j}\Vert_{v_{0}}^{h_{j}}}\right)m\,\mathrm{d}_{v_{0}}(u,V)\\ &\times\max_{\lambda_{0}}{\left\{\left\vert\frac{1}{h_{0}!}P_{\lambda_{0}}^{(h_{0})}(m)\right\vert\right\}}\times\exp{\left\{c_{\theconstante}\sum_{i=1}^{n}{D_{i}(1+m\Vert u_{i}\Vert_{v_{0}})^{\rho_{i}}}\right\}}\end{split}\end{equation*}pourvu que $\mathrm{d}_{v_{0}}(u,V)\le 1$.
\end{lemm} Nous aurons aussi besoin du lemme d'interpolation suivant, d{\^u} {\`a} Waldschmidt~\cite{miw1980}. Si $x$ est un nombre r{\'e}el positif et $f$ une fonction d{\'e}finie sur le disque ferm{\'e} $\overline{D}(0,x)=\{z\in\mathbf{C}\,;\ \vert z\vert \le x\}$, on note $\vert f\vert_{x}$ la borne sup{\'e}rieure des $\vert f(z)\vert$, $z\in\overline{D}(0,x)$.
\begin{lemm}\label{lemmedemiwasept}
Soit $S_{1},T_{1}$ des entiers naturels strictement positifs et $\mathsf{R}\ge\mathsf{r}\ge 2S_{1}$ des nombres r{\'e}els. Soit $f$ une fonction analytique (d'une variable complexe) dans le disque $\overline{D}(0,\mathsf{R})$. Alors on a \begin{equation*}\vert f\vert_{\mathsf{r}}\le 2\vert f\vert_{\mathsf{R}}\left(\frac{2\mathsf{r}}{\mathsf{R}}\right)^{T_{1}S_{1}}+5\left(\frac{9\mathsf{r}}{S_{1}}\right)^{T_{1}S_{1}}\times\max_{\genfrac{}{}{0pt}{}{0\le h<T_{1}}{0\le m<S_{1}}}{\left\{\left\vert\frac{1}{h!}f^{(h)}(m)\right\vert\right\}}\ \cdotp\end{equation*} 
\end{lemm}
De la sorte, pour $m\le S$ et au moyen des in{\'e}galit{\'e}s du lemme~\ref{lemmearchi:choixdesparametres}, le majorant du lemme~\ref{aseptlemme:difference} s'{\'e}crit plus simplement $$\max_{\boldsymbol{\lambda}}{\{\vert p_{\boldsymbol{\lambda}}\vert_{v_{0}}\}}\max_{j}{\{1,\Vert w_{j}\Vert_{v_{0}}\}}^{c_{\cst}T}e^{c_{\theconstante}U}\mathrm{d}_{v_{0}}(u,V)$$pour une certaine constante $c_{\theconstante}\ge 1$ et o\`u $j$ parcourt l'ensemble $\{1,\ldots,g-t\}$.\par Le r{\'e}sultat central de ce paragraphe est le suivant.
\begin{prop}
\label{asept:propextrapolation}
Supposons que $\log\mathrm{d}_{v_{0}}(u,V)\le-C_{0}^{3}U$. Alors $$\vert\alpha\vert_{v_{0}}\le e^{-U}\max_{\boldsymbol{\lambda}}{\{\vert p_{\boldsymbol{\lambda}}\vert_{v_{0}}\}}\max{\{1,\Vert w_{1}\Vert_{v_{0}},\ldots,\Vert w_{g-t}\Vert_{v_{0}}\}}^{C_{0}T}\ .$$
\end{prop}
\begin{proof}Soit $\mathrm{f}:\mathbf{C}\to\mathbf{C}$ la fonction enti{\`e}re d{\'e}finie par $$\mathrm{f}(z)=\frac{1}{\boldsymbol{\tau}!}\cD_{\mathbf{w}}^{\boldsymbol{\tau}}F(z,z\widetilde{u})\ .$$Soit $\mathbf{x}:=(x_{0}=1,x_{1},\ldots,x_{g-t})$ les coordonn{\'e}es de $(1,\widetilde{u})$ dans la base $\mathbf{w}$. La hauteur $\mathrm{L}^{2}$ de $w_{i}$ est projective et quitte {\`a} multiplier $w_{i}$ par un entier assez grand\footnote{Par exemple, $2\max{\{1,\left[\vert x_{1}\vert\right],\ldots,\left[\vert x_{g-t}\vert\right]\}}$ convient. Cette astuce ne d{\'e}pend que de la place $v_{0}$ sur $k$ (et non du choix de la place de $K$ au-dessus de $v_{0}$).} (pour $1\le i\le g-t$), nous pouvons supposer que $\vert x_{i}\vert\le 1$. Pour tout entier $\ell\ge 0$, la d{\'e}riv{\'e}e $\ell^{\text{\`eme}}$ de $\mathrm{f}$ v{\'e}rifie la formule\begin{equation*}\frac{\mathrm{f}^{(\ell)}(z)}{\ell!}=\sum_{\genfrac{}{}{0pt}{}{\mathbf{j}\in\mathbf{N}^{\dim W}}{\vert\mathbf{j}\vert=\ell}}{\binom{\boldsymbol{\tau}+\mathbf{j}}{\mathbf{j}}\,\mathbf{x}^{\mathbf{j}}\,\frac{\cD_{\mathbf{w}}^{\boldsymbol{\tau}+\mathbf{j}}}{(\boldsymbol{\tau}+\mathbf{j})!}F(z,z\widetilde{u})}\end{equation*}et donc il existe une constante $c_{\cst}\ge 1$\newcounter{cdixhuit}\setcounter{cdixhuit}{\value{constante}} telle que \begin{equation}\label{ineq:ilfautunnumeroasept}\max_{\genfrac{}{}{0pt}{}{\ell\le(g+1)T}{s_{0}\le S_{0}}}\left\{\left\vert\frac{\mathrm{f}^{(\ell)}(s_{0})}{\ell!}\right\vert\right\}\le c_{\theconstante}^{T}\max_{(s_{0},\mathbf{j})\in\Upsilon}{\left\vert\frac{\cD_{\mathbf{w}}^{\mathbf{j}}}{\mathbf{j}!}F(s_{0},s_{0}\widetilde{u})\right\vert}\ \cdotp\end{equation}Si l'on remplace $\widetilde{u}$ par $u$ dans le membre de droite de cette in{\'e}galit{\'e}, le terme obtenu est nul par construction de $F$. Par cons{\'e}quent, le lemme de comparaison~\ref{aseptlemme:difference} entra{\^\i}ne la majoration \begin{equation*}\max_{\genfrac{}{}{0pt}{}{\ell\le(g+1)T}{s_{0}\le S_{0}}}\left\{\left\vert\frac{\mathrm{f}^{(\ell)}(s_{0})}{\ell!}\right\vert\right\}\le e^{-C_{0}^{2}U}\max_{\boldsymbol{\lambda}}{\{\vert p_{\boldsymbol{\lambda}}\vert_{v_{0}}\}}\max_{j}{\{1,\Vert w_{j}\Vert_{v_{0}}\}}^{c_{\cst}T}\ .\end{equation*}De m{\^e}me, le lemme~\ref{aseptlemme:schwarz} conduit {\`a} la majoration\begin{equation}\begin{split}\label{asepteq:majdusup}\vert \mathrm{f}\vert_{\mathsf{R}}\le &\max_{\boldsymbol{\lambda}}{\{\vert p_{\boldsymbol{\lambda}}\vert_{v_{0}}\}}e^{\left(U/2+\sum_{i=1}^{n}{C_{0}D_{i}(1+\mathsf{R}\Vert u_{i}\Vert_{v_{0}})^{\rho_{i}}}\right)}\times\max_{\genfrac{}{}{0pt}{}{\lambda_{0}}{\vert z_{0}\vert\le\mathsf{R}}}{\left\{\left\vert\frac{1}{\tau_{0}!}P_{\lambda_{0}}^{(\tau_{0})}(z_{0})\right\vert\right\}}\\ &\times\max_{j}{\{\Vert w_{j}\Vert_{v_{0}}\}}^{c_{\cst}T}\end{split}\end{equation}valide pour tout nombre r{\'e}el $\mathsf{R}\in[1,\mathrm{d}_{v_{0}}(u,V)^{-1}]$. En prenant $\mathsf{R}:=2\mathfrak{e}(g+1)S$, nous constatons que le produit des deux derni{\`e}res quantit{\'e}s de la premi{\`e}re ligne du membre de droite de~\eqref{asepteq:majdusup} est inf{\'e}rieure {\`a} $e^{U}$, ce qui fournit la majoration plus simple \begin{equation*}\vert\mathrm{f}\vert_{\mathsf{R}}\le\max_{\boldsymbol{\lambda}}{\{\vert p_{\boldsymbol{\lambda}}\vert_{v_{0}}\}}\max_{j}{\{1,\Vert w_{j}\Vert_{v_{0}}\}}^{c_{\theconstante}T}\,e^{U}\ .\end{equation*}Mentionnons que pour cela nous avons utilis{\'e} les in{\'e}galit{\'e}s~\ding{195} et~\ding{197} du lemme~\ref{lemmearchi:choixdesparametres}. Choisissons $T_{1}:=(g+1)T$, $S_{1}:=S_{0}$, $\mathsf{r}:=s$ et appliquons le lemme d'interpolation~\ref{lemmedemiwasept} {\`a} ces donn{\'e}es et {\`a} la fonction $\mathrm{f}$. Il vient $$\vert \mathrm{f}(s)\vert\le e^{-C_{0}U}\max_{\boldsymbol{\lambda}}{\{\vert p_{\boldsymbol{\lambda}}\vert\}}\max_{j}{\{1,\Vert w_{j}\Vert_{v_{0}}\}}^{c_{\cst}T}$$(rappelons que $\mathsf{r}/S_{1}=s/S_{0}\le S/S_{0}\le C_{0}^{4}$), puis, par une seconde application du lemme de comparaison~\ref{aseptlemme:difference}, nous obtenons \begin{equation*}\left\vert\frac{1}{\boldsymbol{\tau}!}\cD_{\mathbf{w}}^{\boldsymbol{\tau}}F(s,su)\right\vert\le 2\,e^{-C_{0}U}\,\max_{\boldsymbol{\lambda}}{\{\vert p_{\boldsymbol{\lambda}}\vert\}}\max_{j}{\{1,\Vert w_{j}\Vert_{v_{0}}\}}^{c_{\cst}T}\ \cdotp\end{equation*}La borne pour $\vert\alpha\vert_{v_{0}}$ se d{\'e}duit de cette in{\'e}galit{\'e} \emph{via} la minoration de $\vert\theta_{v_{0},i,\varepsilon_{i}}(su)\vert$ donn{\'e}e par l'in{\'e}galit{\'e}~\eqref{asepteq:inegalitetreize} (p.~\pageref{asepteq:inegalitetreize}) et gr{\^a}ce au choix de $\varepsilon_{i}$.  
\end{proof}
\subsubsection{$v_{0}$ ultram{\'e}trique}\label{subssubsec:voultra}
Consid{\'e}rons le nombre alg{\'e}brique $\alpha$ introduit au d{\'e}but du paragraphe~\ref{aspetsection:extrapolation}. La condition de minimalit{\'e} sur $(s,\vert\boldsymbol{\tau}\vert)$ implique que $\alpha$ est aussi {\'e}gal {\`a} \begin{equation}\label{autreformedealpha}\frac{\cD_{\mathbf{w}}^{\boldsymbol{\tau}}}{\boldsymbol{\tau}!}\left(\sum_{\boldsymbol{\lambda}}{p_{\boldsymbol{\lambda}}\Psi_{v_{0},0}^{\boldsymbol{\lambda}}}\right)(s,su)\end{equation}o\`u, comme dans la d{\'e}monstration de la proposition~\ref{prop:construction}, $\Psi_{v_{0},0}^{\boldsymbol{\lambda}}$ d{\'e}signe $$P_{\lambda_{0}}(z_{0})\prod_{i=1}^{n}{\Psi_{v_{0},i,0}^{\lambda_{i}}}\ .$$Dans la suite, nous noterons $\widetilde{F}$ la somme des $p_{\boldsymbol{\lambda}}\Psi_{v_{0},0}^{\boldsymbol{\lambda}}$. L'{\'e}criture~\eqref{autreformedealpha} pour $\alpha$ s'av{\`e}re plus commode dans le cas ultram{\'e}trique car l'on conna{\^\i}t un d{\'e}veloppement en s{\'e}rie de $\Psi_{v_{0},0}^{\boldsymbol{\lambda}}$, qui converge sur le disque ouvert $D(0,r_{p})$~:\begin{equation*}\forall\boldsymbol{\lambda},\ \exists\,(a_{\mathbf{i},\boldsymbol{\lambda}})_{\mathbf{i}}\in\mathcal{O}_{v_{0}}^{\mathbf{N}^{g}}\,;\quad\Psi_{v_{0},0}^{\boldsymbol{\lambda}}(z)=P_{\lambda_{0}}(z_{0})\left(\sum_{\mathbf{i}\in\mathbf{N}^{g}}{\frac{a_{\mathbf{i},\boldsymbol{\lambda}}}{\mathbf{i}!}\mathbf{z}^{\mathbf{i}}}\right)\end{equation*}lorsque $\mathbf{z}=(z_{1},\ldots,z_{g})$ d{\'e}signe les coordonn{\'e}es de $z\in D(0,r_{p})$ dans la base $\mathbf{e}$ (rappelons que c'est pr{\'e}cis{\'e}ment le choix de cette base qui assure l'int{\'e}gralit{\'e} des coefficients $a_{\mathbf{i},\boldsymbol{\lambda}}$, voir \S~\ref{donnesgeneralesasept}). Nous avons conserv{\'e} $P_{\lambda_{0}}(z_{0})$ intact dans ce d{\'e}veloppement pour des raisons pratiques afin de ne pas m{\'e}langer les coefficients de $P_{\lambda_{0}}$ (qui n'appartiennent pas n{\'e}cessairement {\`a} $\mathcal{O}_{v_{0}}$) avec les $a_{\mathbf{i},\boldsymbol{\lambda}}\in\mathcal{O}_{v_{0}}$. Pour tout multiplet $\mathbf{h}=(h_{0},\mathbf{h}')\in\mathbf{N}\times\mathbf{N}^{g}$ et tout {\'e}l{\'e}ment $(z_{0},z)\in t_{G_{0}}(\mathbf{C}_{v_{0}})\times D(0,r_{p})$, nous obtenons alors \begin{equation}\label{eq:asept:ftilde}\frac{\cD_{\mathbf{e}}^{\mathbf{h}}}{\mathbf{h}!}\widetilde{F}(z_{0},z)=\sum_{\boldsymbol{\lambda},\mathbf{i}}{p_{\boldsymbol{\lambda}}\frac{P_{\lambda_{0}}^{(h_{0})}(z_{0})}{h_{0}!}\frac{a_{\mathbf{i},\boldsymbol{\lambda}}}{(\mathbf{i}-\mathbf{h}')!}\mathbf{z}^{\mathbf{i}-\mathbf{h}'}}\end{equation}($\mathbf{i}$ parcourt $\{(i_{1},\ldots,i_{g})\in\mathbf{N}^{g}\,;\ \forall\,j\in\{1,\ldots,g\},\ i_{j}\ge h'_{j}\}$). Pour $j$ un entier naturel, notons $\sigma_{p}(j)$ la somme des chiffres de $j$ {\'e}crit en base $p$ (le nombre premier $p$ est la caract\'eristique r\'esiduelle de $v_{0}$). On sait que la valuation $p$-adique de $j!$ est $(j-\sigma_{p}(j))/(p-1)$ ce qui entra{\^\i}ne $\vert j!\vert_{v_{0}}\ge r_{p}^{j}$ et, plus g{\'e}n{\'e}ralement, pour $\mathbf{i}\in\mathbf{N}^{g}$, $\vert\mathbf{i}!\vert_{v_{0}}\ge r_{p}^{\vert\mathbf{i}\vert}$ si bien que $\vert\mathbf{z}^{\mathbf{i}}/\mathbf{i}!\vert_{v_{0}}\le 1$ pour $z\in D(0,r_{p})$. De l'{\'e}galit{\'e}~\eqref{eq:asept:ftilde} et de l'in{\'e}galit{\'e} ultram{\'e}trique se d{\'e}duit la majoration \begin{equation}\label{ineqasept:rond}\left\vert\frac{\cD_{\mathbf{e}}^{\mathbf{h}}}{\mathbf{h}!}\widetilde{F}(z_{0},z)\right\vert_{v_{0}}\le\max_{\boldsymbol{\lambda}}{\{\vert p_{\boldsymbol{\lambda}}\vert_{v_{0}}\}}\max_{\lambda_{0}}{\left\{\left\vert\frac{P_{\lambda_{0}}^{(h_{0})}(z_{0})}{h_{0}!}\right\vert_{v_{0}}\right\}}\end{equation}valide pour tout $(z_{0},z)\in t_{G_{0}}(\mathbf{C}_{v_{0}})\times D(0,r_{p})$. Cette in{\'e}galit{\'e} se g{\'e}n{\'e}ralise imm{\'e}diatement {\`a} une base quelconque $\mathbf{e}'=(e_{0}',\ldots,e_{g}')$ de $t_{G_{0}\times G}(\mathbf{C}_{v_{0}})$ en multipliant le membre de droite par $\prod_{j=0}^{g}{\Vert e_{j}'\Vert_{v_{0}}^{h_{j}}}$. De la m{\^e}me fa{\c c}on, l'estimation $$\vert \mathbf{z}^{\mathbf{i}}-\mathbf{z}'^{\mathbf{i}}\vert\le\Vert z-z'\Vert_{v_{0}}\max{\{\Vert z\Vert_{v_{0}},\Vert z'\Vert_{v_{0}}\}}^{\vert\mathbf{i}\vert-1}$$conduit au lemme de comparaison suivant.
\begin{lemm}\label{lemmeultraaseptdifference}
  Supposons que $\mathrm{d}_{v_{0}}(u,V)=\Vert u-\widetilde{u}\Vert_{v_{0}}$ est strictement inf{\'e}rieur {\`a} $r_{p}$. Alors, pour tout couple $(m,\mathbf{h})\in\mathbf{N}\times\mathbf{N}^{\dim W}$, la valeur absolue de la diff{\'e}rence $\frac{\cD_{\mathbf{w}}^{\mathbf{h}}}{\mathbf{h}!}\widetilde{F}(m,mu)-\frac{\cD_{\mathbf{w}}^{\mathbf{h}}}{\mathbf{h}!}\widetilde{F}(m,m\widetilde{u})$ est major{\'e}e par \begin{equation*}\max_{\boldsymbol{\lambda}}{\{\vert p_{\boldsymbol{\lambda}}\vert_{v_{0}}\}}\max_{\lambda_{0}}{\left\{\left\vert\frac{P_{\lambda_{0}}^{(h_{0})}(m)}{h_{0}!}\right\vert_{v_{0}}\right\}}\prod_{j=1}^{g-t}{\Vert w_{j}\Vert_{v_{0}}^{h_{j}}}\times\mathrm{d}_{v_{0}}(u,V)\ .\end{equation*}
\end{lemm}
La condition $\Vert u-\widetilde{u}\Vert_{v_{0}}<r_{p}$ {\'e}quivaut {\`a} $\Vert\widetilde{u}\Vert_{v_{0}}<r_{p}$ (puisque $\Vert u\Vert_{v_{0}}<r_{p}$) et assure de la sorte la coh{\'e}rence de l'{\'e}nonc{\'e}.\par Enfin, comme dans le cas archim{\'e}dien, nous aurons besoin d'un lemme d'interpolation (en une variable), d{\^u} {\`a} Roy~\cite{Roy2001}.
\begin{lemm}\label{lemmedinterpolationpadique}Soit $S_{1},T_{1}$ des entiers $\ge 1$ et $\mathsf{R}\ge\mathsf{r}\ge 1$ des nombres r{\'e}els. Posons $$\kappa:=\frac{S_{1}-\sigma_{p}(S_{1})}{p-1}+\left[\frac{\log S_{1}}{\log p}\right]$$Soit $f:\overline{D}(0,\mathsf{R})\to\mathbf{C}_{p}$ une fonction analytique. Alors \begin{equation*}\frac{\vert f\vert_{\mathsf{r}}}{\mathsf{r}^{(S_{1}+1)T_{1}}}\le p^{\kappa T_{1}}\,\max{\left\{\max_{\genfrac{}{}{0pt}{}{0\le m< S_{1}}{0\le h<T_{1}}}{\left\{\left\vert\frac{f^{(h)}(m)}{h!}\right\vert\right\}},\left(\frac{1}{\mathsf{R}}\right)^{(S_{1}+1)T_{1}}\vert f\vert_{\mathsf{R}}\right\}}\ \cdotp\end{equation*}
\end{lemm}
\begin{proof}
Il s'agit d'un cas tr{\`e}s particulier du corollaire~$1.2$ de~\cite{Roy2001} qui avait {\'e}t{\'e} conjectur{\'e} par P.~Robba en 1978. Avec les notations de cet article, choisissons $K=\mathbf{C}_{p}$, $n=\rho=1$, $E=\Omega=\{0,\ldots,S_{1}\}$, $L=M=(S_{1}+1)T_{1}$. L'{\'e}nonc{\'e} de Roy donne le r{\'e}sultat avec $\Delta(E)\delta(E)$ {\`a} la place de $p^{-\kappa}$, o\`u $\Delta(E)=\min_{x\in E}{\prod_{y\in E\setminus\{x\}}{\vert y-x\vert_{v_{0}}}}$ et $\delta(E)=\min_{x\ne y\in E}{\vert y-x\vert_{v_{0}}}$. En observant que \begin{equation*}\Delta(E)=\vert S_{1}!\vert_{v_{0}}=\vert p\vert_{v_{0}}^{\frac{S_{1}-\sigma_{p}(S_{1})}{p-1}}\quad\text{et}\quad\delta(E)=\vert p\vert_{v_{0}}^{\left[\frac{\log S_{1}}{\log p}\right]},\end{equation*} on a $\Delta(E)\delta(E)=\vert p\vert_{v_{0}}^{\kappa}=p^{-\kappa}$.
\end{proof}
Ces pr{\'e}liminaires {\`a} l'extrapolation {\'e}tant acquis, nous allons {\^e}tre en mesure de d{\'e}montrer la 
\begin{prop}
Supposons que $\log\mathrm{d}_{v_{0}}(u,V)\le-C_{0}^{3}(1+(\log r_{p}^{-1})/\log\mathfrak{r})U$. Alors $$\vert\alpha\vert_{v_{0}}\le e^{-U}\max_{\boldsymbol{\lambda}}{\{\vert p_{\boldsymbol{\lambda}}\vert_{v_{0}}\}}\max_{j}{\{1,\Vert w_{j}\Vert_{v_{0}}\}}^{C_{0}T}\ .$$
\end{prop}
\begin{proof}
Elle suit d'assez pr{\`e}s celle de la proposition~\ref{asept:propextrapolation} en comportant n{\'e}anmoins quelques variantes li{\'e}es, en particulier, {\`a} la finitude du rayon de convergence de l'exponentielle $p$-adique. Soit $D_{u}=\{z\in\mathbf{C}_{p},zu\in\mathscr{T}_{v_{0}}\}$ et $\mathrm{f}:D_{u}\to\mathbf{C}_{p}$ la fonction d{\'e}finie par $\mathrm{f}(z)=\frac{1}{\boldsymbol{\tau}!}\cD_{\mathbf{w}}^{\boldsymbol{\tau}}\widetilde{F}(z,z\widetilde{u})$. L'hypoth{\`e}se $\log\mathrm{d}_{v_{0}}(u,V)\le-C_{0}^{3}U$ entra{\^\i}ne la majoration $\Vert u-\widetilde{u}\Vert_{v_{0}}<\Vert u\Vert_{v_{0}}$ et donc $\Vert u\Vert_{v_{0}}=\Vert\widetilde{u}\Vert_{v_{0}}$, ce qui assure l'analycit{\'e} de $\mathrm{f}$ sur $D_{u}$. Les d{\'e}riv{\'e}es de cette application v{\'e}rifient la m{\^e}me formule que dans le cas archim{\'e}dien et la majoration~\eqref{ineq:ilfautunnumeroasept} est vraie avec $c_{\thecdixhuit}=1$. Le lemme~\ref{lemmeultraaseptdifference} entra{\^\i}ne alors\begin{equation}\begin{split}\label{inequationalphaasept}\max_{\genfrac{}{}{0pt}{}{0\le s_{0}\le S_{0}}{0\le \ell\le(g+1)T}}{\left\{\left\vert\frac{\mathrm{f}^{(\ell)}(s_{0})}{\ell!}\right\vert_{v_{0}}\right\}}\le &\max{\{\vert p_{\boldsymbol{\lambda}}\vert_{v_{0}}\}}\max_{\genfrac{}{}{0pt}{}{\genfrac{}{}{0pt}{}{\lambda_{0}}{h_{0}\le(g+1)T}}{s_{0}\le S_{0}}}{\left\{\left\vert \frac{P_{\lambda_{0}}^{(h_{0})}(s_{0})}{h_{0}!}\right\vert_{v_{0}}\right\}}\,\mathrm{d}_{v_{0}}(u,V)\\ &\times\max{\{1,\Vert w_{1}\Vert_{v_{0}},\ldots,\Vert w_{g-t}\Vert_{v_{0}}\}}^{2(g+1)T}\\ \le &\max{\{\vert p_{\boldsymbol{\lambda}}\vert_{v_{0}}\}}\max_{j}{\{1,\Vert w_{j}\Vert_{v_{0}}\}}^{2(g+1)T}\\ &\times\exp{\left\{-C_{0}^{2}(1+(\log r_{p}^{-1})/\log\mathfrak{r})U\right\}}\ \cdotp\end{split}\end{equation}En outre, de l'in{\'e}galit{\'e}~\eqref{ineqasept:rond} et de la remarque qui suit, nous d{\'e}duisons la majoration \begin{equation*}\vert\mathrm{f}\vert_{\mathsf{R}}\le\max{\{\vert p_{\boldsymbol{\lambda}}\vert_{v_{0}}\}}\max_{\genfrac{}{}{0pt}{}{\lambda_{0}}{\vert z_{0}\vert\le\mathsf{R}}}{\left\{\left\vert\frac{P_{\lambda_{0}}^{(\tau_{0})}(z_{0})}{\tau_{0}!}\right\vert_{v_{0}}\right\}}\prod_{j=1}^{g-t}{\Vert w_{j}\Vert_{v_{0}}^{\tau_{j}}}\end{equation*}valide pour tout $\mathsf{R}$ dans l'intervalle $\left[0,\frac{r_{p}}{\Vert u\Vert_{v_{0}}}\right[$. En appliquant le lemme d'interpolation~\ref{lemmedinterpolationpadique} avec $S_{1}:=S_{0}$, $T_{1}:=(g+1)T$, $\mathsf{r}:=1$ et $\mathsf{R}:=\mathfrak{r}\,r_{p}^{-1}$ {\`a} la fonction $f:=\mathrm{f}$. Le nombre $\kappa=(S_{0}-\sigma_{p}(S_{0}))/(p-1)+[\log S_{0}/\log p]$ de ce lemme est naturellement major\'e par $S_{0}/(p-1)+(\log S_{0})/\log p$ et le choix des param{\`e}tres entra{\^\i}ne \begin{equation*}p^{\kappa(g+1)T}\le (r_{p}^{-1})^{(g+1)TS_{0}}e^{U}\le\exp{\{C_{0}^{3/2}(1+(\log r_{p}^{-1})/\log\mathfrak{r})U\}}\ \cdotp\end{equation*}De la sorte, et gr{\^a}ce {\`a} l'estimation~\eqref{inequationalphaasept}, nous obtenons
\begin{equation*}\begin{split}&p^{\kappa(g+1)T}\max_{\genfrac{}{}{0pt}{}{0\le s_{0}\le S_{0}}{0\le \ell\le(g+1)T}}{\left\{\left\vert\frac{\mathrm{f}^{(\ell)}(s_{0})}{\ell!}\right\vert_{v_{0}}\right\}}\\ &\qquad\le\max{\{\vert p_{\boldsymbol{\lambda}}\vert_{v_{0}}\}}\max_{j}{\{1,\Vert w_{j}\Vert_{v_{0}}\}}^{2(g+1)T}e^{-C_{0}U}\ \cdotp\end{split}\end{equation*}Et par ailleurs nous avons \begin{equation*}\begin{split}&p^{\kappa(g+1)T}\left(\frac{1}{\mathsf{R}}\right)^{(S_{0}+1)(g+1)T}\vert f\vert_{\mathsf{R}}\\ &\qquad\le\max{\{\vert p_{\boldsymbol{\lambda}}\vert_{v_{0}}\}}\max_{j}{\{1,\Vert w_{j}\Vert_{v_{0}}\}}^{2(g+1)T}e^{-C_{0}^{1/2}U}\ \cdotp\end{split}\end{equation*}Le lemme~\ref{lemmedinterpolationpadique} implique alors\begin{equation*}\vert\mathrm{f}(s)\vert\le\vert\mathrm{f}\vert_{1}\le\max_{\boldsymbol{\lambda}}{\{\vert p_{\boldsymbol{\lambda}}\vert_{v_{0}}\}}\max_{j}{\{\Vert w_{j}\Vert_{v_{0}}\}}^{c_{\cst}T}e^{-C_{0}^{1/2}U}\end{equation*}et nous concluons avec le lemme~\ref{lemmeultraaseptdifference}.
\end{proof}
\subsection{Fin de la d{\'e}monstration}\label{subsec:conclusionasept}Comme nous l'avons vu, le nombre alg{\'e}brique $\alpha$ introduit au paragraphe pr{\'e}c{\'e}dent est non nul et satisfait donc {\`a} la formule (du produit)~:\begin{equation}\label{eq:forumuleduproduit}\sum_{\genfrac{}{}{0pt}{}{\text{$v$ place de $K$}}{v\nmid v_{0}}}{\frac{[K_{v}:\mathbf{Q}_{v}]}{[K:\mathbf{Q}]}\log\vert\alpha\vert_{v}}=-\sum_{v\mid v_{0}}{\frac{[K_{v}:\mathbf{Q}_{v}]}{[K:\mathbf{Q}]}\log\vert\alpha\vert_{v}}\ \cdotp\end{equation}D'apr{\`e}s les propositions~\ref{prop:estultrametriqueasept} et~\ref{aseptprop:estiarchimediennes} (appliqu{\'e}es avec $(z_{0},z)=(s,su)$), il existe une constante $c_{\cst}$ telle que le membre de gauche soit major{\'e} par \begin{equation}\label{eqaspet:tresintermediaire}c_{\theconstante}\left(T(\chi_{H}+\sum_{i=1}^{g-t}{h_{\mathrm{L}^{2}}(w_{i})})+h(P)+\sum_{i=1}^{n}{D_{i}h(sp_{i})}+\aleph((P_{\lambda_{0}}))\right)\cdotp\end{equation}Au passage, il faut noter qu'appara{\^\i}t dans le majorant du membre de gauche de~\eqref{eq:forumuleduproduit} la quantit{\'e} \begin{equation}\label{eq:intermed}\frac{1}{[K:\mathbf{Q}]}\sum_{\genfrac{}{}{0pt}{}{v\nmid v_{0}}{i\in\{1,\ldots,m\}}}{\log\frac{1}{\left\vert A_{\varepsilon_{i}}^{(i)}(\psi_{v_{0},i,\varepsilon_{i}}(su),(1:0:\cdots:0))\right\vert_{v}^{D_{i}}}}\end{equation}qui, en vertu de la formule du produit, vaut $$\frac{1}{[k:\mathbf{Q}]}\sum_{i=1}^{m}{\log\left\vert A_{\varepsilon_{i}}^{(i)}(\psi_{v_{0},i,\varepsilon_{i}}(su),(1:0:\cdots:0))\right\vert_{v_{0}}^{D_{i}}}\ ,$$expression qui elle-m{\^e}me est plus petite que $c_{\cst}\sum_{i=1}^{m}{D_{i}(h(sp_{i})+1)}$ pour une certaine constante $c_{\theconstante}\ge 1$. La somme~\eqref{eq:intermed} ci-dessus est donc bien comprise dans le majorant~\eqref{eqaspet:tresintermediaire}. Par construction, nous disposons d'une majoration de la somme des $h_{\mathrm{L}^{2}}(w_{j})$ (voir \S~\ref{subsec:constructionbaseW}) ainsi que de la hauteur de $P$, donn{\'e}e par la proposition~\ref{prop:construction}. En vertu du lemme~\ref{lemme:Raynaud}, la quantit{\'e} $\chi_{H}$ est born{\'e}e. Les lemmes~\ref{lemmearchi:choixdesparametres} et~\ref{lemmeultra:choixdesparametres} montrent alors que la quantit{\'e}~\eqref{eqaspet:tresintermediaire} est major{\'e}e par $U/(C_{0}D)$. Quant au membre de droite de l'{\'e}galit{\'e}~\eqref{eq:forumuleduproduit}, il est minor{\'e} par $U/D-h(P)-C_{0}\sum_{j}{h_{\mathrm{L}^{2}}(w_{j})}\ge U/(2D)$ d{\`e}s lors que $\log\mathrm{d}_{v_{0}}(u,V)\le -C_{0}^{3}U$ (\emph{resp}. $\log\mathrm{d}_{v_{0}}(u,V)\le -C_{0}^{3}(1+(\log r_{p}^{-1})/\log\mathfrak{r})U$) si $v_{0}$ est archim{\'e}dienne (\emph{resp}. ultram{\'e}trique). Nous constatons alors une contradiction avec le majorant $U/(C_{0}D)$. Ce qui conclut la d{\'e}monstration des th{\'e}or{\`e}mes~\ref{theounasept} et~\ref{theodeuxasept} (en observant pour ce dernier que le terme $1+(\log r_{p}^{-1})/\log\mathfrak{r}$ est major{\'e} par $c_{\cst}\log(\mathfrak{r}+1)/\log\mathfrak{r}$ pour une certaine constante absolue $c_{\theconstante}$).

\subsection*{Remerciements} Je remercie G.~Diaz et G.~R{\'e}mond pour leurs remarques et commentaires sur une premi{\`e}re version de ce texte, qui m'ont permis d'am{\'e}liorer la pr{\'e}sentation g{\'e}n{\'e}rale et de corriger de nombreux d{\'e}tails. Je remercie {\'e}galement J.-B.~Bost \`a double titre, \`a la fois pour ses {\'e}claircissements {\`a} propos du lemme~\ref{lemme:Raynaud} et du th{\'e}or{\`e}me de Raynaud sous-jacent, et aussi pour sa lecture critique qui a permis d'all\'eger ce texte.

\bibliographystyle{plain}

\noindent Universit\'e Grenoble I, Institut Fourier.\\ UMR $5582$, BP $74$\\ $38402$ Saint-Martin-d'H{\`e}res Cedex, France.\\ 
Courriel~: \texttt{Eric.Gaudron@ujf-grenoble.fr}\\ 
Page internet~: \texttt{http://www-fourier.ujf-grenoble.fr/\~{}gaudron}

\end{document}